# Differentiability and Covering Constant for Duality Mappings in $l_p$

Jinlu Li

Department of Mathematics
Shawnee State University
Portsmouth, Ohio 45662 USA
jli@shawnee.edu

**Abstract**

In this paper, we investigate the Gâteaux and Fréchet differentiability of some duality mappings in the uniformly convex and uniformly smooth Banach space $l_p$, which includes the normalized duality mapping as a special case and it is denoted by $J$. We will introduce a general duality mapping denoted by $\mathcal{J}$ and a generalized duality mapping denoted by $\mathfrak{J}$ in $l_p$. After the differentiability is proved, we find the explicit Gâteaux, Fréchet and Mordukhovich derivative of $J$, for $2 < p < \infty$, and the duality mappings $\mathcal{J}$ and $\mathfrak{J}$ in $l_p$, for $3 < p < \infty$. By using Mordukhovich derivatives, we find the covering constants for $J$, $\mathcal{J}$ and $\mathfrak{J}$, respectively, in $l_p$. We also prove that, for $2 < p < \infty$, the Gâteaux derivative operator of the normalized duality mapping $J$ has Lipchitz property in $l_p$.

## 1. Introduction

Recall that, in calculus, the development process of the traditional differentiability of functions has been integrated in the following stages: definition of derivatives, properties, explicit formulas of derivatives, and applications. It is well-known that this process has been extended to the differentiation theory in modern analysis.

In analysis of Banach spaces, in contrast with the traditional differentiation in calculus, differentiability of single-valued mappings between Banach spaces have been defined and the most popular concepts of differentiability are Gâteaux (directional) differentiability and Fréchet differentiability. With respect to these definitions, many properties of both Gâteaux and Fréchet differentiability have been achieved. Meanwhile, both Gâteaux and Fréchet differentiability with their properties have been widely applied to both pure and applied mathematics such as, fixed point theory, optimization theory, variational analysis, approximation theory, optimal theory, control theory, and so forth (See [2−4, 7, 28−30, 34]).

In calculus, one can find the explicit derivatives of a differentiable function. It is well-known that the explicit formulas of derivatives of functions play crucial important roles in the applications of differentiation to both pure and applied mathematics, in particular in applied mathematics such as, physics, technology, economics theory, and so forth.

Comparing with the differentiation in calculus, the tasks to prove Gâteaux and Fréchet differentiability, in particular the Fréchet differentiability, of single-valued mappings between Banach spaces are more complicated and more difficult to be accomplished than to test the differentiability of functions in calculus. One more step further, to find the explicit Gâteaux and Fréchet derivatives of differentiable mappings between Banach spaces are extremally complicated and extremally difficult, in general.

It is well-known that the differentiability and the explicit differential formulas (or formulas for derivatives) in calculus are based on special functions. Hence, in order to prove the differentiability and to find some explicit Gâteaux and Fréchet derivatives (formulas) of mappings in Banach spaces, we must have the specifically given mappings and specific underlying Banach spaces. For example, the standard metric projection operator is one of the important operators in the analysis of Banach spaces. The differentiability of the standard metric projection operator has been studied by many authors (See [5, 8−26, 31, 33]). In particular, in [10−12,16−26], the explicit Gâteaux, Fréchet and Mordukhovich derivatives of the standard metric projection operator have been provided, in which the metric projection operator is project to some special nonempty closed and convex subsets of the considered Banach spaces such as, closed balls, closed cylinders and positive cones, and so forth.

Apart from the metric projection operator, the duality mappings, in particular, the normalized duality mapping (denoted by $J$) from Banach spaces to their topological dual spaces, are other ones of important mappings between Banach spaces in the analysis of Banach spaces, in particular in operator theory (See [6, 35]). They have many useful properties, which have been widely applied to fixed point theory, optimization theory, variational analysis, approximation theory, and so forth (See [1, 6, 32, 35]).

In general, the normalized duality mapping $J$ is a set-valued mapping. In [27], the present author studied the Gâteaux, Fréchet and Mordukhovich differentiability of the normalized duality mapping in general Banach spaces, such as $l_p$, for $p \geq 1$, $C[0, 1]$ and so forth. However, there is a gap in the development process of the differentiation studied in [27], which is the lack of explicit solutions for the Gâteaux, Fréchet and Mordukhovich derivatives of the normalized duality mapping.

The purpose of this paper is to enrich the properties of the Gâteaux, Fréchet and Mordukhovich differentiability studied in [27] and to make up for the shortcoming of explicit Gâteaux, Fréchet and Mordukhovich derivatives of the normalized duality mapping $J$. To reach this goal, as we mentioned early, to find the explicit Gâteaux, Fréchet and Mordukhovich derivatives (formulas) of the normalized duality mapping $J$, we must have the specific underlying Banach spaces. Hence, in this paper, we focus in the uniformly convex and uniformly smooth Banach space $(l_p, \|\cdot\|_p)$, for $1 < p < \infty$. More precisely, this paper is organized as follows.

In section 2, we review the concepts and properties of the normalized duality mapping $J$ on Banach space $l_p$ and the concepts of its Gâteaux, Fréchet and Mordukhovich derivatives. We also review the explicit representation of the normalized duality mapping $J$ in $(l_p, \|\cdot\|_p)$, for $1 < p < \infty$, which will be used immediately in the calculation of the derivatives of $J$ in $l_p$.

In section 3, we will accomplish the following goals.

(A) Based on the representation of $J$ on the Banach space $(l_p, \|\cdot\|_p)$, for $2 < p < \infty$, we investigate the Gâteaux differentiability and find the explicit Gâteaux derivatives of $J$ ((ii) of Theorem 3.2);

(B) In part (iii) of Theorem 3.2, we study the Gâteaux directional differentiability and find the explicit Gâteaux directional derivatives of $J$ in $l_p$, for $1 < p < 2$.

In section 4, we first prove that, for $2 < p < \infty$, for any given $\bar{x} \in l_p\backslash\{\theta\}$, the Gâteaux derivative of $J$ at $\bar{x}, J'(\bar{x})(\cdot): l_p \to l_q$ is a linear and Lipchitz mapping from $l_p$ to $l_q$. Then, for $2 < p < \infty$, in Theorem 4.4, we prove the Fréchet differentiability of $J$ and find its explicit Fréchet derivative in a special case that the considered point $\bar{x} \in l_p\backslash\{\theta\}$ has finite number of nonzero entries.

In section 5, by using the connection between Fréchet and Mordukhovich derivatives of single-valued mappings between Banach spaces, as a consequence of Theorem 4.4, we find the explicit Mordukhovich derivative of $J$: $l_p \to l_q$, for $2 < p < \infty$, for point $\bar{x} \in l_p \backslash \{\theta\}$ with finite number of nonzero entries (Theorem 5.1). Then, in this case, we calculate the covering constant for $J$ in $l_p$.

In section 6, we introduce a duality mapping $\mathcal{J}: l_p \to l_q$ and find its explicit Gâteaux, Fréchet and Mordukhovich derivatives, by which, we calculate the covering constant for $\mathcal{J}$ in $l_p$, for $3 < p < \infty$.

In section 7, we introduce a generalized duality mapping $\mathfrak{J}: l_p \to l_q$. Similarly to section 6, we find the explicit Gâteaux, Fréchet and Mordukhovich derivatives of $\mathfrak{J}$, by which, we calculate the covering constant for $\mathfrak{J}$ in $l_p$, for $3 < p < \infty$.

In the conclusion section (section 8), we list some problems for interested readers consideration.

## 2. Review Duality Mappings and Concepts of Differentiation in $l_p$

Throughout this paper, let $p$ and $q$ be positive numbers satisfying $1 < p, q < \infty$ and $\frac{1}{p} + \frac{1}{q} = 1$. In this section, we focus on the real uniformly convex and uniformly smooth Banach spaces $(l_p, \|\cdot\|_p)$ and $(l_q, \|\cdot\|_q)$. The spaces $(l_p, \|\cdot\|_p)$ and $(l_q, \|\cdot\|_q)$ are mutually dual spaces. That is, $l_p^* = l_q$ and $l_p^{**} = l_q^* = l_p$. Both $l_p$ and $l_q$ have the same null element $\theta = (0, 0, \ldots)$. Let $\langle \cdot, \cdot \rangle$ denote the real canonical pairing between $l_q$ and $l_p$.

In general, duality mappings in $l_p$ are defined as follows (See Definitions 3.1 and 3.3 in [6] for the definition of duality mappings on general Banach spaces).

Let $K$: $l_p \to l_q$ be a single-valued mapping. If there is a continuous and strictly increasing function $\varphi: [0, \infty) \to [0, \infty)$ satisfying $\varphi(0) = 0$ and $\varphi(t) \to \infty$, as $t \to \infty$, such that

$$\|K(x)\|_q = \varphi(\|x\|_p), \text{ for any } x \in l_p, \tag{2.1}$$

then, $K$ is called a duality mapping in $l_p$ with gauge function $\varphi$. This definition will be used in sections 6 to define a duality mapping in $l_p$ denoted by $\mathcal{J}: l_p \to l_q$. In this paper, we extend the concept of duality mappings in $l_p$ to generalized duality mappings in $l_p$. Of course, the concept of generalized duality mappings can be extended to general Banach spaces or normed vector spaces.

**Definition 2.1**. Let $K$: $l_p \to l_q$ be a single-valued mapping. If there is a continuous and strictly increasing function $\varphi: [0, \infty) \to [0, \infty)$ satisfying $\varphi(0) = 0$ and $\varphi(t) \to \infty$, as $t \to \infty$, such that

$$\|K(x)\|_q \leq \varphi(\|x\|_p), \text{ for any } x \in l_p, \tag{2.2}$$

then, $K$ is called a generalized duality mapping in $l_p$ with gauge function $\varphi$. This definition will be used in sections 7 to define a generalized duality mapping denoted by $\mathfrak{J}$ in $l_p$.

Especially, the normalized duality mapping $J$: $l_p \to l_q$ is a special duality mapping in $l_p$ with gauge function $\varphi$ defined by $\varphi(t) = t$, for any $t \geq 0$. More precisely, the normalized duality mapping $J$: $l_p \to l_q$ is a continuous single-valued mapping, which is defined, for any $x \in l_p$, by $J(x) \in l_p$ and

$$\langle J(x), x \rangle = \|J(x)\|_q \|x\|_p = \|x\|_p^2 = \|J(x)\|_q^2$$

The normalized duality mapping $J$: $l_p \to l_q$ is has the following explicit representation, which is listed as a lemma without proof.

**Lemma 2.2**. *$J$ and $J^*$ have the following properties.*

(i) $J(\theta) = \theta$;

(ii) *For any $x = (t_1, t_2, \ldots) \in l_p \backslash \{\theta\}$, then*

$$J(x) = \left(\frac{|t_1|^{p-1}\mathrm{sign}(t_1)}{\|x\|_p^{p-2}}, \frac{|t_2|^{p-1}\mathrm{sign}(t_2)}{\|x\|_p^{p-2}}, \ldots\right) = \left(\frac{|t_1|^{p-2}t_1}{\|x\|_p^{p-2}}, \frac{|t_2|^{p-2}t_2}{\|x\|_p^{p-2}}, \ldots\right). \quad (2.3)$$

(iii) $J^*(\theta) = \theta$;

(iv) *For any $y = (s_1, s_2, \ldots) \in l_q \backslash \{\theta\}$, we have*

$$J^*(y) = \left(\frac{|s_1|^{q-1}\mathrm{sign}(s_1)}{\|y\|_q^{q-2}}, \frac{|s_2|^{q-1}\mathrm{sign}(s_2)}{\|y\|_q^{q-2}}, \ldots\right) = \left(\frac{|s_1|^{q-2}s_1}{\|y\|_q^{q-2}}, \frac{|s_2|^{q-2}s_2}{\|y\|_q^{q-2}}, \ldots\right).$$

(v) *In particular, if $p = q = 2$, then $J = J^* = I_{l_2}$, which is the identity mapping in $l_2$.*

(vi) *In general, $J^* \circ J = I_{l_p}$ and $J \circ J^* = I_{l_q}$.*

In particular, when $p > 2$, the indexes $p$ and $q$ of the considered Banach spaces $(l_p, \|\cdot\|_p)$ and $(l_q, \|\cdot\|_q)$ have the following properties, which will be useful in the proofs of the theorems in this section.

**Lemma 2.3**. *Let $p$, $q$ be positive numbers with $2 < p < \infty$ and $\frac{1}{p} + \frac{1}{q} = 1$. Then $p$ and $q$ have the following properties.*

(i) $1 < q < \frac{3}{2}$ *and* $\frac{p}{2q} > 1$;

(ii) $(p-2)q = p - q$, $(p-1)q = p$ *and* $p - 1 - \frac{p}{q} = 0$;

(iii) $0 < \frac{2}{q} - 1 = \frac{1}{q} - \frac{1}{p} < 1$ and $0 < \frac{p-2}{p} < 1$;

(iv) $p - 2q = pq - 3q > 0 \Longrightarrow \frac{p}{2q} > 1$;

(v) $pq - (p - 2q) = 3q \Longrightarrow \frac{pq}{p-2q} > 1$;

(vi) $2p > pq > p - 2q \Longrightarrow \frac{p}{p-2q} > 1$.

*Proof*. The condition $\frac{1}{p} + \frac{1}{q} = 1$ implies $p + q = pq$. This lemma can be straight forwardly checked. □

Now we review the concepts of differentiation of the normalized duality mapping $J$ in $l_p$. Meanwhile, the differentiation of the duality mapping $\mathcal{J}$ and generalized duality mapping $\mathfrak{J}$ in $l_p$ can be similarly defined, which will be given in sections 6 and 7, respectively.

(Gâteaux (directional) differentiability of $J$). Let $\bar{x} \in l_p$ and $v \in l_p \backslash \{\theta\}$. If there is a vector in $l_q$ denoted by $J'(\bar{x}, v)$ such that

$$J'(\bar{x}, v) = \lim_{t \to 0} \frac{J(\bar{x} + tv) - J(\bar{x})}{t},$$

then, $J$ is said to be Gâteaux directionally differentiable at point $\bar{x}$ along direction $v$. The point $J'(\bar{x}, v) \in l_q$ is called the Gâteaux directional derivative of $J$ at point $\bar{x}$ along direction $v$. Furthermore, if $J$ is Gâteaux directionally differentiable at point $\bar{x}$ along every direction $v \in l_p \backslash \{\theta\}$, then $J$ is said to be Gâteaux differentiable at $\bar{x}$ and the Gâteaux derivative of $J$ at point $\bar{x}$ is denoted by

$$J'(\bar{x})(v) = \lim_{t \to 0} \frac{J(\bar{x}+tv)-J(\bar{x})}{t}, \text{ for any } v \in l_p \backslash \{\theta\}.$$

(Fréchet differentiability of $J$). If there is a continuous and linear mapping $\nabla J(\bar{x}): l_p \to l_q$ such that,

$$\lim_{u \xrightarrow{l_p} \theta} \frac{J(\bar{x}+u)-J(\bar{x})-\nabla J(\bar{x})(u)}{\|u\|_p} = \theta,$$

then $J$ is said to be Fréchet differentiable at $\bar{x}$ and $\nabla J(\bar{x})$ is called the Fréchet derivative of $J$ at $\bar{x}$. It is well-known that in Banach spaces (or in normed vector spaces), for single-valued mappings, Fréchet differentiability is stronger than Gâteaux differentiability. That is,

$$J \text{ is Fréchet differentiable at } \bar{x} \quad \Longrightarrow \quad J \text{ is Gâteaux differentiable at } \bar{x}. \tag{2.4}$$

More precisely, the connection between Gâteaux differentiability and Fréchet differentiability in Banach spaces is that if $J$ is Fréchet differentiable at $\bar{x}$, then $J$ is also Gâteaux differentiable at $\bar{x}$ and

$$J'(\bar{x})(v) = \nabla J(\bar{x})(v), \text{ for any } v \in l_p \backslash \{\theta\}. \tag{2.5}$$

It is well-known that the inverse of (2.4) is not true. This paper actually provides some counter examples to demonstrate the following fact (Theorem 3.2 and Theorem 4.1).

$$J \text{ is Gâteaux differentiable at } \bar{x} \quad \nRightarrow \quad J \text{ is Fréchet differentiable at } \bar{x}.$$

(Mordukhovich derivative of $J$, or coderivative, or Mordukhovich coderivative, or Fréchet coderivative of $J$) Notice that, for $1 < p < \infty$, $J$: $l_p \to l_q$ is a single valued continuous mapping. Then, the Mordukhovich derivative of $J$ at point $(\bar{x}, J(\bar{x}))$ is a set-valued mapping $\widehat{D}^*J(\bar{x}, J(\bar{x})): (l_q)^* \rightrightarrows (l_p)^*$. By $(l_q)^* = l_p$ and $(l_p)^* = l_q$, more precisely, we have $\widehat{D}^*J(\bar{x}, J(\bar{x})): l_p \rightrightarrows l_q$. Hence, for any $y \in l_p$, it is defined by

$$\widehat{D}^*J(\bar{x}, J(\bar{x}))(y) \coloneqq \widehat{D}^*J(\bar{x})(y) = \left\{ z \in l_q : \limsup_{\substack{u \to \bar{x} \\ u \in l_p}} \frac{\langle z,\ u-\bar{x}\rangle - \langle y,\ J(u)-J(\bar{x})\rangle}{\|u-\bar{x}\|_p + \|J(u)-J(\bar{x})\|_q} \le 0 \right\}.$$

As a matter of fact, in $l_p$, for $1 < p < \infty$, $\widehat{D}^*J(\bar{x}, J(\bar{x})): l_p \to l_q$ is a single-valued mapping.

The following theorem shows the connection between Fréchet derivatives and Mordukhovich derivatives for single-valued mappings. Its results provide a powerful tool to calculate the Mordukhovich derivatives by the Fréchet derivatives of single-valued mappings.

**Theorem 1.38 in [28]**. *Let X and Y be Banach spaces with dual space $X^*$ and $Y^*$*, respectively. *Let $\varphi$: X → Y be a single-valued mapping. Suppose that $\varphi$ is Fréchet differentiable at $x \in X$ with $y = \varphi(x)$. Then, the Mordukhovich derivative of $\varphi$ at x satisfies the following equation*

$$\widehat{D}^*\varphi(x, y)(y^*) = \{(\nabla\varphi(x))^*(y^*)\}, \ \textit{for all } y^* \in Y^*.$$

One of the most important applications of the Mordukhovich derivatives of single-valued mappings is to define the covering constants, which is recalled as follows. Let $\varphi: X \to Y$ be a single-valued mapping. For given $\bar{x} \in X$ and $\bar{y} \in Y$ with $\bar{y} = \varphi(\bar{x})$, the covering constant for $\varphi$ at $(\bar{x}, \bar{y})$ is defined by

$$\hat{\alpha}(\varphi,\bar{x},\bar{y}) = \sup_{\eta>0} \inf\left\{\|z^*\|_{X^*}: z^* \in \widehat{D}^*\varphi(x,y)(w^*), x \in \mathbb{B}_X(\bar{x},\eta), y = \varphi(x) \in \mathbb{B}_Y(\bar{y},\eta), \|w^*\|_{Y^*} = 1\right\}. \quad (2.6)$$

As a special case of Theorem 1.38 in [28], we have the following connection between Fréchet derivatives and Mordukhovich derivatives for the single-valued mapping $J: l_p \to l_q$. More precisely, let $\bar{x} \in l_p$ and $\bar{y} \in l_q$ with $\bar{y} = J(\bar{x})$. Then, for $x^* \in l_q$, we have

$$x^* = \widehat{D}^*J(\bar{x},\bar{y})(y^*) \quad \Longleftrightarrow \quad \limsup_{\substack{u \to \bar{x} \\ l_p}} \frac{\langle x^*, u-\bar{x}\rangle - \langle y^*, J(u)-\bar{y}\rangle}{\|u-\bar{x}\|_p + \|J(u)-\bar{y}\|_q} \le 0.$$

## 3. The Gâteaux Differentiability of the Normalized Duality Mapping $J$ in $l_p$

As mentioned in section 2, let $1 < p, q < \infty$ with $\frac{1}{p} + \frac{1}{q} = 1$. In this section, we consider the Gâteaux differentiability of the normalized duality mapping $J$ between the real uniformly convex and uniformly smooth Banach spaces $(l_p, \|\cdot\|_p)$ and $(l_q, \|\cdot\|_q)$. Let $\bar{x} = (\bar{t}_1, \bar{t}_2, \dots) \in l_p\backslash\{\theta\}$ and $v = \{s_n\}_{n=1}^{\infty} \in l_p\backslash\{\theta\}$. In case, if $J$ is Gâteaux differentiable at $\bar{x}$ along the direction $v$, we will find the explicit solution for $J'(\bar{x}, v)$.

In this paper, let $\mathbb{N}$ denote the set of positive integers. Let $\mathbb{S}$ denote the linear space of real sequences with term-by-term addition and scalar multiplication. Let $x = (t_1, t_2, \dots) \in \mathbb{S}$. We write

$$\mathcal{N}(x) = \{n \in \mathbb{N}: t_n \neq 0\}.$$

Let $|\mathcal{N}(x)|$ denote the cardinality of $\mathcal{N}(x)$. Next, we define the concept of the ratio sequences.

**Definition 3.1**. Let arbitrarily given $\bar{x} = (\bar{t}_1, \bar{t}_2, \dots) \in \mathbb{S}$. Define an operator $R(\bar{x})(\cdot): \mathbb{S} \to \mathbb{S}$, for $v = \{s_n\}_{n=1}^{\infty} \in \mathbb{S}$ by $R(\bar{x})(v) \coloneqq \{\lambda_n\}_{n=1}^{\infty}$, in which, for each $n = 1, 2, \dots$, $\lambda_n \in \mathbb{R}$ and

$$\lambda_n = \begin{cases} 0, & \text{if } \bar{t}_n = 0, \\ \frac{s_n}{\bar{t}_n}, & \text{if } \bar{t}_n \neq 0. \end{cases}$$

The sequence $R(\bar{x})(v) = \{\lambda_n\}_{n=1}^{\infty} \in \mathbb{S}$ is called the ratio sequence of $v$ with respect to $\bar{x}$. More precisely, $\lambda_n$ satisfies that

$$\lambda_n = \begin{cases} 0, & \text{if } \bar{t}_n = 0, \\ \frac{s_n}{\bar{t}_n}, & \text{if } \bar{t}_n \neq 0 \text{ and } s_n \neq 0, \\ 0, & \text{if } \bar{t}_n \neq 0 \text{ and } s_n = 0. \end{cases} \quad \text{for } n = 1, 2, \dots.$$

This implies that, for some $n \in \mathbb{N}$, if $\bar{t}_n \neq 0$ and $s_n = 0$, then $\lambda_n = 0$. This implies that

$$s_n = \lambda_n \bar{t}_n, \text{ if } \bar{t}_n \neq 0, \text{ for } n \in \mathbb{N}.$$

Notice that, for a given $\bar{x} \in \mathbb{S}$,

$$R(\bar{x})(\cdot): \mathbb{S} \to \mathbb{S} \text{ is a linear operator on } \mathbb{S}. \quad (3.1)$$

In particular, if $v$ satisfies that $s_n\bar{t}_n = 0$, for all $n$ = 1, 2, …, then $R(\bar{x}, v) = \theta$. If $\{\lambda_n\}_{n=1}^{\infty}$ is a bounded sequence and there is a positive number $\Lambda$ such that

$$|\lambda_n| \leq \Lambda, \text{ for } n = 1, 2, \ldots,$$

then $v$ is said to be ratio bounded and with ratio bounds $\Lambda$ with respect to $\bar{x}$.

**Theorem 3.2**. *The normalized duality mapping $J$: $l_p \to l_q$ has the following differentiation properties.*

(i) *For $1 < p < \infty$, $J$ is Gâteaux differentiable at $\theta$ and*

$$J'(\theta)(v) = J(v), \text{ for any } v \in l_p\backslash\{\theta\}. \tag{3.2}$$

(ii) *For $2 < p < \infty$, $J$ is Gâteaux differentiable on $l_p\backslash\{\theta\}$. More precisely, let $\bar{x} = (\bar{t}_1, \bar{t}_2, \ldots) \in l_p\backslash\{\theta\}$. For any $v = \{s_n\}_{n=1}^{\infty} \in l_p\backslash\{\theta\}$ with ratio sequence $\{\lambda_n\}_{n=1}^{\infty}$ with respect to $\bar{x}$, the Gâteaux derivative of $J$ at $\bar{x}$ in direction $v$ is*

$$J'(\bar{x})(v) = J(\bar{x})\left((p-1)\begin{pmatrix}\lambda_1 & 0 & \cdots \\ 0 & \lambda_2 & 0 \\ \vdots & 0 & \ddots\end{pmatrix} - \frac{(p-2)\sum_{i=1}^{\infty}\lambda_i|\bar{t}_i|^p}{\sum_{i=1}^{\infty}|\bar{t}_i|^p}\begin{pmatrix}1 & 0 & \cdots \\ 0 & 1 & 0 \\ \vdots & 0 & \ddots\end{pmatrix}\right). \tag{3.3}$$

*In particular, we have that*

(a) *If $\mathcal{N}(\bar{x}) = 1$ with $\bar{t}_m \neq 0$, then,*

$$J'(\bar{x})(v) = J(\bar{x})\begin{pmatrix}0 & 0 & \cdots \\ 0 & \lambda_m & 0 \\ \vdots & 0 & \ddots\end{pmatrix} = \bar{x}\begin{pmatrix}0 & 0 & \cdots \\ 0 & \lambda_m & 0 \\ \vdots & 0 & \ddots\end{pmatrix}. \tag{3.4}$$

(b) *If $s_n\bar{t}_n = 0$, for all $n$ = 1, 2, …, then*

$$J'(\bar{x})(v) = \theta.$$

(iii) *Let $1 < p < 2$. We have that*

(a) *If $\mathcal{N}(v) \nsubseteq \mathcal{N}(\bar{x})$, then $J$ is not Gâteaux directionally differentiable at $\bar{x}$ in direction $v$.*

$$\mathcal{N}(v) \nsubseteq \mathcal{N}(\bar{x}) \quad \Longrightarrow \quad J'(\bar{x}, v) \text{ does not exist}. \tag{3.5}$$

(b) *If $\mathcal{N}(v) \subseteq \mathcal{N}(\bar{x})$ and $v$ is ratio bounded with respect to $\bar{x}$, then $J$ is Gâteaux directionally differentiable at $\bar{x}$ in direction $v$ and $J'(\bar{x}, v)$ satisfies* (3.3).

*In particular,*

($b_1$) *If $\mathcal{N}(v) \subseteq \mathcal{N}(\bar{x})$ and $0 < |\mathcal{N}(v)| < \infty$, then $J$ is Gâteaux directionally differentiable at $\bar{x}$ in direction $v$ and $J'(\bar{x}, v)$ satisfies* (3.3).

($b_2$) *If $\mathcal{N}(v) \subseteq \mathcal{N}(\bar{x})$ and $0 < |\mathcal{N}(\bar{x})| < \infty$, then $J$ is Gâteaux directionally differentiable at $\bar{x}$ in direction $v$ and $J'(\bar{x}, v)$ satisfies* (3.3).

(iv) *Let $1 < p < \infty$. For any $\bar{x} \in l_p\backslash\{\theta\}$, $J$ is Gâteaux directionally differentiable at $\bar{x}$ in direction $\bar{x}$ and*

$$J'(\bar{x}, \bar{x}) = J(\bar{x}), \textit{ for any } \bar{x} \in l_p \backslash \{\theta\}.$$

*In particular*,

$$J'(\bar{x}, \bar{x}) = \bar{x}, \textit{ for any } \bar{x} \in l_p \textit{ with } |\mathcal{N}(\bar{x})| = 1.$$

*Proof.* Proof of part (i). We first prove (3.2). Let $v = \{s_n\}_{n=1}^{\infty} \in l_p \backslash \{\theta\}$. For any real $t \neq 0$, we calculate

$$J(tv) = \left( \frac{|ts_1|^{p-1}\mathrm{sign}(ts_1)}{\|tv\|_p^{p-2}}, \; \frac{|ts_2|^{p-1}\mathrm{sign}(ts_2)}{\|tv\|_p^{p-2}}, \dots \right)$$

$$= \left( \frac{|t|^{p-1}|s_1|^{p-1}\mathrm{sign}(ts_1)}{|t|^{p-2}\|v\|_p^{p-2}}, \; \frac{|t|^{p-1}|s_2|^{p-1}\mathrm{sign}(ts_2)}{|t|^{p-2}\|v\|_p^{p-2}}, \dots \right)$$

$$= \left( \frac{|t||s_1|^{p-1}\mathrm{sign}(ts_1)}{\|v\|_p^{p-2}}, \; \frac{|t||s_2|^{p-1}\mathrm{sign}(ts_2)}{\|v\|_p^{p-2}}, \dots \right)$$

$$= \begin{cases} \left( \frac{|t||s_1|^{p-1}\mathrm{sign}(s_1)}{\|v\|_p^{p-2}}, \; \frac{|t||s_2|^{p-1}\mathrm{sign}(s_2)}{\|v\|_p^{p-2}}, \dots \right), \text{for } t > 0, \\ -\left( \frac{|t||s_1|^{p-1}\mathrm{sign}(s_1)}{\|v\|_p^{p-2}}, \; \frac{|t||s_2|^{p-1}\mathrm{sign}(s_2)}{\|v\|_p^{p-2}}, \dots \right), \text{for } t < 0 \end{cases}$$

$$= \begin{cases} |t|J(v), \text{for } t > 0, \\ -|t|J(v), \text{for } t < 0. \end{cases}$$

By $J(\theta) = \theta$, this implies that, for any $v \in l_p \backslash \{\theta\}$,

$$J'(\theta)(v) = \lim_{t \to 0} \frac{J(\theta + tv) - J(\theta)}{t} = \lim_{t \to 0} \frac{J(tv)}{t} = \lim_{t \to 0} \frac{|t|J(v)}{|t|} = J(v).$$

By definition of Gâteaux derivative, this proves that *J* is Gâteaux differentiable at $\theta$ and (3.2) is satisfied.

Proof of (ii). $2 < p < \infty$ and $\frac{1}{p} + \frac{1}{q} = 1$. Let $\bar{x} = (\bar{t}_1, \bar{t}_2, \dots) \in l_p \backslash \{\theta\}$. Let $v = \{s_n\}_{n=1}^{\infty} \in l_p \backslash \{\theta\}$ with ratio sequence $\{\lambda_n\}_{n=1}^{\infty}$ with respect to $\bar{x}$. We calculate

$$J'(\bar{x})(v) = \lim_{t \to 0} \frac{J(\bar{x} + tv) - J(\bar{x})}{t}$$

$$= \lim_{t \to 0} \frac{\left\{ \frac{|\bar{t}_n + ts_n|^{p-1}\mathrm{sign}(\bar{t}_n + ts_n)}{\left(\sum_{i=1}^{\infty} |\bar{t}_i + ts_i|^p\right)^{\frac{p-2}{p}}} \right\}_{n=1}^{\infty} - \left\{ \frac{|\bar{t}_n|^{p-1}\mathrm{sign}(\bar{t}_n)}{\left(\sum_{i=1}^{\infty} |\bar{t}_i|^p\right)^{\frac{p-2}{p}}} \right\}_{n=1}^{\infty}}{t} =$$

$$= \lim_{t \to 0} \left\{ \frac{\dots, \; \frac{|ts_n|^{p-1}\mathrm{sign}(ts_n)}{\left(\sum_{i \in \mathcal{N}(\bar{x})} |1 + t\lambda_i|^p |\bar{t}_i|^p + |t|^p \sum_{j \notin \mathcal{N}(\bar{x})} |s_j|^p\right)^{\frac{p-2}{p}}}, \dots,}{t} \right.$$

$$\left.\frac{\frac{|1+t\lambda_m|^{p-1}|\bar{t}_m|^{p-1}\operatorname{sign}(\bar{t}_m+t\lambda_m\bar{t}_m)}{\left(\sum_{i\in\mathcal{N}(\bar{x})}|1+t\lambda_i|^p|\bar{t}_i|^p+|t|^p\sum_{j\notin\mathcal{N}(\bar{x})}|s_j|^p\right)^{\frac{p-2}{p}}}-\frac{|\bar{t}_m|^{p-1}\operatorname{sign}(\bar{t}_m)}{\left(\sum_{i\in\mathcal{N}(\bar{x})}|\bar{t}_i|^p+\sum_{j\notin\mathcal{N}(\bar{x})}|\bar{t}_j|^p\right)^{\frac{p-2}{p}}},\cdots}{t}\right\}$$

$$=\lim_{t\to 0}\frac{\left\{\cdots,\frac{|ts_n|^{p-1}\operatorname{sign}(ts_n)}{\left(\sum_{i\in\mathcal{N}(\bar{x})}|1+t\lambda_i|^p|\bar{t}_i|^p+|t|^p\sum_{j\notin\mathcal{N}(\bar{x})}|s_j|^p\right)^{\frac{p-2}{p}}},\ \cdots,\ \frac{|1+t\lambda_m|^{p-1}|\bar{t}_m|^{p-1}\operatorname{sign}(\bar{t}_m+t\lambda_m\bar{t}_m)}{\left(\sum_{i\in\mathcal{N}(\bar{x})}|1+t\lambda_i|^p|\bar{t}_i|^p+|t|^p\sum_{j\notin\mathcal{N}(\bar{x})}|s_j|^p\right)^{\frac{p-2}{p}}}-\frac{|\bar{t}_m|^{p-1}\operatorname{sign}(\bar{t}_m)}{\left(\sum_{i\in\mathcal{N}(\bar{x})}|\bar{t}_i|^p\right)^{\frac{p-2}{p}}},\cdots\right\}}{t}$$

Here, the *n*th entry of $\bar{x}$ is zero and *n*th entry of $\bar{x}$ is not zero. Then, the above limit becomes that

$$=\lim_{t\to 0}\frac{\left\{\cdots,\frac{|ts_n|^{p-1}\operatorname{sign}(ts_n)}{\left(\sum_{i\in\mathcal{N}(\bar{x})}|1+t\lambda_i|^p|\bar{t}_i|^p+|t|^p\sum_{j\notin\mathcal{N}(\bar{x})}|s_j|^p\right)^{\frac{p-2}{p}}},\ \cdots,\ \frac{|1+t\lambda_m|^{p-1}|\bar{t}_m|^{p-1}\operatorname{sign}((1+t\lambda_m)\bar{t}_m)}{\left(\sum_{i\in\mathcal{N}(\bar{x})}|1+t\lambda_i|^p|\bar{t}_i|^p+|t|^p\sum_{j\notin\mathcal{N}(\bar{x})}|s_j|^p\right)^{\frac{p-2}{p}}}-\frac{|\bar{t}_m|^{p-1}\operatorname{sign}(\bar{t}_m)}{\left(\sum_{i\in\mathcal{N}(\bar{x})}|\bar{t}_i|^p\right)^{\frac{p-2}{p}}},\cdots\right\}}{t}$$

$$=\lim_{t\to 0}\frac{\left\{\cdots,\frac{|ts_n|^{p-1}\operatorname{sign}(ts_n)}{\left(\sum_{i\in\mathcal{N}(\bar{x})}(1+t\lambda_i)^p|\bar{t}_i|^p+|t|^p\sum_{j\notin\mathcal{N}(\bar{x})}|s_j|^p\right)^{\frac{p-2}{p}}},\ \cdots,\left(\frac{(1+t\lambda_m)^{p-1}|\bar{t}_m|^{p-1}\operatorname{sign}((1+t\lambda_m)\bar{t}_m)}{\left(\sum_{i\in\mathcal{N}(\bar{x})}(1+t\lambda_i)^p|\bar{t}_i|^p+|t|^p\sum_{j\notin\mathcal{N}(\bar{x})}|s_j|^p\right)^{\frac{p-2}{p}}}-\frac{|\bar{t}_m|^{p-1}\operatorname{sign}(\bar{t}_m)}{\left(\sum_{i\in\mathcal{N}(\bar{x})}|\bar{t}_i|^p\right)^{\frac{p-2}{p}}}\right),\cdots\right\}}{t}$$

For a given fixed $m\in\mathcal{N}(\bar{x})$, let

$$f_m(t)=\frac{(1+t\lambda_m)^{p-1}|\bar{t}_m|^{p-1}\operatorname{sign}((1+t\lambda_m)\bar{t}_m)}{\left(\sum_{i\in\mathcal{N}(\bar{x})}(1+t\lambda_i)^p|\bar{t}_i|^p+|t|^p\sum_{j\notin\mathcal{N}(\bar{x})}|s_j|^p\right)^{\frac{p-2}{p}}}.$$

$f_m$ is a function of *t* and satisfies that

$$f_m(0)=\frac{|\bar{t}_m|^{p-1}\operatorname{sign}(\bar{t}_m)}{\left(\sum_{i\in\mathcal{N}(\bar{x})}|\bar{t}_i|^p\right)^{\frac{p-2}{p}}}.$$

Notice that

$$\operatorname{sign}((1+t\lambda_m)\bar{t}_m)=\begin{cases}\operatorname{sign}(\bar{t}_m), & \text{if } 1+t\lambda_m>0,\\ 0, & \text{if } 1+t\lambda_m=0,\\ -\operatorname{sign}(\bar{t}_m), & \text{if } 1+t\lambda_m<0.\end{cases}$$

This function $\operatorname{sign}((1+t\lambda_m)\bar{t}_m)$ is a piecewise constant function with possible values $\{-1, 0\ 1\}$. Hence it is considered as a constant when we find the derivative of $f_m$. We calculate

$$f_m'(t)=\left(\frac{(1+t\lambda_m)^p|\bar{t}_m|^{p-1}\operatorname{sign}((1+t\lambda_m)\bar{t}_m)}{\left(\sum_{i\in\mathcal{N}(\bar{x})}(1+t\lambda_i)^p|\bar{t}_i|^p+|t|^p\sum_{j\notin\mathcal{N}(\bar{x})}|s_j|^p\right)^{\frac{p-2}{p}}}\right)'$$

$$=\operatorname{sign}((1+t\lambda_m)\bar{t}_m)\left(\frac{(1+t\lambda_m)^{p-1}|\bar{t}_m|^{p-1}}{\left(\sum_{i\in\mathcal{N}(\bar{x})}(1+t\lambda_i)^p|\bar{t}_i|^p+|t|^p\sum_{j\notin\mathcal{N}(\bar{x})}|s_j|^p\right)^{\frac{p-2}{p}}}\right)'$$

$$= \mathrm{sign}((1+t\lambda_m)\bar{t}_m)\left(\frac{(p-1)\lambda_m(1+t\lambda_m)^{p-2}|\bar{t}_m|^{p-1}\left(\sum_{i\in\mathcal{N}(\bar{x})}(1+t\lambda_i)^p|\bar{t}_i|^p+|t|^p\sum_{j\notin\mathcal{N}(\bar{x})}|s_j|^p\right)}{\left(\sum_{i\in\mathcal{N}(\bar{x})}(1+t\lambda_i)^p|\bar{t}_i|^p+|t|^p\sum_{j\notin\mathcal{N}(\bar{x})}|s_j|^p\right)^{\frac{p-2}{p}+1}}\right.$$

$$\left. - \frac{(1+t\lambda_m)^{p-1}|\bar{t}_m|^{p-1}\frac{p-2}{p}\left(\sum_{i\in\mathcal{N}(\bar{x})}p\lambda_i(1+t\lambda_i)^{p-1}|\bar{t}_i|^p+p|t|^{p-1}\mathrm{sign}(t)\sum_{j\notin\mathcal{N}(\bar{x})}|s_j|^p\right)}{\left(\sum_{i\in\mathcal{N}(\bar{x})}(1+t\lambda_i)^p|\bar{t}_i|^p+|t|^p\sum_{j\notin\mathcal{N}(\bar{x})}|s_j|^p\right)^{\frac{p-2}{p}+1}}\right)$$

$$= \mathrm{sign}\big((1+t\lambda_m)\bar{t}_m\big)\left(\frac{(p-1)\lambda_m(1+t\lambda_m)^{p-2}|\bar{t}_m|^{p-1}\left(\sum_{i\in\mathcal{N}(\bar{x})}(1+t\lambda_i)^p|\bar{t}_i|^p+|t|^p\sum_{j\notin\mathcal{N}(\bar{x})}|s_j|^p\right)}{\left(\sum_{i\in\mathcal{N}(\bar{x})}(1+t\lambda_i)^p|\bar{t}_i|^p+|t|^p\sum_{j\notin\mathcal{N}(\bar{x})}|s_j|^p\right)^{\frac{p-2}{p}+1}}\right.$$

$$\left. - \frac{(p-2)(1+t\lambda_m)^{p-1}|\bar{t}_m|^{p-1}\left(\sum_{i\in\mathcal{N}(\bar{x})}\lambda_i(1+t\lambda_i)^{p-1}|\bar{t}_i|^p+|t|^{p-1}\mathrm{sign}(t)\sum_{j\notin\mathcal{N}(\bar{x})}|s_j|^p\right)}{\left(\sum_{i\in\mathcal{N}(\bar{x})}(1+t\lambda_i)^p|\bar{t}_i|^p+|t|^p\sum_{j\notin\mathcal{N}(\bar{x})}|s_j|^p\right)^{\frac{p-2}{p}+1}}\right)$$

$$= |\bar{t}_m|^{p-1}\mathrm{sign}\big((1+t\lambda_m)\bar{t}_m\big)\left(\frac{(p-1)\lambda_m(1+t\lambda_m)^{p-2}\left(\sum_{i\in\mathcal{N}(\bar{x})}(1+t\lambda_i)^p|\bar{t}_i|^p+|t|^p\sum_{j\notin\mathcal{N}(\bar{x})}|s_j|^p\right)}{\left(\sum_{i\in\mathcal{N}(\bar{x})}(1+t\lambda_i)^p|\bar{t}_i|^p+|t|^p\sum_{j\notin\mathcal{N}(\bar{x})}|s_j|^p\right)^{\frac{p-2}{p}+1}}\right.$$

$$\left. - \frac{(p-2)(1+t\lambda_m)^{p-1}\left(\sum_{i\in\mathcal{N}(\bar{x})}\lambda_i(1+t\lambda_i)^{p-1}|\bar{t}_i|^p+|t|^{p-1}\mathrm{sign}(t)\sum_{j\notin\mathcal{N}(\bar{x})}|s_j|^p\right)}{\left(\sum_{i\in\mathcal{N}(\bar{x})}(1+t\lambda_i)^p|\bar{t}_i|^p+|t|^p\sum_{j\notin\mathcal{N}(\bar{x})}|s_j|^p\right)^{\frac{p-2}{p}+1}}\right).$$

By Taylor's Approximating Theorem for $f_m$ with order 1 around point $t = 0$ (or the Mean Value Theorem) there is some $t_m$ with $|t_m| < |t|$ such that

$$f_m(t) - f_m(0) = f_m'(t_m)t, \text{ for each } m \in \mathcal{N}(\bar{x}).$$

Since we take limit for $t \to 0$, here, we suppose $|t_m| < |t| < 1$. Then $f_m'(t_m)$ satisfies that

$$f_m'(t_m)$$

$$= |\bar{t}_m|^{p-1}\mathrm{sign}\big((1+t_m\lambda_m)\bar{t}_m\big)\left(\frac{(p-1)\lambda_m(1+t_m\lambda_m)^{p-2}\left(\sum_{i\in\mathcal{N}(\bar{x})}(1+t_m\lambda_i)^p|\bar{t}_i|^p+|t_m|^p\sum_{j\notin\mathcal{N}(\bar{x})}|s_j|^p\right)}{\left(\sum_{i\in\mathcal{N}(\bar{x})}(1+t_m\lambda_i)^p|\bar{t}_i|^p+|t_m|^p\sum_{j\notin\mathcal{N}(\bar{x})}|s_j|^p\right)^{\frac{p-2}{p}+1}}\right.$$

$$\left. - \frac{(p-2)(1+t_m\lambda_m)^{p-1}\left(\sum_{i\in\mathcal{N}(\bar{x})}\lambda_i(1+t_m\lambda_i)^{p-1}|\bar{t}_i|^p+|t_m|^{p-1}\mathrm{sign}(t_m)\sum_{j\notin\mathcal{N}(\bar{x})}|s_j|^p\right)}{\left(\sum_{i\in\mathcal{N}(\bar{x})}(1+t_m\lambda_i)^p|\bar{t}_i|^p+|t_m|^p\sum_{j\notin\mathcal{N}(\bar{x})}|s_j|^p\right)^{\frac{p-2}{p}+1}}\right).$$

Substituting this into the last limit, we obtain that

$$J'(\bar{x})(v) = \lim_{t\to 0}\frac{\left\{\ldots, \frac{|ts_n|^{p-1}\mathrm{sign}(ts_n)}{\left(\sum_{i\in\mathcal{N}(\bar{x})}(1+t\lambda_i)^p|\bar{t}_i|^p+|t|^p\sum_{j\notin\mathcal{N}(\bar{x})}|s_j|^p\right)^{\frac{p-2}{p}}}, \ldots, (f_m(t)-f_m(0)), \ldots\right\}}{t}$$

$$= \lim_{t\to 0} \frac{\left\{\ldots, \ \frac{|ts_n|^{p-1}\mathrm{sign}(ts_n)}{\left(\sum_{i\in\mathcal{N}(\bar{x})}(1+t\lambda_i)^p|\bar{t}_i|^p+|t|^p\sum_{j\notin\mathcal{N}(\bar{x})}|s_j|^p\right)^{\frac{p-2}{p}}}, \ \ldots, f_m'(t_m)t, \ldots\right\}}{t}$$

$$= \lim_{t\to 0}\left\{\ldots, \ \frac{|t|^{p-2}|s_n|^{p-1}\mathrm{sign}(ts_n)}{\left(\sum_{i\in\mathcal{N}(\bar{x})}(1+t\lambda_i)^p|\bar{t}_i|^p+|t|^p\sum_{j\notin\mathcal{N}(\bar{x})}|s_j|^p\right)^{\frac{p-2}{p}}}, \ \ldots, \ f'(t_m), \ldots\right\}. \tag{3.6}$$

By $|t_m| < |t| < 1$, we have

$$\lim_{t\to 0} f_m'(t_m)$$

$$= \lim_{t\to 0} |\bar{t}_m|^{p-1}\mathrm{sign}\big((1+t_m\lambda_m)\bar{t}_m\big)\left(\frac{(p-1)\lambda_m(1+t_m\lambda_m)^{p-2}\left(\sum_{i\in\mathcal{N}(\bar{x})}(1+t_m\lambda_i)^p|\bar{t}_i|^p+|t_m|^p\sum_{j\notin\mathcal{N}(\bar{x})}|s_j|^p\right)}{\left(\sum_{i\in\mathcal{N}(\bar{x})}(1+t_m\lambda_i)^p|\bar{t}_i|^p+|t_m|^p\sum_{j\notin\mathcal{N}(\bar{x})}|s_j|^p\right)^{\frac{p-2}{p}+1}}\right.$$

$$\left. - \frac{(p-2)(1+t_m\lambda_m)^{p-1}\left(\sum_{i\in\mathcal{N}(\bar{x})}\lambda_i(1+t_i\lambda_i)^{p-1}|\bar{t}_i|^p+|t_m|^{p-1}\mathrm{sign}(t_m)\sum_{j\notin\mathcal{N}(\bar{x})}|s_j|^p\right)}{\left(\sum_{i\in\mathcal{N}(\bar{x})}(1+t_m\lambda_i)^p|\bar{t}_i|^p+|t_m|^p\sum_{j\notin\mathcal{N}(\bar{x})}|s_j|^p\right)^{\frac{p-2}{p}+1}}\right)$$

$$= |\bar{t}_m|^{p-1}\mathrm{sign}(\bar{t}_m)\left(\frac{(p-1)\lambda_m\sum_{i\in\mathcal{N}(\bar{x})}|\bar{t}_i|^p}{\left(\sum_{i\in\mathcal{N}(\bar{x})}|\bar{t}_i|^p\right)^{\frac{p-2}{p}+1}} - \frac{(p-2)\sum_{i\in\mathcal{N}(\bar{x})}\lambda_i|\bar{t}_i|^p}{\left(\sum_{i\in\mathcal{N}(\bar{x})}|\bar{t}_i|^p\right)^{\frac{p-2}{p}+1}}\right)$$

$$= |\bar{t}_m|^{p-1}\mathrm{sign}(\bar{t}_m)\left(\frac{(p-1)\lambda_m}{\left(\sum_{i\in\mathcal{N}(\bar{x})}|\bar{t}_i|^p\right)^{\frac{p-2}{p}}} - \frac{(p-2)\sum_{i\in\mathcal{N}(\bar{x})}\lambda_i|\bar{t}_i|^p}{\left(\sum_{i\in\mathcal{N}(\bar{x})}|\bar{t}_i|^p\right)^{\frac{p-2}{p}+1}}\right)$$

$$= \frac{|\bar{t}_m|^{p-1}\mathrm{sign}(\bar{t}_m)}{\left(\sum_{i=1}^{\infty}|\bar{t}_i|^p\right)^{\frac{p-2}{p}}}\left((p-1)\lambda_m - \frac{(p-2)\sum_{i=1}^{\infty}\lambda_i|\bar{t}_i|^p}{\sum_{i=1}^{\infty}|\bar{t}_i|^p}\right).$$

Notice that in order to ensure that the above solution of $\lim_{t\to 0} f_m'(t_m)$ is well-defined, we need to show the following two statements (I) and (II), for any $\bar{x} = (\bar{t}_1, \bar{t}_2, \ldots) \in l_p\backslash\{\theta\}$ and $v = \{s_n\}_{n=1}^{\infty} \in l_p\backslash\{\theta\}$ with ratio sequence $\{\lambda_n\}_{n=1}^{\infty}$ with respect to $\bar{x}$:

(I) $\quad |\sum_{i=1}^{\infty}\lambda_i|\bar{t}_i|^p| < \infty$.

We estimate it as follows.

$$|\textstyle\sum_{i=1}^{\infty}\lambda_i|\bar{t}_i|^p| \le \sum_{i=1}^{\infty}|\lambda_i|\bar{t}_i|^p| = \sum_{i=1}^{\infty}|\lambda_i\bar{t}_i|\,|\bar{t}_i|^{p-1} = \sum_{i=1}^{\infty}|s_i|\,|\bar{t}_i|^{p-1}$$

$$\le (\textstyle\sum_{i=1}^{\infty}|s_i|^p)^{\frac{1}{p}}(\sum_{i=1}^{\infty}|\bar{t}_i|^{(p-1)q})^{\frac{1}{q}} = \|v\|_p(\sum_{i=1}^{\infty}|\bar{t}_i|^p)^{\frac{1}{q}} = \|v\|_p\|\bar{x}\|_p^{\frac{p}{q}} < \infty.$$

(II) $\quad \left\{\frac{|\bar{t}_n|^{p-1}\mathrm{sign}(\bar{t}_n)(p-1)\lambda_n}{\left(\sum_{i=1}^{\infty}|\bar{t}_i|^p\right)^{\frac{p-2}{p}}}\right\}_{n=1}^{\infty} \in l_q.$

To this end, for $2 < p < \infty$, we only estimate its numerator as follows without the constant $p-1$.

$$\sum_{n=1}^{\infty}\left||\bar{t}_n|^{p-1}\mathrm{sign}(\bar{t}_n)\lambda_n\right|^q$$

$$= \sum_{n=1}^{\infty}\left||\bar{t}_n|^{p-2}|\bar{t}_n|\mathrm{sign}(\bar{t}_n)\lambda_n\right|^q$$

$$= \sum_{n=1}^{\infty}\left||\bar{t}_n|^{p-2}\bar{t}_n\lambda_n\right|^q$$

$$= \sum_{n=1}^{\infty}\left|s_n|\bar{t}_n|^{p-2}\right|^q$$

$$= \sum_{n=1}^{\infty}|s_n|^q|\bar{t}_n|^{(p-2)q}$$

$$\leq \left(\sum_{n=1}^{\infty}(|s_n|^q)^{\frac{p}{q}}\right)^{\frac{q}{p}}\left(\sum_{n=1}^{\infty}\left(|\bar{t}_n|^{(p-2)q}\right)^{\frac{p}{p-q}}\right)^{\frac{p-q}{p}} \quad (2 < p < \infty \Longrightarrow 1 < q < 2)$$

$$\leq \left(\sum_{n=1}^{\infty}|s_n|^p\right)^{\frac{q}{p}}\left(\sum_{n=1}^{\infty}|\bar{t}_n|^p\right)^{\frac{p-q}{p}} \qquad ((p-2)q = p-q)$$

$$= \|v\|_p^q\ \|\bar{x}\|_p^{p-q}.$$

By (I), we have that

$$\left\{\frac{|\bar{t}_n|^{p-1}\mathrm{sign}(\bar{t}_n)}{\left(\sum_{i=1}^{\infty}|\bar{t}_i|^p\right)^{\frac{p-2}{p}}}\right\}_{n=1}^{\infty} = -\frac{(p-2)\sum_{i=1}^{\infty}\lambda_i|\bar{t}_i|^p}{\sum_{i=1}^{\infty}|\bar{t}_i|^p}\left\{\frac{|\bar{t}_n|^{p-1}\mathrm{sign}(\bar{t}_n)}{\left(\sum_{i=1}^{\infty}|\bar{t}_i|^p\right)^{\frac{p-2}{p}}}\right\}_{n=1}^{\infty} = -\frac{(p-2)\sum_{i=1}^{\infty}\lambda_i|\bar{t}_i|^p}{\sum_{i=1}^{\infty}|\bar{t}_i|^p}J(\bar{x}) \in l_q.$$

By (I), we have that

$$\left\{\frac{|\bar{t}_n|^{p-1}\mathrm{sign}(\bar{t}_n)}{\left(\sum_{i=1}^{\infty}|\bar{t}_i|^p\right)^{\frac{p-2}{p}}}\left((p-1)\lambda_n\right)\right\}_{n=1}^{\infty} = (p-1)\left\{\frac{|\bar{t}_n|^{p-1}\mathrm{sign}(\bar{t}_n)\lambda_n}{\left(\sum_{i=1}^{\infty}|\bar{t}_i|^p\right)^{\frac{p-2}{p}}}\right\}_{n=1}^{\infty} \in l_q. \tag{3.7}$$

By (3.5) and (3.7), we obtain that

$$\left\{\frac{|\bar{t}_n|^{p-1}\mathrm{sign}(\bar{t}_n)}{\left(\sum_{i=1}^{\infty}|\bar{t}_i|^p\right)^{\frac{p-2}{p}}}\left((p-1)\lambda_n - \frac{(p-2)\sum_{i=1}^{\infty}\lambda_i|\bar{t}_i|^p}{\sum_{i=1}^{\infty}|\bar{t}_i|^p}\right)\right\}_{n=1}^{\infty} \in l_q. \tag{3.8}$$

Considering the above estimations and (3.8), it allows us to take the last limit (3.6).

$$J'(\bar{x})(v) = \lim_{t\to 0}\frac{\left\{\dots,\ \frac{|ts_n|^{p-1}\mathrm{sign}(ts_n)}{\left(\sum_{i\in\mathcal{N}(\bar{x})}(1+t\lambda_i)^p|\bar{t}_i|^p+|t|^p\sum_{j\notin\mathcal{N}(\bar{x})}|s_j|^p\right)^{\frac{p-2}{p}}},\ \dots,(f_m(t)-f_m(0)),\dots\right\}}{t}$$

$$= \lim_{t\to 0}\frac{\left\{\dots,\ \frac{|ts_n|^{p-1}\mathrm{sign}(ts_n)}{\left(\sum_{i\in\mathcal{N}(\bar{x})}|1+t\lambda_i|^p|\bar{t}_i|^p+|t|^p\sum_{j\notin\mathcal{N}(\bar{x})}|s_j|^p\right)^{\frac{p-2}{p}}},\ \dots,f_m'(t_m)t,\dots\right\}}{t}$$

$$= \lim_{t\to 0}\left\{\dots,\ \frac{|t|^{p-2}|s_n|^{p-1}\mathrm{sign}(ts_n)}{\left(\sum_{i\in\mathcal{N}(\bar{x})}|1+t\lambda_i|^p|\bar{t}_i|^p+|t|^p\sum_{j\notin\mathcal{N}(\bar{x})}|s_j|^p\right)^{\frac{p-2}{p}}},\ \dots,\ f'(t_m),\dots\right\}$$

$$= \lim_{t\to 0}\left\{\ldots,\ \frac{|t|^{p-2}|s_n|^{p-1}\operatorname{sign}(ts_n)}{\left(\sum_{i\in\mathcal{N}(\bar{x})}|1+t\lambda_i|^p|\bar{t}_i|^p+|t|^p\sum_{j\notin\mathcal{N}(\bar{x})}|s_j|^p\right)^{\frac{p-2}{p}}},\ldots\right.$$

$$\left.\ldots,|\bar{t}_m|^{p-1}\operatorname{sign}(\bar{t}_m)\frac{(p-1)\lambda_m\sum_{i\in\mathcal{N}(\bar{x})}|\bar{t}_i|^p-(p-2)\sum_{i\in\mathcal{N}(\bar{x})}\lambda_i|\bar{t}_i|^p}{\left(\sum_{i\in\mathcal{N}(\bar{x})}|\bar{t}_i|^p\right)^{\frac{p-2}{p}+1}},\ldots\right\}$$

$$= \left\{\ldots,\ 0,\ \ \ldots,|\bar{t}_m|^{p-1}\operatorname{sign}(\bar{t}_m)\frac{(p-1)\lambda_m\sum_{i\in\mathcal{N}(\bar{x})}|\bar{t}_i|^p-(p-2)\sum_{i\in\mathcal{N}(\bar{x})}\lambda_i|\bar{t}_i|^p}{\left(\sum_{i\in\mathcal{N}(\bar{x})}|\bar{t}_i|^p\right)^{\frac{p-2}{p}+1}},\ldots\right\}$$

$$= \left\{\ldots,\ 0,\ \ \ldots,\frac{|\bar{t}_m|^{p-1}\operatorname{sign}(\bar{t}_m)}{\left(\sum_{i\in\mathcal{N}(\bar{x})}|\bar{t}_i|^p\right)^{\frac{p-2}{p}}}\frac{(p-1)\lambda_m\sum_{i\in\mathcal{N}(\bar{x})}|\bar{t}_i|^p-(p-2)\sum_{i\in\mathcal{N}(\bar{x})}\lambda_i|\bar{t}_i|^p}{\sum_{i\in\mathcal{N}(\bar{x})}|\bar{t}_i|^p},\ldots\right\}$$

$$= \left\{\ldots,\ 0,\ \ \ldots,\frac{|\bar{t}_m|^{p-1}\operatorname{sign}(\bar{t}_m)}{\left(\sum_{i\in\mathcal{N}(\bar{x})}|\bar{t}_i|^p\right)^{\frac{p-2}{p}}}\left((p-1)\lambda_m-\frac{(p-2)\sum_{i\in\mathcal{N}(\bar{x})}\lambda_i|\bar{t}_i|^p}{\sum_{i\in\mathcal{N}(\bar{x})}|\bar{t}_i|^p}\right),\ldots\right\}$$

$$= \left\{\ldots,\ 0,\ \ \ldots,\frac{|\bar{t}_m|^{p-1}\operatorname{sign}(\bar{t}_m)}{\left(\sum_{i=1}^{\infty}|\bar{t}_i|^p\right)^{\frac{p-2}{p}}}\left((p-1)\lambda_m-\frac{(p-2)\sum_{i=1}^{\infty}\lambda_i|\bar{t}_i|^p}{\sum_{i=1}^{\infty}|\bar{t}_i|^p}\right),\ldots\right\}\quad (\text{If } \bar{t}_i=0, \text{ then } \lambda_i=0)$$

$$= \left\{\frac{|\bar{t}_n|^{p-1}\operatorname{sign}(\bar{t}_n)}{\left(\sum_{i=1}^{\infty}|\bar{t}_i|^p\right)^{\frac{p-2}{p}}}\right\}_{n=1}^{\infty}\begin{pmatrix}(p-1)\lambda_1-\frac{(p-2)\sum_{i=1}^{\infty}\lambda_i|\bar{t}_i|^p}{\sum_{i=1}^{\infty}|\bar{t}_i|^p} & 0 & \ldots\\ 0 & (p-1)\lambda_2-\frac{(p-2)\sum_{i=1}^{\infty}\lambda_i|\bar{t}_i|^p}{\sum_{i=1}^{\infty}|\bar{t}_i|^p} & 0\\ \vdots & 0 & \ddots\end{pmatrix}$$

$$= J(\bar{x})\begin{pmatrix}(p-1)\lambda_1-\frac{(p-2)\sum_{i=1}^{\infty}\lambda_i|\bar{t}_i|^p}{\sum_{i=1}^{\infty}|\bar{t}_i|^p} & 0 & \ldots\\ 0 & (p-1)\lambda_2-\frac{(p-2)\sum_{i=1}^{\infty}\lambda_i|\bar{t}_i|^p}{\sum_{i=1}^{\infty}|\bar{t}_i|^p} & 0\\ \vdots & 0 & \ddots\end{pmatrix}$$

$$= J(\bar{x})\left((p-1)\begin{pmatrix}\lambda_1 & 0 & \ldots\\ 0 & \lambda_2 & 0\\ \vdots & 0 & \ddots\end{pmatrix}-\frac{(p-2)\sum_{i=1}^{\infty}\lambda_i|\bar{t}_i|^p}{\sum_{i=1}^{\infty}|\bar{t}_i|^p}\begin{pmatrix}1 & 0 & \ldots\\ 0 & 1 & 0\\ \vdots & 0 & \ddots\end{pmatrix}\right).$$

Proof of (a) in (ii). In particular, if $\bar{x}$ satisfies $|\mathcal{N}(\bar{x})| < \infty$, then, for any $v = \{s_n\}_{n=1}^{\infty} \in l_p\backslash\{\theta\}$, with respect to $\bar{x}$ with ratio sequence $\{\lambda_n\}_{n=1}^{\infty}$ that has finite number of nonzero entries. Then, *J* is Gâteaux differentiable at $\bar{x}$ along direction *v*. Let $2 < p < \infty$. For any $\bar{x} \in l_p\backslash\{\theta\}$, it is clear that the ratio sequence $\{\lambda_n\}_{n=1}^{\infty}$ of $\bar{x}$ with respect to itself satisfies that the entries are either 1 or 0:

$$\lambda_n = \begin{cases}0, & \text{if } \bar{t}_n = 0,\\ 1, & \text{if } \bar{t}_n \neq 0,\end{cases}\quad \text{for each } n = 1, 2, \ldots.$$

By the solution in (ii), we obtain

$$J'(\bar{x})(\bar{x}) = J(\bar{x}).$$

Then, (a) is proved by (i) and (ii).

Proof of (b) in (ii). Let $\bar{x} = \{\bar{t}_n\}_{n=1}^{\infty} \in l_p\backslash\{\theta\}$ any $v = \{s_n\}_{n=1}^{\infty} \in l_p\backslash\{\theta\}$. If, $s_n\bar{t}_n = 0$, for each $n$ = 1, 2, …, then ratio sequence $\{\lambda_n\}_{n=1}^{\infty}$ of $v$ respect to $\bar{x}$ satisfies $\{\lambda_n\}_{n=1}^{\infty} = \theta$, which implies

$$J'(\bar{x}, v) = J(\bar{x})\left((p-1)\begin{pmatrix} 0 & 0 & \dots \\ 0 & 0 & 0 \\ \vdots & 0 & \ddots \end{pmatrix} - \frac{(p-2)0}{\sum_{i=1}^{\infty}|\bar{t}_i|^p}\begin{pmatrix} 1 & 0 & \dots \\ 0 & 1 & 0 \\ \vdots & 0 & \ddots \end{pmatrix}\right) = \theta.$$

This completes the proof of part (ii). Next, we prove part (iii) for$1 < p < 2$ and $2 < q < \infty$.

Proof of (a) in part (iii). Suppose that $\mathcal{N}(v) \nsubseteq \mathcal{N}(\bar{x})$. Then there is a positive integer $m$ such that $\bar{t}_m = 0$ and $s_m \neq 0$. Then, for any $t \neq 0$, the $m$th entry of $\frac{J(\bar{x}+tv)-J(\bar{x})}{t}$ satisfies the following properties.

$$\frac{\frac{|\bar{t}_m+ts_m|^{p-1}\operatorname{sign}(\bar{t}_m+ts_m)}{\left(\sum_{i=1}^{\infty}|\bar{t}_i+ts_i|^p\right)^{\frac{p-2}{p}}} - \frac{|\bar{t}_m|^{p-1}\operatorname{sign}(\bar{t}_m)}{\left(\sum_{i=1}^{\infty}|\bar{t}_i|^p\right)^{\frac{p-2}{p}}}}{t} = \frac{\frac{|ts_m|^{p-1}\operatorname{sign}(ts_m)}{\left(\sum_{i=1}^{\infty}|\bar{t}_i+ts_i|^p\right)^{\frac{p-2}{p}}}}{t}$$

$$= \frac{\frac{|ts_m|^{p-1}\operatorname{sign}(s_m)\operatorname{sign}(t)}{\left(\sum_{i=1}^{\infty}|\bar{t}_i+ts_i|^p\right)^{\frac{p-2}{p}}}}{t} = \frac{\frac{|ts_m|^{p-1}\operatorname{sign}(s_m)}{\left(\sum_{i=1}^{\infty}|\bar{t}_i+ts_i|^p\right)^{\frac{p-2}{p}}}}{t\operatorname{sign}(t)}$$

$$= \frac{\frac{|t|^{p-1}|s_m|^{p-1}\operatorname{sign}(s_m)}{\left(\sum_{i=1}^{\infty}|\bar{t}_i+ts_i|^p\right)^{\frac{p-2}{p}}}}{|t|} = \frac{|t|^{p-2}|s_m|^{p-1}\operatorname{sign}(s_m)}{\left(\sum_{i=1}^{\infty}|\bar{t}_i+ts_i|^p\right)^{\frac{p-2}{p}}}.$$

By $p - 2 < 0$, this implies that

$$\lim_{t\to 0}\frac{\frac{|\bar{t}_m+ts_m|^{p-1}\operatorname{sign}(\bar{t}_m+ts_m)}{\left(\sum_{i=1}^{\infty}|\bar{t}_i+ts_i|^p\right)^{\frac{p-2}{p}}} - \frac{|\bar{t}_m|^{p-1}\operatorname{sign}(\bar{t}_m)}{\left(\sum_{i=1}^{\infty}|\bar{t}_i|^p\right)^{\frac{p-2}{p}}}}{t} = \lim_{t\to 0}\frac{|t|^{p-2}|s_m|^{p-1}\operatorname{sign}(s_m)}{\left(\sum_{i=1}^{\infty}|\bar{t}_i+ts_i|^p\right)^{\frac{p-2}{p}}} = \infty.$$

We obtain that if there is a positive integer $m$ such that $\bar{t}_m = 0$ and $s_m \neq 0$, then $J$ is not Gâteaux directionally differentiable at $\bar{x}$ along direction $v = \{s_n\}_{n=1}^{\infty} \in l_p\backslash\{\theta\}$. Hence,

$$\mathcal{N}(v) \nsubseteq \mathcal{N}(\bar{x}) \quad \Longrightarrow \quad J'(\bar{x}, v) \text{ does not exist.}$$

Proof (b) of part (iii). *Let* $1 < p < 2$. Let $\bar{x} = (\bar{t}_1, \bar{t}_2, \dots) \in l_p\backslash\{\theta\}$. Let $v = \{s_n\}_{n=1}^{\infty} \in l_p\backslash\{\theta\}$.

Suppose that $\mathcal{N}(v) \subseteq \mathcal{N}(\bar{x})$ and $v$ is ratio bounded and with ratio bounds $\Lambda$ with respect to $\bar{x}$, for some positive number $\Lambda$ satisfying

$$|\lambda_n| \leq \Lambda, \text{ for } n = 1, 2, \dots . \tag{3.9}$$

By the condition that $\mathcal{N}(v) \subseteq \mathcal{N}(\bar{x})$, for each positive integer $m$, there are 3 possible cases of the $m$th entry of $\frac{J(\bar{x}+tv)-J(\bar{x})}{t}$.

Case 1. For a positive integer $m$, suppose $\bar{t}_m = 0$ and $s_m = 0$. In this case, $\lambda_m = 0$. And for any $t \neq 0$,

the $m$th entry of $\frac{J(\bar{x}+tv)-J(\bar{x})}{t}$ has the following properties.

$$\frac{\frac{|\bar{t}_m+ts_m|^{p-1}\text{sign}(\bar{t}_m+ts_m)}{\left(\sum_{i=1}^{\infty}|\bar{t}_i+ts_i|^p\right)^{\frac{p-2}{p}}}-\frac{|\bar{t}_m|^{p-1}\text{sign}(\bar{t}_m)}{\left(\sum_{i=1}^{\infty}|\bar{t}_i|^p\right)^{\frac{p-2}{p}}}}{t}=\frac{0-0}{t}=0.$$

Case 2. For a positive integer $m$, suppose $\bar{t}_m \neq 0$ and $s_m = 0$. Then, for any $t \neq 0$, the $m$th entry of $\frac{J(\bar{x}+tv)-J(\bar{x})}{t}$ has the following properties.

$$\frac{\frac{|\bar{t}_m+ts_m|^{p-1}\text{sign}(\bar{t}_m+ts_m)}{\left(\sum_{i=1}^{\infty}|\bar{t}_i+ts_i|^p\right)^{\frac{p-2}{p}}}-\frac{|\bar{t}_m|^{p-1}\text{sign}(\bar{t}_m)}{\left(\sum_{i=1}^{\infty}|\bar{t}_i|^p\right)^{\frac{p-2}{p}}}}{t}$$

$$=\frac{\frac{|\bar{t}_m|^{p-1}\text{sign}(\bar{t}_m)}{\left(\sum_{i=1}^{\infty}|\bar{t}_i+ts_i|^p\right)^{\frac{p-2}{p}}}-\frac{|\bar{t}_m|^{p-1}\text{sign}(\bar{t}_m)}{\left(\sum_{i=1}^{\infty}|\bar{t}_i|^p\right)^{\frac{p-2}{p}}}}{t}$$

$$=|\bar{t}_m|^{p-1}\text{sign}(\bar{t}_m)\left(\frac{1}{\left(\sum_{i=1}^{\infty}|\bar{t}_i+ts_i|^p\right)^{\frac{p-2}{p}}}-\frac{1}{\left(\sum_{i=1}^{\infty}|\bar{t}_i|^p\right)^{\frac{p-2}{p}}}\right)$$

$$=|\bar{t}_m|^{p-1}\text{sign}(\bar{t}_m)\left(\left(\sum_{i=1}^{\infty}|\bar{t}_i+ts_i|^p\right)^{\frac{2-p}{p}}-\left(\sum_{i=1}^{\infty}|\bar{t}_i|^p\right)^{\frac{2-p}{p}}\right)\quad (1<p<2).$$

Define a function $h$ of $t$ by

$$h(t)=\left(\sum_{i=1}^{\infty}|\bar{t}_i+ts_i|^p\right)^{\frac{2-p}{p}},\ \text{for } t\in\mathbb{R}.$$

Then,

$$h'(t)=\frac{2-p}{p}\left(\sum_{i=1}^{\infty}|\bar{t}_i+ts_i|^p\right)^{\frac{2-p}{p}-1}\left(\sum_{i=1}^{\infty}p|\bar{t}_i+ts_i|^{p-1}s_i\text{sign}(\bar{t}_i+ts_i)\right)$$

$$=(2-p)\left(\sum_{i=1}^{\infty}|\bar{t}_i+ts_i|^p\right)^{\frac{2-p}{p}-1}\left(\sum_{i=1}^{\infty}|\bar{t}_i+ts_i|^{p-1}s_i\text{sign}(\bar{t}_i+ts_i)\right).$$

If $t$ satisfies that $|t|\leq 1$, we have

$$\left|\sum_{i=1}^{\infty}|\bar{t}_i+ts_i|^{p-1}s_i\text{sign}(\bar{t}_i+ts_i)\right|$$

$$\leq\left(\sum_{i=1}^{\infty}(|\bar{t}_i+ts_i|^{p-1})^q\right)^{\frac{1}{q}}\left(\sum_{i=1}^{\infty}|s_i\text{sign}(\bar{t}_i+ts_i)|^p\right)^{\frac{1}{p}}$$

$$=\left(\sum_{i=1}^{\infty}|\bar{t}_i+ts_i|^p\right)^{\frac{p}{pq}}\left(\sum_{i=1}^{\infty}|s_i|^p\right)^{\frac{1}{p}}$$

$$=\|\bar{x}+tv\|_p^{\frac{p}{q}}\|v\|_p$$

$$\leq\left(\|\bar{x}\|_p+\|tv\|_p\right)^{\frac{p}{q}}\|v\|_p$$

$$\leq\left(\|\bar{x}\|_p+\|v\|_p\right)^{\frac{p}{q}}\|v\|_p. \quad (3.10)$$

(3.10) ensures that (3.8) is well-defined. By Mean Value Theorem for this function $h$ into (3.7), we have

$$\left(\sum_{i=1}^{\infty}|\bar{t}_i + ts_i|^p\right)^{\frac{2-p}{p}} - \left(\sum_{i=1}^{\infty}|\bar{t}_i|^p\right)^{\frac{2-p}{p}}$$

$$= h(t) - h(0)$$

$$= h'(s)t \text{ with } |s| \leq |t|. \tag{3.11}$$

Since $\bar{t}_m \neq 0$ and $s_m = 0$, it yields that $\lambda_m = 0$. By (3.8) and substituting (3.11) to (3.7), for the $m$th entry of $\frac{J(\bar{x}+tv)-J(\bar{x})}{t}$, we have

$$\lim_{t\to 0}\frac{\frac{|\bar{t}_m+ts_m|^{p-1}\operatorname{sign}(\bar{t}_m+ts_m)}{\left(\sum_{i=1}^{\infty}|\bar{t}_i+ts_i|^p\right)^{\frac{p-2}{p}}} - \frac{|\bar{t}_m|^{p-1}\operatorname{sign}(\bar{t}_m)}{\left(\sum_{i=1}^{\infty}|\bar{t}_i|^p\right)^{\frac{p-2}{p}}}}{t}$$

$$= \lim_{t\to 0}\frac{|\bar{t}_m|^{p-1}\operatorname{sign}(\bar{t}_m)\left((2-p)\left(\sum_{i=1}^{\infty}|\bar{t}_i+ss_i|^p\right)^{\frac{2-p}{p}-1}\left(\sum_{i=1}^{\infty}|\bar{t}_i+ss_i|^{p-1}s_i\operatorname{sign}(\bar{t}_i+ss_i)\right)t\right)}{t}$$

$$= \lim_{t\to 0}|\bar{t}_m|^{p-1}\operatorname{sign}(\bar{t}_m)\left((2-p)\left(\sum_{i=1}^{\infty}|\bar{t}_i + ss_i|^p\right)^{\frac{2-p}{p}-1}\left(\sum_{i=1}^{\infty}|\bar{t}_i + ss_i|^{p-1} s_i\operatorname{sign}(\bar{t}_i + ss_i)\right)\right)$$

$$= \frac{|\bar{t}_m|^{p-1}\operatorname{sign}(\bar{t}_m)}{\left(\sum_{i=1}^{\infty}|\bar{t}_i|^p\right)^{\frac{p-2}{p}}}\left(-\frac{(p-2)\left(\sum_{i=1}^{\infty}|\bar{t}_i|^{p-1}s_i\operatorname{sign}(\bar{t}_i)\right)}{\sum_{i=1}^{\infty}|\bar{t}_i|^p}\right)$$

$$= \frac{|\bar{t}_m|^{p-1}\operatorname{sign}(\bar{t}_m)}{\left(\sum_{i=1}^{\infty}|\bar{t}_i|^p\right)^{\frac{p-2}{p}}}\left((p-1)\lambda_m - \frac{(p-2)\left(\sum_{i=1}^{\infty}|\bar{t}_i|^{p-1}s_i\operatorname{sign}(\bar{t}_i)\right)}{\sum_{i=1}^{\infty}|\bar{t}_i|^p}\right).$$

$$= \frac{|\bar{t}_m|^{p-1}\operatorname{sign}(\bar{t}_m)}{\left(\sum_{i=1}^{\infty}|\bar{t}_i|^p\right)^{\frac{p-2}{p}}}\left((p-1)\lambda_m - \frac{(p-2)\left(\sum_{i=1}^{\infty}|\bar{t}_i|^p\frac{s_i}{|\bar{t}_i|}\operatorname{sign}(\bar{t}_i)\right)}{\sum_{i=1}^{\infty}|\bar{t}_i|^p}\right).$$

$$= \frac{|\bar{t}_m|^{p-1}\operatorname{sign}(\bar{t}_m)}{\left(\sum_{i=1}^{\infty}|\bar{t}_i|^p\right)^{\frac{p-2}{p}}}\left((p-1)\lambda_m - \frac{(p-2)\left(\sum_{i=1}^{\infty}|\bar{t}_i|^p\frac{s_i}{|\bar{t}_i|\operatorname{sign}(\bar{t}_i)}\right)}{\sum_{i=1}^{\infty}|\bar{t}_i|^p}\right).$$

$$= \frac{|\bar{t}_m|^{p-1}\operatorname{sign}(\bar{t}_m)}{\left(\sum_{i=1}^{\infty}|\bar{t}_i|^p\right)^{\frac{p-2}{p}}}\left((p-1)\lambda_m - \frac{(p-2)\left(\sum_{i=1}^{\infty}|\bar{t}_i|^p\frac{s_i}{\bar{t}_i}\right)}{\sum_{i=1}^{\infty}|\bar{t}_i|^p}\right).$$

$$= \frac{|\bar{t}_m|^{p-1}\operatorname{sign}(\bar{t}_m)}{\left(\sum_{i=1}^{\infty}|\bar{t}_i|^p\right)^{\frac{p-2}{p}}}\left((p-1)\lambda_m - \frac{(p-2)\left(\sum_{i=1}^{\infty}|\bar{t}_i|^p\lambda_i\right)}{\sum_{i=1}^{\infty}|\bar{t}_i|^p}\right).$$

Case 3. For some positive integer $m$, suppose $\bar{t}_m \neq 0$ and $s_m \neq 0$. Then, for any $t \neq 0$, the $m$th entry of $\frac{J(\bar{x}+tv)-J(\bar{x})}{t}$ has the following properties.

$$\frac{\frac{|\bar{t}_m+ts_m|^{p-1}\operatorname{sign}(\bar{t}_m+ts_m)}{\left(\sum_{i=1}^{\infty}|\bar{t}_i+ts_i|^p\right)^{\frac{p-2}{p}}} - \frac{|\bar{t}_m|^{p-1}\operatorname{sign}(\bar{t}_m)}{\left(\sum_{i=1}^{\infty}|\bar{t}_i|^p\right)^{\frac{p-2}{p}}}}{t}$$

$$= \frac{|\bar{t}_m+ts_m|^{p-1}\text{sign}(\bar{t}_m+ts_m)\left(\sum_{i=1}^{\infty}|\bar{t}_i+ts_i|^p\right)^{\frac{2-p}{p}} - |\bar{t}_m|^{p-1}\text{sign}(\bar{t}_m)\left(\sum_{i=1}^{\infty}|\bar{t}_i|^p\right)^{\frac{2-p}{p}}}{t} \quad (1 < p < 2). \quad (3.12)$$

Define a function $k$ of $t$ by

$$k(t) = |\bar{t}_m + ts_m|^{p-1}\text{sign}(\bar{t}_m + ts_m)\left(\sum_{i=1}^{\infty}|\bar{t}_i + ts_i|^p\right)^{\frac{2-p}{p}}, \text{ for } t \in \mathbb{R}.$$

Then,

$$k'(t) = (p-1)|\bar{t}_m + ts_m|^{p-2}s_m\text{sign}(\bar{t}_m + ts_m)\text{sign}(\bar{t}_m + ts_m)\left(\sum_{i=1}^{\infty}|\bar{t}_i + ts_i|^p\right)^{\frac{2-p}{p}}$$

$$+|\bar{t}_m + ts_m|^{p-1}\text{sign}(\bar{t}_m + ts_m)\frac{2-p}{p}\left(\sum_{i=1}^{\infty}|\bar{t}_i + ts_i|^p\right)^{\frac{2-p}{p}-1}\left(\sum_{i=1}^{\infty}p|\bar{t}_i + ts_i|^{p-1}s_i\text{sign}(\bar{t}_i + ts_i)\right)$$

$$= (p-1)|\bar{t}_m + ts_m|^{p-2}s_m\left(\sum_{i=1}^{\infty}|\bar{t}_i + ts_i|^p\right)^{\frac{2-p}{p}}$$

$$+(2-p)|\bar{t}_m + ts_m|^{p-1}\text{sign}(\bar{t}_m + ts_m)\left(\sum_{i=1}^{\infty}|\bar{t}_i + ts_i|^p\right)^{\frac{2-p}{p}-1}\left(\sum_{i=1}^{\infty}|\bar{t}_i + ts_i|^{p-1}s_i\text{sign}(\bar{t}_i + ts_i)\right). \quad (3.13)$$

Similarly to (3.8) and (3.10), one can show that (3.13) is well-defined. By Mean Value Theorem for this function $k$ into (3.12), we have

$$|\bar{t}_m + ts_m|^{p-1}\text{sign}(\bar{t}_m + ts_m)\left(\sum_{i=1}^{\infty}|\bar{t}_i + ts_i|^p\right)^{\frac{2-p}{p}} - |\bar{t}_m|^{p-1}\text{sign}(\bar{t}_m)\left(\sum_{i=1}^{\infty}|\bar{t}_i|^p\right)^{\frac{2-p}{p}}$$

$$= k(t) - k(0)$$

$$= k'(s)t \text{ with } |s| \le |t|. \quad (3.14)$$

Since $\bar{t}_m \neq 0$ and $s_m \neq 0$, it yields that $\lambda_m \neq 0$. Substituting (3.14) to (3.12), for the $m$th entry of $\frac{J(\bar{x}+tv)-J(\bar{x})}{t}$, we have

$$\lim_{t\to 0}\frac{\frac{|\bar{t}_m+ts_m|^{p-1}\text{sign}(\bar{t}_m+ts_m)}{\left(\sum_{i=1}^{\infty}|\bar{t}_i+ts_i|^p\right)^{\frac{p-2}{p}}} - \frac{|\bar{t}_m|^{p-1}\text{sign}(\bar{t}_m)}{\left(\sum_{i=1}^{\infty}|\bar{t}_i|^p\right)^{\frac{p-2}{p}}}}{t}$$

$$= \lim_{t\to 0}\left(\frac{(p-1)|\bar{t}_m+ss_m|^{p-2}s_m\left(\sum_{i=1}^{\infty}|\bar{t}_i+ss_i|^p\right)^{\frac{2-p}{p}}}{t}\right.$$

$$\left.+\frac{(2-p)|\bar{t}_m+ss_m|^{p-1}\text{sign}(\bar{t}_m+ss_m)\left(\sum_{i=1}^{\infty}|\bar{t}_i+ss_i|^p\right)^{\frac{2-p}{p}-1}\left(\sum_{i=1}^{\infty}|\bar{t}_i+ss_i|^{p-1}s_i\text{sign}(\bar{t}_i+ss_i)\right)}{t}\right)t$$

$$= \lim_{t\to 0}\left((p-1)|\bar{t}_m + ss_m|^{p-2}s_m\left(\sum_{i=1}^{\infty}|\bar{t}_i + ss_i|^p\right)^{\frac{2-p}{p}}\right.$$

$$\left.+ (2-p)|\bar{t}_m + ss_m|^{p-1}\text{sign}(\bar{t}_m + ss_m)\left(\sum_{i=1}^{\infty}|\bar{t}_i + ss_i|^p\right)^{\frac{2-p}{p}-1}\left(\sum_{i=1}^{\infty}|\bar{t}_i + ss_i|^{p-1}s_i\text{sign}(\bar{t}_i + ss_i)\right)\right)$$

$$= (p-1)|\bar{t}_m|^{p-2}s_m\left(\sum_{i=1}^{\infty}|\bar{t}_i|^p\right)^{\frac{2-p}{p}}$$

$$+(2-p)|\bar{t}_m|^{p-1}\text{sign}(\bar{t}_m)(\sum_{i=1}^{\infty}|\bar{t}_i|^p)^{\frac{2-p}{p}-1}(\sum_{i=1}^{\infty}|\bar{t}_i|^{p-1}\, s_i\text{sign}(\bar{t}_i))$$

$$= (p-1)|\bar{t}_m|^{p-1}\frac{s_m}{|\bar{t}_m|}\text{sign}(\bar{t}_m)\text{sign}(\bar{t}_m)(\sum_{i=1}^{\infty}|\bar{t}_i|^p)^{\frac{2-p}{p}}$$

$$+(2-p)|\bar{t}_m|^{p-1}\text{sign}(\bar{t}_m)(\sum_{i=1}^{\infty}|\bar{t}_i|^p)^{\frac{2-p}{p}-1}(\sum_{i=1}^{\infty}|\bar{t}_i|^{p-1}\, s_i\text{sign}(\bar{t}_i))$$

$$= (p-1)|\bar{t}_m|^{p-1}\frac{s_m}{|\bar{t}_m|\text{sign}(\bar{t}_m)}\text{sign}(\bar{t}_m)(\sum_{i=1}^{\infty}|\bar{t}_i|^p)^{\frac{2-p}{p}}$$

$$+(2-p)|\bar{t}_m|^{p-1}\text{sign}(\bar{t}_m)(\sum_{i=1}^{\infty}|\bar{t}_i|^p)^{\frac{2-p}{p}}(\sum_{i=1}^{\infty}|\bar{t}_i|^p)^{-1}(\sum_{i=1}^{\infty}|\bar{t}_i|^{p-1}\, s_i\text{sign}(\bar{t}_i))$$

$$= (p-1)|\bar{t}_m|^{p-1}\frac{s_m}{\bar{t}_m}\text{sign}(\bar{t}_m)(\sum_{i=1}^{\infty}|\bar{t}_i|^p)^{\frac{2-p}{p}}$$

$$+(2-p)|\bar{t}_m|^{p-1}\text{sign}(\bar{t}_m)(\sum_{i=1}^{\infty}|\bar{t}_i|^p)^{\frac{2-p}{p}}(\sum_{i=1}^{\infty}|\bar{t}_i|^p)^{-1}(\sum_{i=1}^{\infty}|\bar{t}_i|^{p-1}\, s_i\text{sign}(\bar{t}_i))$$

$$= \frac{|\bar{t}_m|^{p-1}\text{sign}(\bar{t}_m)}{(\sum_{i=1}^{\infty}|\bar{t}_i|^p)^{\frac{p-2}{p}}}\left((p-1)\lambda_m - \frac{(p-2)\left(\sum_{i=1}^{\infty}|\bar{t}_i|^{p-1}s_i\text{sign}(\bar{t}_i)\right)}{\sum_{i=1}^{\infty}|\bar{t}_i|^p}\right).$$

$$= \frac{|\bar{t}_m|^{p-1}\text{sign}(\bar{t}_m)}{(\sum_{i=1}^{\infty}|\bar{t}_i|^p)^{\frac{p-2}{p}}}\left((p-1)\lambda_m - \frac{(p-2)\left(\sum_{i=1}^{\infty}|\bar{t}_i|^p\lambda_i\right)}{\sum_{i=1}^{\infty}|\bar{t}_i|^p}\right).$$

For the case in part (iii) that $1 < p < 2$, similarly to part (ii) that is the case $2 < p < \infty$, in order to ensure the following statement

$$\left\{\frac{|\bar{t}_n|^{p-1}\text{sign}(\bar{t}_n)}{(\sum_{i=1}^{\infty}|\bar{t}_i|^p)^{\frac{p-2}{p}}}\left((p-1)\lambda_m - \frac{(p-2)\left(\sum_{i=1}^{\infty}|\bar{t}_i|^p\lambda_i\right)}{\sum_{i=1}^{\infty}|\bar{t}_i|^p}\right)\right\}_{n=1}^{\infty} \in l_q.$$

For any $\bar{x} = (\bar{t}_1, \bar{t}_2, \dots) \in l_p\backslash\{\theta\}$ and $v = \{s_n\}_{n=1}^{\infty} \in l_p\backslash\{\theta\}$ with ratio sequence $\{\lambda_n\}_{n=1}^{\infty}$ with respect to $\bar{x}$, we need to show that $\bar{x}$ and $v$ satisfy the following two properties (III) and (IV). By the assumption in (3.9), $v$ is ratio bounded and with ratio bounds $\Lambda$ with respect to $\bar{x}$, for some positive number $\Lambda$ satisfying

$$|\lambda_n| \le \Lambda, \text{ for } n = 1, 2, \dots .$$

(III) $\quad |\sum_{i=1}^{\infty}\lambda_i|\bar{t}_i|^p| < \infty$. We estimate

$$|\sum_{i=1}^{\infty}\lambda_i|\bar{t}_i|^p| \le \Lambda\sum_{i=1}^{\infty}|\bar{t}_i|^p < \infty.$$

(IV) $\quad \left\{\frac{|\bar{t}_n|^{p-1}\text{sign}(\bar{t}_n)(p-1)\lambda_n}{(\sum_{i=1}^{\infty}|\bar{t}_i|^p)^{\frac{p-2}{p}}}\right\}_{n=1}^{\infty} \in l_q.$

By (3.9), we only estimate its numerator as follows without the constant $p-1$.

$$\sum_{n=1}^{\infty}||\bar{t}_n|^{p-1}\text{sign}(\bar{t}_n)\lambda_n|^q \le \Lambda\sum_{n=1}^{\infty}|\bar{t}_n|^{(p-1)q}$$

$$= \Lambda\sum_{n=1}^{\infty}|\bar{t}_n|^p = \Lambda\|\bar{x}\|_p^p < \infty.$$

By (III), we have that

$$\left\{\frac{|\bar{t}_n|^{p-1}\text{sign}(\bar{t}_n)}{\left(\sum_{i=1}^{\infty}|\bar{t}_i|^p\right)^{\frac{p-2}{p}}}\right\}_{n=1}^{\infty} = -\frac{(p-2)\sum_{i=1}^{\infty}\lambda_i|\bar{t}_i|^p}{\sum_{i=1}^{\infty}|\bar{t}_i|^p}\left\{\frac{|\bar{t}_n|^{p-1}\text{sign}(\bar{t}_n)}{\left(\sum_{i=1}^{\infty}|\bar{t}_i|^p\right)^{\frac{p-2}{p}}}\right\}_{n=1}^{\infty} = -\frac{(p-2)\sum_{i=1}^{\infty}\lambda_i|\bar{t}_i|^p}{\sum_{i=1}^{\infty}|\bar{t}_i|^p}J(\bar{x}) \in l_q.$$

By (IV), we have that

$$\left\{\frac{|\bar{t}_n|^{p-1}\text{sign}(\bar{t}_n)}{\left(\sum_{i=1}^{\infty}|\bar{t}_i|^p\right)^{\frac{p-2}{p}}}\left((p-1)\lambda_n\right)\right\}_{n=1}^{\infty} = (p-1)\left\{\frac{|\bar{t}_n|^{p-1}\text{sign}(\bar{t}_n)\lambda_n}{\left(\sum_{i=1}^{\infty}|\bar{t}_i|^p\right)^{\frac{p-2}{p}}}\right\}_{n=1}^{\infty} \in l_q.$$

By (3.5) and (3.7), under the conditions that $\mathcal{N}(v) \subseteq \mathcal{N}(\bar{x})$ and $v$ is ratio bounded with respect to $\bar{x}$, we obtain that

$$\left\{\frac{|\bar{t}_n|^{p-1}\text{sign}(\bar{t}_n)}{\left(\sum_{i=1}^{\infty}|\bar{t}_i|^p\right)^{\frac{p-2}{p}}}\left((p-1)\lambda_n - \frac{(p-2)\sum_{i=1}^{\infty}\lambda_i|\bar{t}_i|^p}{\sum_{i=1}^{\infty}|\bar{t}_i|^p}\right)\right\}_{n=1}^{\infty} \in l_q.$$

The above inclusion and cases 1, 2 and 3 imply that $J$ is Gâteaux directionally differentiable at $\bar{x}$ along direction $v$ and

$$J'(\bar{x}, v) = J(\bar{x})\left((p-1)\begin{pmatrix}\lambda_1 & 0 & \dots \\ 0 & \lambda_2 & 0 \\ \vdots & 0 & \ddots\end{pmatrix} - \frac{(p-2)\sum_{i=1}^{\infty}\lambda_i|\bar{t}_i|^p}{\sum_{i=1}^{\infty}|\bar{t}_i|^p}\begin{pmatrix}1 & 0 & \dots \\ 0 & 1 & 0 \\ \vdots & 0 & \ddots\end{pmatrix}\right).$$

Then, we show the special cases of part (b) in (ii).

Proof of ($b_1$) of part (b) in (ii). By $|\mathcal{N}(v)| < \infty$, the ratio sequence $\{\lambda_n\}_{n=1}^{\infty}$ of $v$ with respect to $\bar{x}$ has finite number of nonzero entries that immediately imply that $v$ is ratio bounded with respect to $\bar{x}$. Then, ($b_1$) is induced by part (b) in (ii).

Proof of ($b_2$) of part (b) in (ii). Recall that the ratio sequence $\{\lambda_n\}_{n=1}^{\infty}$ of $v$ with respect to $\bar{x}$ is defined by

$$\lambda_n = \begin{cases} 0, & \text{if } \bar{t}_n = 0, \\ \frac{s_n}{\bar{t}_n}, & \text{if } \bar{t}_n \neq 0, \end{cases} \quad \text{for } n = 1, 2, \dots.$$

Then the conditions that $\mathcal{N}(v) \subseteq \mathcal{N}(\bar{x})$ and $0 < |\mathcal{N}(\bar{x})| < \infty$ imply $\mathcal{N}(v) \subseteq \mathcal{N}(\bar{x})$ and $|\mathcal{N}(v)| < \infty$. So, ($b_2$) is a special case of ($b_1$).

Proof of part (iv). Let $1 < p, q < \infty$ and $\frac{1}{p} + \frac{1}{q} = 1$. Part (iv) can be implied by parts (ii) and (iii). However, we provide a direct proof of part (iv). Let $\bar{x} = (\bar{t}_1, \bar{t}_2, \dots) \in l_p \backslash \{\theta\}$. For real number $t$ with $|t| < 1$, we calculate

$$J'(\bar{x})(\bar{x}) = \lim_{t\to 0}\frac{J(\bar{x}+t\bar{x}) - J(\bar{x})}{t}$$

$$= \lim_{t\to 0}\frac{\left\{\frac{|\bar{t}_n+t\bar{t}_n|^{p-1}\text{sign}(\bar{t}_n+t\bar{t}_n)}{\left(\sum_{i=1}^{\infty}|\bar{t}_i+t\bar{t}_i|^p\right)^{\frac{p-2}{p}}}\right\}_{n=1}^{\infty} - \left\{\frac{|\bar{t}_n|^{p-1}\text{sign}(\bar{t}_n)}{\left(\sum_{i=1}^{\infty}|\bar{t}_i|^p\right)^{\frac{p-2}{p}}}\right\}_{n=1}^{\infty}}{t}$$

$$= \lim_{t\to 0} \frac{\left\{ \frac{|\bar{t}_n + t\bar{t}_n|^{p-1}\text{sign}(\bar{t}_n + t\bar{t}_n)}{\left(\sum_{i=1}^{\infty}|\bar{t}_i + t\bar{t}_i|^p\right)^{\frac{p-2}{p}}} - \frac{|\bar{t}_n|^{p-1}\text{sign}(\bar{t}_n)}{\left(\sum_{i=1}^{\infty}|\bar{t}_i|^p\right)^{\frac{p-2}{p}}} \right\}_{n=1}^{\infty}}{t}$$

$$= \lim_{t\to 0} \frac{\left\{ \frac{|(1+t)\bar{t}_n|^{p-1}\text{sign}((1+t)\bar{t}_n)}{\left(\sum_{i=1}^{\infty}|(1+t)\bar{t}_i|^p\right)^{\frac{p-2}{p}}} - \frac{|\bar{t}_n|^{p-1}\text{sign}(\bar{t}_n)}{\left(\sum_{i=1}^{\infty}|\bar{t}_i|^p\right)^{\frac{p-2}{p}}} \right\}_{n=1}^{\infty}}{t}$$

$$= \lim_{\substack{t\to 0\\ |t|<1}} \frac{\left\{ \text{sign}(\bar{t}_n)\left( (1+t)\frac{|\bar{t}_n|^{p-1}}{\left(\sum_{i=1}^{\infty}|\bar{t}_i|^p\right)^{\frac{p-2}{p}}} - \frac{|\bar{t}_n|^{p-1}}{\left(\sum_{i=1}^{\infty}|\bar{t}_i|^p\right)^{\frac{p-2}{p}}} \right) \right\}_{n=1}^{\infty}}{t}$$

$$= \lim_{\substack{t\to 0\\ |t|<1}} \frac{\left\{ \text{sign}(\bar{t}_n)\left( \frac{t|\bar{t}_n|^{p-1}}{\left(\sum_{i=1}^{\infty}|\bar{t}_i|^p\right)^{\frac{p-2}{p}}} \right) \right\}_{n=1}^{\infty}}{t}$$

$$= \left\{ \text{sign}(\bar{t}_n)\left( \frac{|\bar{t}_n|^{p-1}}{\left(\sum_{i=1}^{\infty}|\bar{t}_i|^p\right)^{\frac{p-2}{p}}} \right) \right\}_{n=1}^{\infty}$$

$$= J(\bar{x}).$$

Hence, we obtain that

$$J'(\bar{x})(\bar{x}) = J(\bar{x}), \text{ for any } \bar{x} \in l_p\backslash\{\theta\}. \tag{3.15}$$

In particular, suppose $|\mathcal{N}(\bar{x})| = 1$ that there is a unique nonzero entry $\bar{t}_m$ of $\bar{x} = (\bar{t}_1, \bar{t}_2, \dots) \in l_p\backslash\{\theta\}$ with $\bar{t}_m \neq 0$ and $\bar{t}_n = 0$, for all $n \neq m$. Then, one can check that

$$J(\bar{x}) = \bar{x}, \text{ for any } \bar{x} \in l_p \text{ with } |\mathcal{N}(\bar{x})| \leq 1.$$

By (3.15), this implies

$$J'(\bar{x})(\bar{x}) = \bar{x}, \text{ for any } \bar{x} \in l_p \text{ with } |\mathcal{N}(\bar{x})| = 1.$$

This completes the proof of this theorem. □

**Remarks 3.3.** Let $1 < p < 2$. Let $\bar{x} = (\bar{t}_1, \bar{t}_2, \dots) \in l_p\backslash\{\theta\}$. Let $v = \{s_n\}_{n=1}^{\infty} \in l_p\backslash\{\theta\}$ with ratio sequence $\{\lambda_n\}_{n=1}^{\infty}$ with respect to $\bar{x}$. By part (b) of (ii) in Theorem 3.2, the conditions that $\mathcal{N}(v) \subseteq \mathcal{N}(\bar{x})$ and $v$ is ratio bounded with respect to $\bar{x}$ ensure the Gâteaux directionally differentiability of $J$ at $\bar{x}$ along direction $v$ and $J'(\bar{x}, v)$ is explicitly represented by (3.3). However, we will provide a counterexample below to demonstrate that the ratio boundness of $v$ with respect to $\bar{x}$ is also necessarily for the Gâteaux directionally differentiability of $J$ at $\bar{x}$ in direction $v$ satisfying (3.3).

**Example 3.4**. Let $1 < p < 2$, $2 < q < \infty$ with $\frac{1}{p} + \frac{1}{q} = 1$. In this example, we construct $\bar{x}, v \in l_p\backslash\{\theta\}$ satisfying that $\mathcal{N}(v)=\mathcal{N}(\bar{x})$ and $|\mathcal{N}(\bar{x})|=\infty$ and the ratio sequence of $v$ with respect to $\bar{x}$ is unbounded. Then, we show that $J$ is not Gâteaux directionally differentiability at $\bar{x}$ in direction $v$ with that (3.3) is not satisfied.

The condition that $\frac{1}{p}+\frac{1}{q}=1$ implies $\frac{p}{(p-1)q}=1$. Let $\beta>1=\frac{p}{(p-1)q}$. Since $0<\frac{p}{q}<1$ and $\frac{p}{q}=p-1$, by $1<p<2$, we calculate

$$(\beta-1)-\left((p-1)\beta-\frac{p}{q}\right)=2\beta-1-p\beta+\frac{p}{q}$$

$$=2\beta-1-p\beta+p-1=(2-p)(\beta-1)>0.$$

And $(p-1)\beta-\frac{p}{q}>(p-1)\,\frac{p}{(p-1)q}-\frac{p}{q}=0$. Then, we take $\alpha$ with

$$\beta-1>\alpha\geq(p-1)\beta-\frac{p}{q}.$$

Now, we define $\bar{x}=\left\{\frac{1}{n^{\frac{\beta}{p}}}\right\}_{n=1}^{\infty}$ and $v=\left\{\frac{1}{n^{\frac{\beta-\alpha}{p}}}\right\}_{n=1}^{\infty}$. By the selection of $\alpha$, we have that $\bar{x}\in l_p\backslash\{\theta\}$ and $v\in l_p\backslash\{\theta\}$ and $\mathcal{N}(v)=\mathcal{N}(\bar{x})$ with $|\mathcal{N}(\bar{x})|=\infty$. The ratio sequence of $v$ with respect to $\bar{x}$ is $\left\{n^{\frac{\alpha}{p}}\right\}_{n=1}^{\infty}$ that is unbounded. We calculate

$$\left\{\frac{\left|\frac{1}{n^{\frac{\beta}{p}}}\right|^{p-1}\operatorname{sign}\left(\frac{1}{n^{\frac{\beta}{p}}}\right)}{\left(\sum_{i=1}^{\infty}\left|\frac{1}{n^{\frac{\beta}{p}}}\right|^{p}\right)^{\frac{p-2}{p}}}(p-1)n^{\frac{\alpha}{p}}\right\}_{n=1}^{\infty}=\frac{p-1}{\left(\sum_{i=1}^{\infty}\frac{1}{i^{\beta}}\right)^{\frac{p-2}{p}}}\left\{\frac{1}{n^{\frac{\beta}{p}(p-1)}}n^{\frac{\alpha}{p}}\right\}_{n=1}^{\infty}=\frac{p-1}{\left(\sum_{i=1}^{\infty}\frac{1}{i^{\beta}}\right)^{\frac{p-2}{p}}}\left\{\frac{1}{n^{\frac{(p-1)\beta-\alpha}{p}}}\right\}_{n=1}^{\infty}.$$

By $\alpha\geq(p-1)\beta-\frac{p}{q}$, we have that $\frac{(p-1)\beta-\alpha}{p}q\leq\frac{\frac{p}{q}}{p}q=1$. This implies that

$$\frac{p-1}{\left(\sum_{i=1}^{\infty}\frac{1}{i^{\beta}}\right)^{\frac{p-2}{p}}}\left\{\frac{1}{n^{\frac{(p-1)\beta-\alpha}{p}}}\right\}_{n=1}^{\infty}\notin l_q.$$

This implies that for $\bar{x}=\left\{\frac{1}{n^{\frac{\beta}{p}}}\right\}_{n=1}^{\infty}$ and $v=\left\{\frac{1}{n^{\frac{\beta-\alpha}{p}}}\right\}_{n=1}^{\infty}$,

$$J(\bar{x})\left((p-1)\begin{pmatrix}1&0&\cdots\\0&2^{\frac{\alpha}{p}}&0\\\vdots&0&\ddots\end{pmatrix}-\frac{(p-2)\sum_{i=1}^{\infty}i^{\frac{\alpha}{p}}\frac{1}{i^{\beta}}}{\sum_{i=1}^{\infty}\frac{1}{i^{\beta}}}\begin{pmatrix}1&0&\cdots\\0&1&0\\\vdots&0&\ddots\end{pmatrix}\right)\notin l_q.$$

Hence, for $\bar{x}=\left\{\frac{1}{n^{\frac{\beta}{p}}}\right\}_{n=1}^{\infty}$ and $v=\left\{\frac{1}{n^{\frac{\beta-\alpha}{p}}}\right\}_{n=1}^{\infty}$, (3.3) is not defined in $l_q$.

## 4. The Fréchet Differentiability of the Normalized Duality Mapping $J$ in $l_p$

To start this section for considering the Fréchet differentiability of $J$ in $l_p$, we prove the following proposition with the Fréchet differentiability of $J$ at $\theta$.

**Proposition 4.1**. *Let* $1 < p < \infty$. *$J$ is not Fréchet differentiable at $\theta$; that is,* $\nabla J(\theta)$ *does not exist.*

*Proof.* Assume, by the way of contradiction, that $J$ is Fréchet differentiable at $\theta$, which induces that $\nabla J(\theta)$ exists and $\nabla J(\theta)\colon l_p \to l_q$ is a continuous and linear operator. However, by the connection between Gâteaux differentiability and Fréchet differentiability in Banach spaces and by (3.2) in part (i) of Theorem 3.2, this assumption implies that

$$\nabla J(\theta)(v) = J'(\theta)(v) = J(v)\text{, for any } v \in l_p\backslash\{\theta\}.$$

It is clear that $J$ is not linear, which gets a contradiction and proves this proposition. □

Recall that in case if $2 < p < \infty$, for any given $\bar{x} = (\bar{t}_1, \bar{t}_2, \ldots) \in l_p\backslash\{\theta\}$, by (3.3) in part (ii) of Theorem 3.2, $J$ is Gâteaux differentiable at $\bar{x}$ and the Gâteaux derivative of $J$ at $\bar{x}$ satisfies

$$J'(\bar{x})(v) = J(\bar{x})\left((p-1)\begin{pmatrix}\lambda_1 & 0 & \ldots \\ 0 & \lambda_2 & 0 \\ \vdots & 0 & \ddots\end{pmatrix} - \frac{(p-2)\sum_{i=1}^{\infty}\lambda_i|\bar{t}_i|^p}{\sum_{i=1}^{\infty}|\bar{t}_i|^p}\begin{pmatrix}1 & 0 & \ldots \\ 0 & 1 & 0 \\ \vdots & 0 & \ddots\end{pmatrix}\right),$$

for any $v = \{s_n\}_{n=1}^{\infty} \in l_p\backslash\{\theta\}$ with ratio sequence $\{\lambda_n\}_{n=1}^{\infty}$ with respect to $\bar{x}$. More precisely, $J'(\bar{x})$ has the following properties.

**Lemma 4.2**. *Let* $2 < p < \infty$. *For any given* $\bar{x} = (\bar{t}_1, \bar{t}_2, \ldots) \in l_p\backslash\{\theta\}$, $J'(\bar{x})(\cdot)\colon l_p \to l_q$ *represented by* (3.3) *is a continuous and linear operator from $l_p$ to $l_q$*. Moreover, $J'(\bar{x})(\cdot)\colon l_p \to l_q$ *is a linear and Lipchitz operator such that*

$$\|J'(\bar{x})(v-u)\|_q \le \left(p-1+(p-2)\|\bar{x}\|_p^q\right)\|v-u\|_p, \text{ for any } u, v \in l_p\backslash\{\theta\} \text{ with } u \neq v. \tag{4.1}$$

*Proof.* Let $\bar{x} = (\bar{t}_1, \bar{t}_2, \ldots) \in l_p\backslash\{\theta\}$. By the property (3.1) of ratio sequences of points in $l_p$ with respect to $\bar{x}$, the operator $J'(\bar{x})(\cdot)\colon l_p \to l_q$ is linear. Then, let $v = \{s_n\}_{n=1}^{\infty} \in l_p\backslash\{\theta\}$ and $u = \{t_n\}_{n=1}^{\infty} \in l_p\backslash\{\theta\}$ with ratio sequences $\{\lambda_n\}_{n=1}^{\infty}$ and $\{\gamma_n\}_{n=1}^{\infty}$, respectively. By the linearity of $J'(\bar{x})(\cdot)\colon l_p \to l_q$, we have

$$J'(\bar{x})(v-u) = J(\bar{x})\left((p-1)\begin{pmatrix}\lambda_1-\gamma_1 & 0 & \ldots \\ 0 & \lambda_2-\gamma_2 & 0 \\ \vdots & 0 & \ddots\end{pmatrix} - \frac{(p-2)\sum_{i=1}^{\infty}(\lambda_i-\gamma_i)|\bar{t}_i|^p}{\sum_{i=1}^{\infty}|\bar{t}_i|^p}\begin{pmatrix}1 & 0 & \ldots \\ 0 & 1 & 0 \\ \vdots & 0 & \ddots\end{pmatrix}\right).$$

Similarly to the proof of (I) in the proof of Theorem 3.2, we first estimate the sum inside (4.2).

$$\begin{aligned} &\left|\sum_{i=1}^{\infty}(\lambda_i-\gamma_i)|\bar{t}_i|^p\right| \\ &\le \sum_{i=1}^{\infty}\left|\lambda_i-\gamma_i|\bar{t}_i|^p\right| \\ &= \sum_{i=1}^{\infty}|(\lambda_i-\gamma_i)\bar{t}_i|\,|\bar{t}_i|^{p-1} \\ &= \sum_{i=1}^{\infty}|s_i-t_i|\,|\bar{t}_i|^{p-1} \\ &\le \left(\sum_{i=1}^{\infty}|s_i-t_i|^p\right)^{\frac{1}{p}}\left(\sum_{i=1}^{\infty}|\bar{t}_i|^{(p-1)q}\right)^{\frac{1}{q}} \end{aligned}$$

$$= \|v-u\|_p \left(\sum_{i=1}^{\infty} |\bar{t}_i|^p\right)^{\frac{1}{q}}$$

$$= \|\bar{x}\|_p^{\frac{p}{q}} \|v-u\|_p < \infty. \qquad (4.2)$$

Similarly to the proof of (I) in the proof of Theorem 3.2, the we estimate the numerator of the first part inside (4.1) without the constant $p-1$.

$$\sum_{n=1}^{\infty} \left| |\bar{t}_n|^{p-1} \operatorname{sign}(\bar{t}_n)(\lambda_n - \gamma_n) \right|^q$$

$$= \sum_{n=1}^{\infty} \left| |\bar{t}_n|^{p-2} |\bar{t}_n| \operatorname{sign}(\bar{t}_n)(\lambda_n - \gamma_n) \right|^q$$

$$= \sum_{n=1}^{\infty} \left| (s_n - t_n) |\bar{t}_n|^{p-2} \right|^q$$

$$= \sum_{n=1}^{\infty} |s_n - t_n|^q |\bar{t}_n|^{(p-2)q}$$

$$\leq \left( \sum_{n=1}^{\infty} (|s_n - t_n|^q)^{\frac{p}{q}} \right)^{\frac{q}{p}} \left( \sum_{n=1}^{\infty} \left( |\bar{t}_n|^{(p-2)q} \right)^{\frac{p}{p-q}} \right)^{\frac{p-q}{p}} \quad (2 < p < \infty \Longrightarrow 1 < q < 2)$$

$$\leq \left( \sum_{n=1}^{\infty} |s_n - t_n|^p \right)^{\frac{q}{p}} \left( \sum_{n=1}^{\infty} |\bar{t}_n|^p \right)^{\frac{p-q}{p}} \qquad ((p-2)q = p-q)$$

$$= \|\bar{x}\|_p^{p-q} \|v-u\|_p^q. \qquad (4.3)$$

By (4.2) and (4.3), we estimate $\|J'(\bar{x})(v)\|_q$, by (3.3), as follows.

$$\|J'(\bar{x})(v-u)\|_q$$

$$= \left\| J(\bar{x}) \left( (p-1) \begin{pmatrix} \lambda_1 - \gamma_1 & 0 & \dots \\ 0 & \lambda_2 - \gamma_2 & 0 \\ \vdots & 0 & \ddots \end{pmatrix} - \frac{(p-2)\sum_{i=1}^{\infty} (\lambda_i - \gamma_i) |\bar{t}_i|^p}{\sum_{i=1}^{\infty} |\bar{t}_i|^p} \begin{pmatrix} 1 & 0 & \dots \\ 0 & 1 & 0 \\ \vdots & 0 & \ddots \end{pmatrix} \right) \right\|_q$$

$$\leq \left\| \left\{ \frac{|\bar{t}_n|^{p-1} \operatorname{sign}(\bar{t}_n)(p-1)(\lambda_n - \gamma_n)}{\left(\sum_{i=1}^{\infty} |\bar{t}_i|^p\right)^{\frac{p-2}{p}}} \right\}_{n=1}^{\infty} \right\|_q + \left\| J(\bar{x}) \frac{(p-2)\sum_{i=1}^{\infty} (\lambda_i - \gamma_i) |\bar{t}_i|^p}{\sum_{i=1}^{\infty} |\bar{t}_i|^p} \begin{pmatrix} 1 & 0 & \dots \\ 0 & 1 & 0 \\ \vdots & 0 & \ddots \end{pmatrix} \right\|_q$$

$$= \frac{p-1}{\|\bar{x}\|_p^{p-2}} \left\| \{ |\bar{t}_n|^{p-1} \operatorname{sign}(\bar{t}_n)(\lambda_n - \gamma_n) \}_{n=1}^{\infty} \right\|_q + \left| \frac{(p-2)\sum_{i=1}^{\infty} (\lambda_i - \gamma_i) |\bar{t}_i|^p}{\sum_{i=1}^{\infty} |\bar{t}_i|^p} \right| \|J(\bar{x})\|_q$$

$$\leq \frac{p-1}{\|\bar{x}\|_p^{p-2}} \left( \|\bar{x}\|_p^{p-q} \|v-u\|_p^q \right)^{\frac{1}{q}} + \frac{p-2}{\|\bar{x}\|_p^p} \|\bar{x}\|_p^{\frac{p}{q}} \|v-u\|_p \|\bar{x}\|_p \qquad (\|J(\bar{x})\|_q = \|\bar{x}\|_p)$$

$$= \left( \frac{p-1}{\|\bar{x}\|_p^{p-2}} \|\bar{x}\|_p^{\frac{p-q}{q}} + \frac{p-2}{\|\bar{x}\|_p^p} \|\bar{x}\|_p^{\frac{p}{q}+1} \right) \|v-u\|_p$$

$$= \left( p - 1 + (p-2) \|\bar{x}\|_p^{\frac{p}{q}+1-p} \right) \|v-u\|_p$$

$$= \left( p - 1 + (p-2) \|\bar{x}\|_p^q \right) \|v-u\|_p.$$

This proves (4.1). □

By the connection between Gâteaux and Fréchet differentiability, the results (4.1) of Lemma 4.2 leads us to raise the following question:

**Question 4.3**. Let $2 < p < \infty$. Let $\bar{x} = (\bar{t}_1, \bar{t}_2, \dots) \in l_p \backslash \{\theta\}$. Let $v = \{s_n\}_{n=1}^{\infty} \in l_p \backslash \{\theta\}$ with ratio sequence $\{\lambda_n\}_{n=1}^{\infty}$ with respect to $\bar{x}$, Does $\nabla J(\bar{x})(v)$ exist with

$$\nabla J(\bar{x})(v) = J'(\bar{x})(v) = J(\bar{x})\left( (p-1)\begin{pmatrix} \lambda_1 & 0 & \dots \\ 0 & \lambda_2 & 0 \\ \vdots & 0 & \ddots \end{pmatrix} - \frac{(p-2)\sum_{i=1}^{\infty} \lambda_i |\bar{t}_i|^p}{\sum_{i=1}^{\infty} |\bar{t}_i|^p} \begin{pmatrix} 1 & 0 & \dots \\ 0 & 1 & 0 \\ \vdots & 0 & \ddots \end{pmatrix} \right) ?$$

We believe the answer for Question 4.3 is yes. But we will only prove a special case. Next, we prove the Fréchet differentiability of $J$ at $\bar{x} \in l_p \backslash \{\theta\}$ for $|\mathcal{N}(\bar{x})| < \infty$.

**Theorem. 4.4**. *Let $2 < p < \infty$. Let $\bar{x} = (\bar{t}_1, \bar{t}_2, \dots) \in l_p \backslash \{\theta\}$ with $1 \le |\mathcal{N}(\bar{x})| < \infty$. Then $J$ is Fréchet differentiable at $\bar{x}$. For any $v = \{s_n\}_{n=1}^{\infty} \in l_p \backslash \{\theta\}$ with ratio sequence $\{\lambda_n\}_{n=1}^{\infty}$ with respect to $\bar{x}$, one has*

$$\nabla J(\bar{x})(v) = J'(\bar{x})(v) = J(\bar{x})\left( (p-1)\begin{pmatrix} \lambda_1 & 0 & \dots \\ 0 & \lambda_2 & 0 \\ \vdots & 0 & \ddots \end{pmatrix} - \frac{(p-2)\sum_{i=1}^{\infty} \lambda_i |\bar{t}_i|^p}{\sum_{i=1}^{\infty} |\bar{t}_i|^p} \begin{pmatrix} 1 & 0 & \dots \\ 0 & 1 & 0 \\ \vdots & 0 & \ddots \end{pmatrix} \right).$$

*Proof*. We divide the proof of this theorem to two cases with respect to the value of $|\mathcal{N}(\bar{x})|$.

Case 1. $|\mathcal{N}(\bar{x})| = 1$. Let $\bar{x} = (\bar{t}_1, \bar{t}_2, \dots) \in l_p \backslash \{\theta\}$ with $|\mathcal{N}(\bar{x})| = 1$. without loss of generality, we suppose that $\bar{t}_1 \neq 0$ and $\bar{t}_n = 0$, for all $n = 2, 3, \dots$. Let $v = \{s_n\}_{n=1}^{\infty} \in l_p \backslash \{\theta\}$. Since $\bar{x} = (\bar{t}_1, 0, 0, \dots)$, then, by definition, the ratio sequence of $v$ with respect to $\bar{x}$ is $(\lambda_1, 0, 0, \dots)$ and $s_1 = \lambda_1 \bar{t}_1$, in which $\lambda_1$ is a real number. Since $|\mathcal{N}(\bar{x})| = 1$, $J(\bar{x}) = |\bar{t}_1| \mathrm{sign}(\bar{t}_1) = \bar{t}_1$. By (3.4), we have

$$\frac{J(\bar{x}+v) - J(\bar{x}) - J'(\bar{x})(v)}{\|v\|_p} = \frac{\left\{ \frac{|\bar{t}_1+s_1|^{p-1}\mathrm{sign}(\bar{t}_1+s_1)}{\left(\sum_{i=1}^{\infty}|\bar{t}_i+s_i|^p\right)^{\frac{p-2}{p}}} - |\bar{t}_1|\mathrm{sign}(\bar{t}_1) - \lambda_1 \bar{t}_1,\ \frac{|s_2|^{p-1}\mathrm{sign}(s_2)}{\left(\sum_{i=1}^{\infty}|\bar{t}_i+s_i|^p\right)^{\frac{p-2}{p}}},\ \frac{|s_3|^{p-1}\mathrm{sign}(s_3)}{\left(\sum_{i=1}^{\infty}|\bar{t}_i+s_i|^p\right)^{\frac{p-2}{p}}},\ \dots \right\}}{\|v\|_p}.$$

The first entry of $\frac{J(\bar{x}+v) - J(\bar{x}) - J'(\bar{x})(v)}{\|v\|_p}$ is

$$\frac{\frac{|\bar{t}_1+s_1|^{p-1}\mathrm{sign}(\bar{t}_1+s_1)}{\left(\sum_{i=1}^{\infty}|\bar{t}_i+s_i|^p\right)^{\frac{p-2}{p}}} - |\bar{t}_1|\mathrm{sign}(\bar{t}_1) - \lambda_1 \bar{t}_1}{\|v\|_p}$$

$$= \frac{\frac{|\bar{t}_1+s_1|^{p-1}\mathrm{sign}(\bar{t}_1+s_1)}{\left(|\bar{t}_1+s_1|^p + \sum_{i=2}^{\infty}|\bar{t}_i+s_i|^p\right)^{\frac{p-2}{p}}} - |\bar{t}_1|\mathrm{sign}(\bar{t}_1) - \lambda_1 \bar{t}_1}{\|v\|_p}$$

$$= \frac{\frac{|\bar{t}_1+s_1|^{p-1}\mathrm{sign}(\bar{t}_1+s_1)}{\left(|\bar{t}_1+s_1|^p + \sum_{i=2}^{\infty}|s_i|^p\right)^{\frac{p-2}{p}}} - |\bar{t}_1|\mathrm{sign}(\bar{t}_1) - s_1}{\left(|s_1|^p + \sum_{i=2}^{\infty}|s_i|^p\right)^{\frac{1}{p}}}$$

$$= \frac{\frac{|\bar{t}_1+s_1|^{p-1}\mathrm{sign}(\bar{t}_1+s_1)}{\left(|\bar{t}_1+s_1|^p+\sum_{i=2}^{\infty}|s_i|^p\right)^{\frac{p-2}{p}}} - |\bar{t}_1|\mathrm{sign}(\bar{t}_1) - s_1}{\left(|s_1|^p+\sum_{i=2}^{\infty}|s_i|^p\right)^{\frac{1}{p}}}. \tag{4.4}$$

It is considered as a function of two variables $s_1$ and $\sum_{i=2}^{\infty}|s_i|^p$. The first term of the numerator in this function is written by

$$f\left(s_1, \sum_{i=2}^{\infty}|s_i|^p\right) = \frac{|\bar{t}_1+s_1|^{p-1}\mathrm{sign}(\bar{t}_1+s_1)}{\left(|\bar{t}_1+s_1|^p+\sum_{i=2}^{\infty}|s_i|^p\right)^{\frac{p-2}{p}}}.$$

We calculate the partial derivatives of $f$ with respect to two variables $s_1$ and $\alpha \coloneqq \sum_{i=2}^{\infty}|s_i|^p$.

Subcase 1.1. $\bar{t}_1 + s_1 \geq 0$.

$$\frac{\partial f(s_1, \alpha)}{\partial s_1} = \mathrm{sign}(\bar{t}_1+s_1)\frac{(p-1)|\bar{t}_1+s_1|^{p-2}\left(|\bar{t}_1+s_1|^p+\sum_{i=2}^{\infty}|s_i|^p\right) - |\bar{t}_1+s_1|^{p-1}\left(\frac{p-2}{p}p|\bar{t}_1+s_1|^{p-1}\right)}{\left(|\bar{t}_1+s_1|^p+\sum_{i=2}^{\infty}|s_i|^p\right)^{\frac{p-2}{p}+1}},$$

$$\frac{\partial f(s_1, \alpha)}{\partial \alpha} = -\mathrm{sign}(\bar{t}_1)\frac{p|\bar{t}_1+s_1|^{p-1}\frac{p-2}{p}}{\left(|\bar{t}_1+s_1|^p+\sum_{i=2}^{\infty}|s_i|^p\right)^{\frac{p-2}{p}+1}}.$$

Subcase 1.2. $\bar{t}_1 + s_1 < 0$.

$$\frac{\partial f(s_1, \alpha)}{\partial s_1} = \mathrm{sign}(\bar{t}_1+s_1)\frac{(-1)(p-1)|\bar{t}_1+s_1|^{p-2}\left(|\bar{t}_1+s_1|^p+\sum_{i=2}^{\infty}|s_i|^p\right) - |\bar{t}_1+s_1|^{p-1}\left(\frac{p-2}{p}p(-1)|\bar{t}_1+s_1|^{p-1}\right)}{\left(|\bar{t}_1+s_1|^p+\sum_{i=2}^{\infty}|s_i|^p\right)^{\frac{p-2}{p}+1}},$$

$$\frac{\partial f(s_1, \alpha)}{\partial \alpha} = -\mathrm{sign}(\bar{t}_1)\frac{p(-1)|\bar{t}_1+s_1|^{p-1}\frac{p-2}{p}}{\left(|\bar{t}_1+s_1|^p+\sum_{i=2}^{\infty}|s_i|^p\right)^{\frac{p-2}{p}+1}}.$$

Combining the above two cases, we obtain that

$$\frac{\partial f(s_1, \alpha)}{\partial s_1} = \mathrm{sign}(\bar{t}_1+s_1)\frac{(p-1)|\bar{t}_1+s_1|^{p-2}\left(|\bar{t}_1+s_1|^p+\sum_{i=2}^{\infty}|s_i|^p\right) - (p-2)|\bar{t}_1+s_1|^{p-1}|\bar{t}_1+s_1|^{p-1}}{\left(|\bar{t}_1+s_1|^p+\sum_{i=2}^{\infty}|s_i|^p\right)^{\frac{p-2}{p}+1}}\mathrm{sign}(\bar{t}_1+s_1)$$

$$= \frac{(p-1)|\bar{t}_1+s_1|^{p-2}\left(|\bar{t}_1+s_1|^p+\sum_{i=2}^{\infty}|s_i|^p\right) - (p-2)|\bar{t}_1+s_1|^{p-1}|\bar{t}_1+s_1|^{p-1}}{\left(|\bar{t}_1+s_1|^p+\sum_{i=2}^{\infty}|s_i|^p\right)^{\frac{p-2}{p}+1}}.$$

$$\frac{\partial f(s_1, \alpha)}{\partial \alpha} = -\mathrm{sign}(\bar{t}_1)\frac{(p-2)|\bar{t}_1+s_1|^{p-1}\,\mathrm{sign}(\bar{t}_1+s_1)}{\left(|\bar{t}_1+s_1|^p+\sum_{i=2}^{\infty}|s_i|^p\right)^{\frac{p-2}{p}+1}} = -\frac{(p-2)|\bar{t}_1+s_1|^{p-1}}{\left(|\bar{t}_1+s_1|^p+\sum_{i=2}^{\infty}|s_i|^p\right)^{\frac{p-2}{p}+1}}. \tag{4.5}$$

Since $\bar{t}_1 \neq 0$, when $s_1$ is small enough, $\mathrm{sign}(\bar{t}_1+s_1) = \mathrm{sign}(\bar{t}_1)$.

Since $f(0,0) = |\bar{t}_1|\mathrm{sign}(\bar{t}_1)$. By (4.5), applying the Taylor's Approximation Theorem to the numerator in the limit (4.4) with respect to two variables $s_1$ and $\sum_{i=2}^{\infty}|s_i|^p$, there is $t_1$ with $|t_1| \leq |s_1|$ and $0 \leq \beta \leq \alpha \coloneqq \sum_{i=2}^{\infty}|s_i|^p$ such that

$$\frac{|\bar{t}_1+s_1|^{p-1}\mathrm{sign}(\bar{t}_1+s_1)}{\left(|\bar{t}_1+s_1|^p+\sum_{i=2}^{\infty}|s_i|^p\right)^{\frac{p-2}{p}}} - |\bar{t}_1|\mathrm{sign}(\bar{t}_1) - s_1$$

$$= f(s_1, \textstyle\sum_{i=2}^{\infty}|s_i|^p) - f(0,0) - s_1$$

$$= \frac{\partial f(s_1,\ \alpha)}{\partial s_1}(t_1,\beta)s_1 + \frac{\partial f(s_1,\ \alpha)}{\partial \alpha}(t_1,\beta)\textstyle\sum_{i=2}^{\infty}|s_i|^p - s_1$$

$$= \frac{(p-1)|\bar{t}_1+s_1|^{p-2}\left(|\bar{t}_1+s_1|^p+\sum_{i=2}^{\infty}|s_i|^p\right)-(p-2)|\bar{t}_1+s_1|^{p-1}|\bar{t}_1+s_1|^{p-1}}{\left(|\bar{t}_1+s_1|^p+\sum_{i=2}^{\infty}|s_i|^p\right)^{\frac{p-2}{p}+1}}s_1 - \frac{(p-2)|\bar{t}_1+s_1|^{p-1}}{\left(|\bar{t}_1+s_1|^p+\sum_{i=2}^{\infty}|s_i|^p\right)^{\frac{p-2}{p}+1}}\textstyle\sum_{i=2}^{\infty}|s_i|^p - s_1$$

$$= \left(\frac{(p-1)|\bar{t}_1+t_1|^{p-2}(|\bar{t}_1+t_1|^p+\beta)-|\bar{t}_1+s_1|^{p-1}\left(\frac{p-2}{p}p|\bar{t}_1+t_1|^{p-1}\right)}{(|\bar{t}_1+t_1|^p+\beta)^{\frac{p-2}{p}+1}} - 1\right)s_1 - \frac{|\bar{t}_1+t_1|^{p-1}\frac{p-2}{p}}{(|\bar{t}_1+t_1|^p+\beta)^{\frac{p-2}{p}+1}}\textstyle\sum_{i=2}^{\infty}|s_i|^p.$$

Notice that

$$\textstyle\sum_{i=1}^{\infty}|s_i|^p \to 0 \implies \sum_{i=2}^{\infty}|s_i|^p \to 0 \text{ and } s_1 \to 0.$$

Let $\sum_{i=1}^{\infty}|s_i|^p \to 0$ in the first entry in the limit $\lim_{v\to 0}\frac{J(\bar{x}+v)-J(\bar{x})-J'(\bar{x})(v)}{\|v\|_p}$. We have

$$\lim_{v\to 0}\frac{\frac{|\bar{t}_1+s_1|^{p-1}\text{sign}(\bar{t}_1+s_1)}{\left(|\bar{t}_1+s_1|^p+\sum_{i=2}^{\infty}|s_i|^p\right)^{\frac{p-2}{p}}} - |\bar{t}_1|\text{sign}(\bar{t}_1) - s_1}{\left(|s_1|^p+\sum_{i=2}^{\infty}|s_i|^p\right)^{\frac{1}{p}}}$$

$$= \lim_{v\to 0}\frac{\left(\frac{(p-1)|\bar{t}_1+t_1|^{p-2}(|\bar{t}_1+t_1|^p+\beta)-|\bar{t}_1+s_1|^{p-1}\left(\frac{p-2}{p}p|\bar{t}_1+t_1|^{p-1}\right)}{(|\bar{t}_1+t_1|^p+\beta)^{\frac{p-2}{p}+1}} - 1\right)s_1 - \frac{|\bar{t}_1+t_1|^{p-1}\frac{p-2}{p}}{(|\bar{t}_1+t_1|^p+\beta)^{\frac{p-2}{p}+1}}\sum_{i=2}^{\infty}|s_i|^p}{\left(\sum_{i=1}^{\infty}|s_i|^p\right)^{\frac{1}{p}}}$$

$$= \lim_{v\to 0}\left[\left(\frac{(p-1)|\bar{t}_1+t_1|^{p-2}(|\bar{t}_1+t_1|^p+\beta)-|\bar{t}_1+s_1|^{p-1}\left(\frac{p-2}{p}p|\bar{t}_1+t_1|^{p-1}\right)}{(|\bar{t}_1+t_1|^p+\beta)^{\frac{p-2}{p}+1}} - 1\right)\frac{s_1}{\left(\sum_{i=1}^{\infty}|s_i|^p\right)^{\frac{1}{p}}} - \frac{|\bar{t}_1+t_1|^{p-1}\frac{p-2}{p}}{(|\bar{t}_1+t_1|^p+\beta)^{\frac{p-2}{p}+1}}\frac{\sum_{i=2}^{\infty}|s_i|^p}{\left(\sum_{i=1}^{\infty}|s_i|^p\right)^{\frac{1}{p}}}\right]$$

$$= \lim_{v\to 0}\left[\left(\frac{(p-1)|\bar{t}_1+t_1|^{p-2}(|\bar{t}_1+t_1|^p+\beta)-|\bar{t}_1+s_1|^{p-1}\left(\frac{p-2}{p}p|\bar{t}_1+t_1|^{p-1}\right)}{(|\bar{t}_1+t_1|^p+\beta)^{\frac{p-2}{p}+1}} - 1\right)\frac{s_1}{\left(\sum_{i=1}^{\infty}|s_i|^p\right)^{\frac{1}{p}}} - \frac{|\bar{t}_1+t_1|^{p-1}\frac{p-2}{p}}{(|\bar{t}_1+t_1|^p+\beta)^{\frac{p-2}{p}+1}}\left(\sum_{i=1}^{\infty}|s_i|^p\right)^{\frac{p-1}{p}}\right]$$

$$= \left[\left(\frac{(p-1)|\bar{t}_1|^{p-2}(|\bar{t}_1|^p)-|\bar{t}_1|^{p-1}\left(\frac{p-2}{p}p|\bar{t}_1|^{p-1}\right)}{(|\bar{t}_1|^p)^{\frac{p-2}{p}+1}} - 1\right)1 - \frac{|\bar{t}_1|^{p-1}\frac{p-2}{p}}{(|\bar{t}_1|^p)^{\frac{p-2}{p}+1}}\,0\right]$$

$$= \left[\left(\frac{(p-1)|\bar{t}_1|^{2p-2}-(p-2)|\bar{t}_1|^{2p-2}}{(|\bar{t}_1|^p)^{\frac{p-2}{p}+1}} - 1\right)1 - \frac{|\bar{t}_1|^{p-1}\frac{p-2}{p}}{(|\bar{t}_1|^p)^{\frac{p-2}{p}+1}}\,0\right]$$

$$= \left[(1 - 1)1 - \frac{|\bar{t}_1|^{p-1}\frac{p-2}{p}}{(|\bar{t}_1|^p)^{\frac{p-2}{p}+1}}\ 0\right]$$

$$= 0. \tag{4.6}$$

Recall that

$$\frac{J(\bar{x}+v)-J(\bar{x})-J'(\bar{x})(v)}{\|v\|_p}$$

$$= \frac{\left\{\frac{|\bar{t}_1+s_1|^{p-1}\mathrm{sign}(\bar{t}_1+s_1)}{\left(\sum_{i=1}^{\infty}|\bar{t}_i+s_i|^p\right)^{\frac{p-2}{p}}} - |\bar{t}_1|\mathrm{sign}(\bar{t}_1) - \lambda_1\bar{t}_1,\ \frac{|s_2|^{p-1}\mathrm{sign}(s_2)}{\left(\sum_{i=1}^{\infty}|\bar{t}_i+s_i|^p\right)^{\frac{p-2}{p}}},\ \frac{|s_3|^{p-1}\mathrm{sign}(s_3)}{\left(\sum_{i=1}^{\infty}|\bar{t}_i+s_i|^p\right)^{\frac{p-2}{p}}},\ \ldots\right\}}{\|v\|_p}.$$

Let $\sum_{i=1}^{\infty}|s_i|^p \to 0$ in the rest entries except the first one in the limit $\lim_{v\to 0}\frac{J(\bar{x}+v)-J(\bar{x})-J'(\bar{x})(v)}{\|v\|_p}$.

$$\lim_{v\to 0}\frac{\left\|\left\{0,\ \frac{|s_2|^{p-1}\mathrm{sign}(s_2)}{\left(\sum_{i=1}^{\infty}|\bar{t}_i+s_i|^p\right)^{\frac{p-2}{p}}},\ \frac{|s_3|^{p-1}\mathrm{sign}(s_3)}{\left(\sum_{i=1}^{\infty}|\bar{t}_i+s_i|^p\right)^{\frac{p-2}{p}}},\ \ldots\right\}\right\|_q}{\|v\|_p} = \lim_{v\to 0}\frac{\left(\sum_{i=2}^{\infty}\left|\frac{|s_i|^{p-1}\mathrm{sign}(s_i)}{\left(|\bar{t}_1+s_1|^p+\sum_{i=2}^{\infty}|s_i|^p\right)^{\frac{p-2}{p}}}\right|^q\right)^{\frac{1}{q}}}{\|v\|_p}$$

$$= \lim_{v\to 0}\frac{\frac{\left(\sum_{i=2}^{\infty}|s_i|^p\right)^{\frac{1}{q}}}{\left(|\bar{t}_1+s_1|^p+\sum_{i=2}^{\infty}|s_i|^p\right)^{\frac{p-2}{p}}}}{\left(\sum_{i=1}^{\infty}|s_i|^p\right)^{\frac{1}{p}}} \le \lim_{v\to 0}\frac{\frac{\left(\sum_{i=1}^{\infty}|s_i|^p\right)^{\frac{1}{q}}}{\left(|\bar{t}_1+s_1|^p+\sum_{i=2}^{\infty}|s_i|^p\right)^{\frac{p-2}{p}}}}{\left(\sum_{i=1}^{\infty}|s_i|^p\right)^{\frac{1}{p}}} = \lim_{v\to 0}\frac{\left(\sum_{i=1}^{\infty}|s_i|^p\right)^{\frac{1}{q}-\frac{1}{p}}}{\left(|\bar{t}_1+s_1|^p+\sum_{i=2}^{\infty}|s_i|^p\right)^{\frac{p-2}{p}}}.$$

Since $\lim_{v\to 0}\frac{1}{\left(|\bar{t}_1+s_1|^p+\sum_{i=2}^{\infty}|s_i|^p\right)^{\frac{p-2}{p}}} = \lim_{v\to 0}\frac{1}{|\bar{t}_1|^{p-2}} > 0$ and $\frac{1}{q} - \frac{1}{p} > 0$ (by $p > 2$), the above limit satisfies that

$$\lim_{v\to 0}\frac{\left(\sum_{i=1}^{\infty}|s_i|^p\right)^{\frac{1}{q}-\frac{1}{p}}}{\left(|\bar{t}_1+s_1|^p+\sum_{i=2}^{\infty}|s_i|^p\right)^{\frac{p-2}{p}}} = 0. \tag{4.7}$$

By (4.6) and (4.7), we obtain that

$$\lim_{v\to 0}\frac{J(\bar{x}+v)-J(\bar{x})-J'(\bar{x})(v)}{\|v\|_p} = \theta.$$

This completes the proof of case 1. Then, we prove case 2.

Case 2. Suppose that $|\mathcal{N}(\bar{x})| = k$, for some positive integer $k > 1$. Without loss of generality, we suppose that $\bar{t}_n \neq 0$, for $n = 1, 2, \ldots, k$ and $\bar{t}_n = 0$, for all $n > k$. Let $v = \{s_n\}_{n=1}^{\infty} \in l_p\backslash\{\theta\}$. Then, by definition, the ratio sequence of $v$ with respect to $\bar{x}$ is $(\lambda_1, \ldots \lambda_k, 0, 0, \ldots)$. For $n = 1, 2, \ldots, k$, $\lambda_n = \frac{s_n}{\bar{t}_n}$, which is a real number. In this case, (3.3) becomes

$$J'(\bar{x})(v) = J(\bar{x})\left((p-1)\begin{pmatrix}\lambda_1 & 0 & \ldots\\ 0 & \lambda_2 & 0\\ \vdots & 0 & \ddots\end{pmatrix} - \frac{(p-2)\sum_{i=1}^{k}\lambda_i|\bar{t}_i|^p}{\sum_{i=1}^{k}|\bar{t}_i|^p}\begin{pmatrix}1 & 0 & \ldots\\ 0 & 1 & 0\\ \vdots & 0 & \ddots\end{pmatrix}\right).$$

More precisely,

$$\text{the } n\text{th-entry of } J'(\bar{x})(v) = \begin{cases} \frac{|\bar{t}_n|^{p-1}\operatorname{sign}(\bar{t}_n)}{\left(\sum_{i=1}^{k}|\bar{t}_i|^p\right)^{\frac{p-2}{p}}}\left((p-1)\lambda_n - \frac{(p-2)\sum_{i=1}^{k}\lambda_i|\bar{t}_i|^p}{\sum_{i=1}^{k}|\bar{t}_i|^p}\right), & \text{if } n \le k, \\ 0, & \text{if } n > k. \end{cases}$$

This implies that

$J(\bar{x}+v) - J(\bar{x}) - J'(\bar{x})(v)$

$$= \Bigg\{ \frac{|\bar{t}_1+s_1|^{p-1}\operatorname{sign}(\bar{t}_1+s_1)}{\left(\sum_{i=1}^{\infty}|\bar{t}_i+s_i|^p\right)^{\frac{p-2}{p}}} - \frac{|\bar{t}_1|^{p-1}\operatorname{sign}(\bar{t}_1)}{\left(\sum_{i=1}^{\infty}|\bar{t}_i|^p\right)^{\frac{p-2}{p}}} - \frac{|\bar{t}_1|^{p-1}\operatorname{sign}(\bar{t}_1)}{\left(\sum_{i=1}^{\infty}|\bar{t}_i|^p\right)^{\frac{p-2}{p}}}\left((p-1)\lambda_1 - \frac{(p-2)\sum_{i=1}^{k}\lambda_i|\bar{t}_i|^p}{\sum_{i=1}^{k}|\bar{t}_i|^p}\right), \ldots$$

$$\ldots, \frac{|\bar{t}_k+s_k|^{p-1}\operatorname{sign}(\bar{t}_n+s_n)}{\left(\sum_{i=1}^{\infty}|\bar{t}_i+s_i|^p\right)^{\frac{p-2}{p}}} - \frac{|\bar{t}_k|^{p-1}\operatorname{sign}(\bar{t}_k)}{\left(\sum_{i=1}^{\infty}|\bar{t}_i|^p\right)^{\frac{p-2}{p}}} - \frac{|\bar{t}_k|^{p-1}\operatorname{sign}(\bar{t}_k)}{\left(\sum_{i=1}^{\infty}|\bar{t}_i|^p\right)^{\frac{p-2}{p}}}\left((p-1)\lambda_k - \frac{(p-2)\sum_{i=1}^{k}\lambda_i|\bar{t}_i|^p}{\sum_{i=1}^{k}|\bar{t}_i|^p}\right),$$

$$\frac{|s_{k+1}|^{p-1}\operatorname{sign}(s_{k+1})}{\left(\sum_{i=1}^{\infty}|\bar{t}_i+s_i|^p\right)^{\frac{p-2}{p}}}, \frac{|s_{k+2}|^{p-1}\operatorname{sign}(s_{k+2})}{\left(\sum_{i=1}^{\infty}|\bar{t}_i+s_i|^p\right)^{\frac{p-2}{p}}}, \ldots \Bigg\}$$

$$= \Bigg\{ \frac{|\bar{t}_1+s_1|^{p-1}\operatorname{sign}(\bar{t}_1+s_1)}{\left(\sum_{i=1}^{k}|\bar{t}_i+s_i|^p+\sum_{i=k+1}^{\infty}|s_i|^p\right)^{\frac{p-2}{p}}} - \frac{|\bar{t}_1|^{p-1}\operatorname{sign}(\bar{t}_1)}{\left(\sum_{i=1}^{k}|\bar{t}_i|^p\right)^{\frac{p-2}{p}}} - \frac{|\bar{t}_1|^{p-1}\operatorname{sign}(\bar{t}_1)}{\left(\sum_{i=1}^{k}|\bar{t}_i|^p\right)^{\frac{p-2}{p}}}\left((p-1)\lambda_1 - \frac{(p-2)\sum_{i=1}^{k}\lambda_i|\bar{t}_i|^p}{\sum_{i=1}^{k}|\bar{t}_i|^p}\right), \ldots$$

$$\ldots, \frac{|\bar{t}_k+s_k|^{p-1}\operatorname{sign}(\bar{t}_n+s_n)}{\left(\sum_{i=1}^{k}|\bar{t}_i+s_i|^p+\sum_{i=k+1}^{\infty}|s_i|^p\right)^{\frac{p-2}{p}}} - \frac{|\bar{t}_k|^{p-1}\operatorname{sign}(\bar{t}_k)}{\left(\sum_{i=1}^{k}|\bar{t}_i|^p\right)^{\frac{p-2}{p}}} - \frac{|\bar{t}_k|^{p-1}\operatorname{sign}(\bar{t}_k)}{\left(\sum_{i=1}^{k}|\bar{t}_i|^p\right)^{\frac{p-2}{p}}}\left((p-1)\lambda_k - \frac{(p-2)\sum_{i=1}^{k}\lambda_i|\bar{t}_i|^p}{\sum_{i=1}^{k}|\bar{t}_i|^p}\right),$$

$$\frac{|s_{k+1}|^{p-1}\operatorname{sign}(s_{k+1})}{\left(\sum_{i=1}^{k}|\bar{t}_i+s_i|^p+\sum_{i=k+1}^{\infty}|s_i|^p\right)^{\frac{p-2}{p}}}, \frac{|s_{k+2}|^{p-1}\operatorname{sign}(s_{k+2})}{\left(\sum_{i=1}^{k}|\bar{t}_i+s_i|^p+\sum_{i=k+1}^{\infty}|s_i|^p\right)^{\frac{p-2}{p}}}, \ldots \Bigg\}$$

For $n$ = 1, 2, …, $k$, the $n$th-entry of $J(\bar{x}+v) - J(\bar{x}) - J'(\bar{x})(v)$ satisfies that

$$\frac{|\bar{t}_n+s_n|^{p-1}\operatorname{sign}(\bar{t}_n+s_n)}{\left(\sum_{i=1}^{k}|\bar{t}_i+s_i|^p+\sum_{i=k+1}^{\infty}|s_i|^p\right)^{\frac{p-2}{p}}} - \frac{|\bar{t}_n|^{p-1}\operatorname{sign}(\bar{t}_n)}{\left(\sum_{i=1}^{k}|\bar{t}_i|^p\right)^{\frac{p-2}{p}}} - \frac{|\bar{t}_n|^{p-1}\operatorname{sign}(\bar{t}_n)}{\left(\sum_{i=1}^{k}|\bar{t}_i|^p\right)^{\frac{p-2}{p}}}\left((p-1)\lambda_n - \frac{(p-2)\sum_{i=1}^{k}\lambda_i|\bar{t}_i|^p}{\sum_{i=1}^{k}|\bar{t}_i|^p}\right)$$

$$= \frac{|\bar{t}_n+s_n|^{p-1}\operatorname{sign}(\bar{t}_n+s_n)}{\left(\sum_{i=1}^{k}|\bar{t}_i+s_i|^p+\sum_{i=k+1}^{\infty}|s_i|^p\right)^{\frac{p-2}{p}}} - \frac{|\bar{t}_n|^{p-1}\operatorname{sign}(\bar{t}_n)}{\left(\sum_{i=1}^{k}|\bar{t}_i|^p\right)^{\frac{p-2}{p}}} - \frac{(p-1)\lambda_n|\bar{t}_n|^{p-1}\operatorname{sign}(\bar{t}_n)}{\left(\sum_{i=1}^{k}|\bar{t}_i|^p\right)^{\frac{p-2}{p}}} + \frac{(p-2)|\bar{t}_n|^{p-1}\operatorname{sign}(\bar{t}_n)\sum_{i=1}^{k}\lambda_i|\bar{t}_i|^p}{\left(\sum_{i=1}^{k}|\bar{t}_i|^p\right)^{\frac{p-2}{p}+1}}$$

For $n$ = 1, 2, …, $k$, we write a function

$$f_n(s_1, s_2, \ldots, s_k, \alpha) = \frac{|\bar{t}_n+s_n|^{p-1}\operatorname{sign}(\bar{t}_n+s_n)}{\left(\sum_{i=1}^{k}|\bar{t}_i+s_i|^p+\sum_{i=k+1}^{\infty}|s_i|^p\right)^{\frac{p-2}{p}}} = \frac{|\bar{t}_n+s_n|^{p-1}\operatorname{sign}(\bar{t}_n+s_n)}{\left(\sum_{i=1}^{\infty}|\bar{t}_i+s_i|^p+\alpha\right)^{\frac{p-2}{p}}}$$

Taking the partial derivative of $f_n$ with respect to $s_1, s_2, \ldots, s_k, \alpha$, we have

Subcase 2.1. The partial derivative of $f_n$ with respect to $s_j$, for $j$ = 1, 2, …, $k$ and $j \neq n$.

$$\frac{\partial f_n(s_1,s_2,\dots,s_k,\alpha)}{\partial s_j} = -\frac{|\bar{t}_n+s_n|^{p-1}\operatorname{sign}(\bar{t}_n+s_n)\frac{p-2}{p}p\left|\bar{t}_j+s_j\right|^{p-1}\operatorname{sign}\left(\bar{t}_j+s_j\right)}{\left(\sum_{i=1}^{k}|\bar{t}_i+s_i|^p+\sum_{i=k+1}^{\infty}|s_i|^p\right)^{\frac{p-2}{p}+1}}$$

$$= -\frac{(p-2)|\bar{t}_n+s_n|^{p-1}\operatorname{sign}(\bar{t}_n+s_n)\left|\bar{t}_j+s_j\right|^{p-1}\operatorname{sign}\left(\bar{t}_j+s_j\right)}{\left(\sum_{i=1}^{k}|\bar{t}_i+s_i|^p+\alpha\right)^{\frac{p-2}{p}+1}}$$

$$= -\frac{(p-2)|\bar{t}_n+s_n|^{p-1}\left|\bar{t}_j+s_j\right|^{p-1}\operatorname{sign}(\bar{t}_n+s_n)\operatorname{sign}\left(\bar{t}_j+s_j\right)}{\left(\sum_{i=1}^{k}|\bar{t}_i+s_i|^p+\alpha\right)^{\frac{p-2}{p}+1}}$$

Subcase 2.2. The partial derivative of $f_n$ with respect to $s_j$, for $j = n$.

$$\frac{\partial f_n(s_1,s_2,\dots,s_k,\alpha)}{\partial s_n} = \frac{(p-1)|\bar{t}_n+s_n|^{p-2}\operatorname{sign}(\bar{t}_n+s_n)\operatorname{sign}(\bar{t}_n+s_n)\left(\sum_{i=1}^{k}|\bar{t}_i+s_i|^p+\sum_{i=k+1}^{\infty}|s_i|^p\right)}{\left(\sum_{i=1}^{k}|\bar{t}_i+s_i|^p+\sum_{i=k+1}^{\infty}|s_i|^p\right)^{\frac{p-2}{p}+1}}$$

$$-\frac{|\bar{t}_n+s_n|^{p-1}\operatorname{sign}(\bar{t}_n+s_n)\frac{p-2}{p}p|\bar{t}_n+s_n|^{p-1}\operatorname{sign}(\bar{t}_n+s_n)}{\left(\sum_{i=1}^{k}|\bar{t}_i+s_i|^p+\sum_{i=k+1}^{\infty}|s_i|^p\right)^{\frac{p-2}{p}+1}}$$

$$= \frac{(p-1)|\bar{t}_n+s_n|^{p-2}\left(\sum_{i=1}^{k}|\bar{t}_i+s_i|^p+\sum_{i=k+1}^{\infty}|s_i|^p\right)-(p-2)|\bar{t}_n+s_n|^{p-1}|\bar{t}_n+s_n|^{p-1}}{\left(\sum_{i=1}^{k}|\bar{t}_i+s_i|^p+\sum_{i=k+1}^{\infty}|s_i|^p\right)^{\frac{p-2}{p}+1}}$$

$$= \frac{(p-1)|\bar{t}_n+s_n|^{p-2}\left(\sum_{i=1}^{k}|\bar{t}_i+s_i|^p+\sum_{i=k+1}^{\infty}|s_i|^p\right)-(p-2)|\bar{t}_n+s_n|^{2p-2}}{\left(\sum_{i=1}^{k}|\bar{t}_i+s_i|^p+\alpha\right)^{\frac{p-2}{p}+1}}$$

Subcase 2.3. The partial derivative of $f_n$ with respect to $\alpha$.

$$\frac{\partial f_n(s_1,s_2,\dots,s_k,\alpha)}{\partial \alpha} = -\frac{|\bar{t}_n+s_n|^{p-1}\operatorname{sign}(\bar{t}_n+s_n)\frac{p-2}{p}}{\left(\sum_{i=1}^{k}|\bar{t}_i+s_i|^p+\sum_{i=k+1}^{\infty}|s_i|^p\right)^{\frac{p-2}{p}+1}} = -\frac{\frac{p-2}{p}|\bar{t}_n+s_n|^{p-1}\operatorname{sign}(\bar{t}_n+s_n)}{\left(\sum_{i=1}^{k}|\bar{t}_i+s_i|^p+\alpha\right)^{\frac{p-2}{p}+1}}.$$

One sees that

$$f_n(0,0,\dots,0,0) = \frac{|\bar{t}_n|^{p-1}\operatorname{sign}(\bar{t}_n)}{\left(\sum_{i=1}^{k}|\bar{t}_i|^p\right)^{\frac{p-2}{p}}}.$$

By the above 3 subcases in (3.6), applying the Taylor's Approximation Theorem to $f$ with respect to $k + 1$ variables $s_1, s_2, \dots, s_k$ and $\alpha \coloneqq \sum_{i=2}^{\infty}|s_i|^p$, there are $t_j$ with $\left|t_j\right| \le \left|s_j\right|$ and $0 \le \beta \le \alpha \coloneqq \sum_{i=2}^{\infty}|s_i|^p$ such that

$$\frac{|\bar{t}_n+s_n|^{p-1}\operatorname{sign}(\bar{t}_n+s_n)}{\left(\sum_{i=1}^{k}|\bar{t}_i+s_i|^p+\sum_{i=k+1}^{\infty}|s_i|^p\right)^{\frac{p-2}{p}}} - \frac{|\bar{t}_n|^{p-1}\operatorname{sign}(\bar{t}_n)}{\left(\sum_{i=1}^{k}|\bar{t}_i|^p\right)^{\frac{p-2}{p}}}$$

$$= f_n(s_1, s_2, \dots, s_k, \alpha) - f_n(0,0,\dots,0,0)$$

$$= \sum_{1\le j\le k, j\ne n}\left(-\frac{(p-2)|\bar{t}_n+t_n|^{p-1}\left|\bar{t}_j+t_j\right|^{p-1}\operatorname{sign}(\bar{t}_n+t_n)\operatorname{sign}\left(\bar{t}_j+t_j\right)}{\left(\sum_{i=1}^{k}|\bar{t}_i+t_i|^p+\beta\right)^{\frac{p-2}{p}+1}}\right)s_j$$

$$+\frac{(p-1)|\bar{t}_n+t_n|^{p-2}\left(\sum_{i=1}^{k}|\bar{t}_i+t_i|^p+\sum_{i=k+1}^{\infty}|t_i|^p\right)-(p-2)|\bar{t}_n+t_n|^{2p-2}}{\left(\sum_{i=1}^{k}|\bar{t}_i+t_i|^p+\beta\right)^{\frac{p-2}{p}+1}}s_n-\frac{\frac{p-2}{p}|\bar{t}_n+t_n|^{p-1}\mathrm{sign}(\bar{t}_n+t_n)}{\left(\sum_{i=1}^{k}|\bar{t}_i+t_i|^p+\beta\right)^{\frac{p-2}{p}+1}}\alpha$$

This implies that

$$\frac{|\bar{t}_n+s_n|^{p-1}\mathrm{sign}(\bar{t}_n+s_n)}{\left(\sum_{i=1}^{k}|\bar{t}_i+s_i|^p+\sum_{i=k+1}^{\infty}|s_i|^p\right)^{\frac{p-2}{p}}}-\frac{|\bar{t}_n|^{p-1}\mathrm{sign}(\bar{t}_n)}{\left(\sum_{i=1}^{k}|\bar{t}_i|^p\right)^{\frac{p-2}{p}}}-\frac{(p-1)\lambda_n|\bar{t}_n|^{p-1}\mathrm{sign}(\bar{t}_n)}{\left(\sum_{i=1}^{k}|\bar{t}_i|^p\right)^{\frac{p-2}{p}}}+\frac{(p-2)|\bar{t}_n|^{p-1}\mathrm{sign}(\bar{t}_n)\sum_{i=1}^{k}\lambda_i|\bar{t}_i|^p}{\left(\sum_{i=1}^{k}|\bar{t}_i|^p\right)^{\frac{p-2}{p}+1}}$$

$$=\sum_{1\le j\le k,j\ne n}\left(-\frac{(p-2)|\bar{t}_n+t_n|^{p-1}|\bar{t}_j+t_j|^{p-1}\mathrm{sign}(\bar{t}_n+t_n)\mathrm{sign}(\bar{t}_j+t_j)}{\left(\sum_{i=1}^{k}|\bar{t}_i+t_i|^p+\beta\right)^{\frac{p-2}{p}+1}}\right)s_j$$

$$+\frac{(p-1)|\bar{t}_n+t_n|^{p-2}\left(\sum_{i=1}^{k}|\bar{t}_i+t_i|^p+\sum_{i=k+1}^{\infty}|t_i|^p\right)-(p-2)|\bar{t}_n+t_n|^{2p-2}}{\left(\sum_{i=1}^{k}|\bar{t}_i+t_i|^p+\beta\right)^{\frac{p-2}{p}+1}}s_n-\frac{\frac{p-2}{p}|\bar{t}_n+t_n|^{p-1}\mathrm{sign}(\bar{t}_n+t_n)}{\left(\sum_{i=1}^{k}|\bar{t}_i+t_i|^p+\beta\right)^{\frac{p-2}{p}+1}}\alpha$$

$$-\frac{(p-1)\lambda_n|\bar{t}_n|^{p-1}\mathrm{sign}(\bar{t}_n)}{\left(\sum_{i=1}^{k}|\bar{t}_i|^p\right)^{\frac{p-2}{p}}}+\frac{(p-2)|\bar{t}_n|^{p-1}\mathrm{sign}(\bar{t}_n)\sum_{i=1}^{k}\lambda_i|\bar{t}_i|^p}{\left(\sum_{i=1}^{k}|\bar{t}_i|^p\right)^{\frac{p-2}{p}+1}}$$

$$=\sum_{1\le j\le k,j\ne n}\left(-\frac{(p-2)|\bar{t}_n+t_n|^{p-1}|\bar{t}_j+t_j|^{p-1}\mathrm{sign}(\bar{t}_n+t_n)\mathrm{sign}(\bar{t}_j+t_j)}{\left(\sum_{i=1}^{k}|\bar{t}_i+t_i|^p+\beta\right)^{\frac{p-2}{p}+1}}\right)s_j$$

$$+\frac{(p-1)|\bar{t}_n+t_n|^{p-2}\left(\sum_{i=1}^{k}|\bar{t}_i+t_i|^p+\sum_{i=k+1}^{\infty}|t_i|^p\right)-(p-2)|\bar{t}_n+t_n|^{2p-2}}{\left(\sum_{i=1}^{k}|\bar{t}_i+t_i|^p+\beta\right)^{\frac{p-2}{p}+1}}s_n-\frac{\frac{p-2}{p}|\bar{t}_n+t_n|^{p-1}\mathrm{sign}(\bar{t}_n+t_n)}{\left(\sum_{i=1}^{k}|\bar{t}_i+t_i|^p+\beta\right)^{\frac{p-2}{p}+1}}\alpha$$

$$-\frac{(p-1)\lambda_n|\bar{t}_n|^{p-1}\mathrm{sign}(\bar{t}_n)}{\left(\sum_{i=1}^{k}|\bar{t}_i|^p\right)^{\frac{p-2}{p}}}+\frac{(p-2)|\bar{t}_n|^{p-1}\mathrm{sign}(\bar{t}_n)\sum_{1\le j\le k,j\ne n}\lambda_j|\bar{t}_j|^p}{\left(\sum_{i=1}^{k}|\bar{t}_i|^p\right)^{\frac{p-2}{p}+1}}+\frac{(p-2)|\bar{t}_n|^{p-1}\mathrm{sign}(\bar{t}_n)\lambda_n|\bar{t}_n|^p}{\left(\sum_{i=1}^{k}|\bar{t}_i|^p\right)^{\frac{p-2}{p}+1}}$$

$$=\sum_{1\le j\le k,j\ne n}\left(\left(-\frac{(p-2)|\bar{t}_n+t_n|^{p-1}|\bar{t}_j+t_j|^{p-1}\mathrm{sign}(\bar{t}_n+t_n)\mathrm{sign}(\bar{t}_j+t_j)}{\left(\sum_{i=1}^{k}|\bar{t}_i+t_i|^p+\beta\right)^{\frac{p-2}{p}+1}}\right)s_j+\frac{(p-2)|\bar{t}_n|^{p-1}\mathrm{sign}(\bar{t}_n)\lambda_j|\bar{t}_j|^p}{\left(\sum_{i=1}^{k}|\bar{t}_i|^p\right)^{\frac{p-2}{p}+1}}\right)$$

$$+\frac{(p-1)|\bar{t}_n+t_n|^{p-2}\left(\sum_{i=1}^{k}|\bar{t}_i+t_i|^p+\sum_{i=k+1}^{\infty}|t_i|^p\right)-(p-2)|\bar{t}_n+t_n|^{2p-2}}{\left(\sum_{i=1}^{k}|\bar{t}_i+t_i|^p+\beta\right)^{\frac{p-2}{p}+1}}s_n-\frac{\frac{p-2}{p}|\bar{t}_n+t_n|^{p-1}\mathrm{sign}(\bar{t}_n+t_n)}{\left(\sum_{i=1}^{k}|\bar{t}_i+t_i|^p+\beta\right)^{\frac{p-2}{p}+1}}\alpha$$

$$-\frac{(p-1)\lambda_n|\bar{t}_n|^{p-1}\mathrm{sign}(\bar{t}_n)}{\left(\sum_{i=1}^{k}|\bar{t}_i|^p\right)^{\frac{p-2}{p}}}+\frac{(p-2)|\bar{t}_n|^{p-1}\mathrm{sign}(\bar{t}_n)\lambda_n|\bar{t}_n|^p}{\left(\sum_{i=1}^{k}|\bar{t}_i|^p\right)^{\frac{p-2}{p}+1}}$$

$$=\sum_{1\le j\le k,j\ne n}\left(\left(-\frac{(p-2)|\bar{t}_n+t_n|^{p-1}|\bar{t}_j+t_j|^{p-1}\mathrm{sign}(\bar{t}_n+t_n)\mathrm{sign}(\bar{t}_j+t_j)}{\left(\sum_{i=1}^{k}|\bar{t}_i+t_i|^p+\beta\right)^{\frac{p-2}{p}+1}}\right)s_j+\frac{(p-2)|\bar{t}_n|^{p-1}|\bar{t}_j|^{p-1}\mathrm{sign}(\bar{t}_n)\lambda_j|\bar{t}_j|}{\left(\sum_{i=1}^{k}|\bar{t}_i|^p\right)^{\frac{p-2}{p}+1}}\right)$$

$$+\frac{(p-1)|\bar{t}_n+t_n|^{p-2}\left(\sum_{i=1}^{k}|\bar{t}_i+t_i|^p+\sum_{i=k+1}^{\infty}|t_i|^p\right)-(p-2)|\bar{t}_n+t_n|^{2p-2}}{\left(\sum_{i=1}^{k}|\bar{t}_i+t_i|^p+\beta\right)^{\frac{p-2}{p}+1}}s_n-\frac{\frac{p-2}{p}|\bar{t}_n+t_n|^{p-1}\mathrm{sign}(\bar{t}_n+t_n)}{\left(\sum_{i=1}^{k}|\bar{t}_i+t_i|^p+\beta\right)^{\frac{p-2}{p}+1}}\alpha$$

$$-\frac{(p-1)|\bar{t}_n|^{p-2}\left(\sum_{i=1}^{k}|\bar{t}_i|^p\right)\lambda_n \mathrm{sign}(\bar{t}_n)|\bar{t}_n|}{\left(\sum_{i=1}^{k}|\bar{t}_i|^p\right)^{\frac{p-2}{p}+1}}+\frac{(p-2)|\bar{t}_n|^{p-1}|\bar{t}_n|^{p-1}\mathrm{sign}(\bar{t}_n)\lambda_n|\bar{t}_n|}{\left(\sum_{i=1}^{k}|\bar{t}_i|^p\right)^{\frac{p-2}{p}+1}}$$

$$=\sum_{1\le j\le k,j\neq n}\left(\left(-\frac{(p-2)|\bar{t}_n+t_n|^{p-1}|\bar{t}_j+t_j|^{p-1}\mathrm{sign}(\bar{t}_n+t_n)\mathrm{sign}(\bar{t}_j+t_j)}{\left(\sum_{i=1}^{k}|\bar{t}_i+t_i|^p+\beta\right)^{\frac{p-2}{p}+1}}\right)s_j+\frac{(p-2)|\bar{t}_n|^{p-1}|\bar{t}_j|^{p-1}\mathrm{sign}(\bar{t}_n)\lambda_j|\bar{t}_j|}{\left(\sum_{i=1}^{k}|\bar{t}_i|^p\right)^{\frac{p-2}{p}+1}}\right)$$

$$+\frac{(p-1)|\bar{t}_n+t_n|^{p-2}\left(\sum_{i=1}^{k}|\bar{t}_i+t_i|^p+\sum_{i=k+1}^{\infty}|t_i|^p\right)-(p-2)|\bar{t}_n+t_n|^{2p-2}}{\left(\sum_{i=1}^{k}|\bar{t}_i+t_i|^p+\beta\right)^{\frac{p-2}{p}+1}}s_n-\frac{\frac{p-2}{p}|\bar{t}_n+t_n|^{p-1}\mathrm{sign}(\bar{t}_n+t_n)}{\left(\sum_{i=1}^{k}|\bar{t}_i+t_i|^p+\beta\right)^{\frac{p-2}{p}+1}}\alpha$$

$$-\frac{(p-1)|\bar{t}_n|^{p-2}\left(\sum_{i=1}^{k}|\bar{t}_i|^p\right)\lambda_n \mathrm{sign}(\bar{t}_n)|\bar{t}_n|}{\left(\sum_{i=1}^{k}|\bar{t}_i|^p\right)^{\frac{p-2}{p}+1}}+\frac{(p-2)|\bar{t}_n|^{2p-2}\mathrm{sign}(\bar{t}_n)\lambda_n|\bar{t}_n|}{\left(\sum_{i=1}^{k}|\bar{t}_i|^p\right)^{\frac{p-2}{p}+1}}$$

$$=\sum_{1\le j\le k,j\neq n}\left(\left(-\frac{(p-2)|\bar{t}_n+t_n|^{p-1}|\bar{t}_j+t_j|^{p-1}\mathrm{sign}(\bar{t}_n+t_n)\mathrm{sign}(\bar{t}_j+t_j)}{\left(\sum_{i=1}^{k}|\bar{t}_i+t_i|^p+\beta\right)^{\frac{p-2}{p}+1}}\right)s_j+\frac{(p-2)|\bar{t}_n|^{p-1}|\bar{t}_j|^{p-1}\mathrm{sign}(\bar{t}_n)\lambda_j\mathrm{sign}(\bar{t}_j)\bar{t}_j}{\left(\sum_{i=1}^{k}|\bar{t}_i|^p\right)^{\frac{p-2}{p}+1}}\right)$$

$$+\frac{(p-1)|\bar{t}_n+t_n|^{p-2}\left(\sum_{i=1}^{k}|\bar{t}_i+t_i|^p+\sum_{i=k+1}^{\infty}|t_i|^p\right)-(p-2)|\bar{t}_n+t_n|^{2p-2}}{\left(\sum_{i=1}^{k}|\bar{t}_i+t_i|^p+\beta\right)^{\frac{p-2}{p}+1}}s_n-\frac{\frac{p-2}{p}|\bar{t}_n+t_n|^{p-1}\mathrm{sign}(\bar{t}_n+t_n)}{\left(\sum_{i=1}^{k}|\bar{t}_i+t_i|^p+\beta\right)^{\frac{p-2}{p}+1}}\alpha$$

$$-\frac{(p-1)|\bar{t}_n|^{p-2}\left(\sum_{i=1}^{k}|\bar{t}_i|^p\right)\lambda_n\bar{t}_n}{\left(\sum_{i=1}^{k}|\bar{t}_i|^p\right)^{\frac{p-2}{p}+1}}+\frac{(p-2)|\bar{t}_n|^{2p-2}\mathrm{sign}(\bar{t}_n)\lambda_n|\bar{t}_n|}{\left(\sum_{i=1}^{k}|\bar{t}_i|^p\right)^{\frac{p-2}{p}+1}}$$

$$=\sum_{1\le j\le k,j\neq n}\left(-\frac{(p-2)|\bar{t}_n+t_n|^{p-1}|\bar{t}_j+t_j|^{p-1}\mathrm{sign}(\bar{t}_n+t_n)\mathrm{sign}(\bar{t}_j+t_j)}{\left(\sum_{i=1}^{k}|\bar{t}_i+t_i|^p+\beta\right)^{\frac{p-2}{p}+1}}+\frac{(p-2)|\bar{t}_n|^{p-1}|\bar{t}_j|^{p-1}\mathrm{sign}(\bar{t}_n)\mathrm{sign}(\bar{t}_j)}{\left(\sum_{i=1}^{k}|\bar{t}_i|^p\right)^{\frac{p-2}{p}+1}}\right)s_j$$

$$+\frac{(p-1)|\bar{t}_n+t_n|^{p-2}\left(\sum_{i=1}^{k}|\bar{t}_i+t_i|^p+\sum_{i=k+1}^{\infty}|t_i|^p\right)-(p-2)|\bar{t}_n+t_n|^{2p-2}}{\left(\sum_{i=1}^{k}|\bar{t}_i+t_i|^p+\beta\right)^{\frac{p-2}{p}+1}}s_n-\frac{\frac{p-2}{p}|\bar{t}_n+t_n|^{p-1}\mathrm{sign}(\bar{t}_n+t_n)}{\left(\sum_{i=1}^{k}|\bar{t}_i+t_i|^p+\beta\right)^{\frac{p-2}{p}+1}}\alpha$$

$$-\frac{(p-1)|\bar{t}_n|^{p-2}\left(\sum_{i=1}^{k}|\bar{t}_i|^p\right)}{\left(\sum_{i=1}^{k}|\bar{t}_i|^p\right)^{\frac{p-2}{p}+1}}s_n+\frac{(p-2)|\bar{t}_n|^{2p-2}}{\left(\sum_{i=1}^{k}|\bar{t}_i|^p\right)^{\frac{p-2}{p}+1}}s_n$$

$$=\sum_{1\le j\le k,j\neq n}\left(-\frac{(p-2)|\bar{t}_n+t_n|^{p-1}|\bar{t}_j+t_j|^{p-1}\mathrm{sign}(\bar{t}_n+t_n)\mathrm{sign}(\bar{t}_j+t_j)}{\left(\sum_{i=1}^{k}|\bar{t}_i+t_i|^p+\beta\right)^{\frac{p-2}{p}+1}}+\frac{(p-2)|\bar{t}_n|^{p-1}|\bar{t}_j|^{p-1}\mathrm{sign}(\bar{t}_n)\mathrm{sign}(\bar{t}_j)}{\left(\sum_{i=1}^{k}|\bar{t}_i|^p\right)^{\frac{p-2}{p}+1}}\right)s_j$$

$$+\frac{(p-1)|\bar{t}_n+t_n|^{p-2}\left(\sum_{i=1}^{k}|\bar{t}_i+t_i|^p+\sum_{i=k+1}^{\infty}|t_i|^p\right)-(p-1)|\bar{t}_n|^{p-2}\left(\sum_{i=1}^{k}|\bar{t}_i|^p\right)-(p-2)|\bar{t}_n+t_n|^{2p-2}+(p-2)|\bar{t}_n|^{2p-2}}{\left(\sum_{i=1}^{k}|\bar{t}_i+t_i|^p+\beta\right)^{\frac{p-2}{p}+1}}s_n$$

$$-\frac{\frac{p-2}{p}|\bar{t}_n+t_n|^{p-1}\mathrm{sign}(\bar{t}_n+t_n)}{\left(\sum_{i=1}^{k}|\bar{t}_i+t_i|^p+\beta\right)^{\frac{p-2}{p}+1}}\alpha. \tag{4.8}$$

By (4.8), for $1\le n\le k$, the *n*th-entry of the limit $\lim\limits_{\substack{v\to\theta\\ l_p}}\frac{J(\bar{x}+v)-J(\bar{x})-J'(\bar{x})(v)}{\|v\|_p}$ satisfies

$$\lim_{\substack{v\to\theta\\ l_p}} \frac{\left| \frac{|\bar{t}_n+s_n|^{p-1}\operatorname{sign}(\bar{t}_n+s_n)}{\left(\sum_{i=1}^{k}|\bar{t}_i+s_i|^p+\sum_{i=k+1}^{\infty}|s_i|^p\right)^{\frac{p-2}{p}}} - \frac{|\bar{t}_n|^{p-1}\operatorname{sign}(\bar{t}_n)}{\left(\sum_{i=1}^{k}|\bar{t}_i|^p\right)^{\frac{p-2}{p}}} - \frac{(p-1)\lambda_n|\bar{t}_n|^{p-1}\operatorname{sign}(\bar{t}_n)}{\left(\sum_{i=1}^{k}|\bar{t}_i|^p\right)^{\frac{p-2}{p}}} + \frac{(p-2)|\bar{t}_n|^{p-1}\operatorname{sign}(\bar{t}_n)\sum_{i=1}^{k}\lambda_i|\bar{t}_i|^p}{\left(\sum_{i=1}^{k}|\bar{t}_i|^p\right)^{\frac{p-2}{p}+1}} \right|}{\|v\|_p}$$

$$\leq \lim_{\substack{v\to\theta\\ l_p}} \left( \frac{\left| \sum_{1\leq j\leq k, j\neq n}\left( -\frac{(p-2)|\bar{t}_n+t_n|^{p-1}\left|\bar{t}_j+t_j\right|^{p-1}\operatorname{sign}(\bar{t}_n+t_n)\operatorname{sign}\left(\bar{t}_j+t_j\right)}{\left(\sum_{i=1}^{k}|\bar{t}_i+t_i|^p+\beta\right)^{\frac{p-2}{p}+1}} + \frac{(p-2)|\bar{t}_n|^{p-1}\left|\bar{t}_j\right|^{p-1}\operatorname{sign}(\bar{t}_n)\operatorname{sign}\left(\bar{t}_j\right)}{\left(\sum_{i=1}^{k}|\bar{t}_i|^p\right)^{\frac{p-2}{p}+1}} \right) s_j \right|}{\|v\|_p} \right.$$

$$+ \frac{\left| \frac{(p-1)|\bar{t}_n+t_n|^{p-2}\left(\sum_{i=1}^{k}|\bar{t}_i+t_i|^p+\sum_{i=k+1}^{\infty}|t_i|^p\right)-(p-1)|\bar{t}_n|^{p-2}\left(\sum_{i=1}^{k}|\bar{t}_i|^p\right) - (p-2)|\bar{t}_n+t_n|^{2p-2}+(p-2)|\bar{t}_n|^{2p-2}}{\left(\sum_{i=1}^{k}|\bar{t}_i+t_i|^p+\beta\right)^{\frac{p-2}{p}+1}} s_n \right|}{\|v\|_p}$$

$$\left. + \frac{\left| \frac{\frac{p-2}{p}|\bar{t}_n+t_n|^{p-1}\operatorname{sign}(\bar{t}_n+t_n)}{\left(\sum_{i=1}^{k}|\bar{t}_i+t_i|^p+\beta\right)^{\frac{p-2}{p}+1}} \sum_{i=k+1}^{\infty}|t_i|^p \right|}{\|v\|_p} \right)$$

$$\leq \lim_{\substack{v\to\theta\\ l_p}} \left( \left| \sum_{1\leq j\leq k, j\neq n}\left( -\frac{(p-2)|\bar{t}_n+t_n|^{p-1}\left|\bar{t}_j+t_j\right|^{p-1}\operatorname{sign}(\bar{t}_n+t_n)\operatorname{sign}\left(\bar{t}_j+t_j\right)}{\left(\sum_{i=1}^{k}|\bar{t}_i+t_i|^p+\beta\right)^{\frac{p-2}{p}+1}} + \frac{(p-2)|\bar{t}_n|^{p-1}\left|\bar{t}_j\right|^{p-1}\operatorname{sign}(\bar{t}_n)\operatorname{sign}\left(\bar{t}_j\right)}{\left(\sum_{i=1}^{k}|\bar{t}_i|^p\right)^{\frac{p-2}{p}+1}} \right) \right| \frac{|s_j|}{\left(\sum_{i=k+1}^{\infty}|t_i|^p\right)^{\frac{1}{p}}} \right.$$

$$+ \left| \frac{(p-1)|\bar{t}_n+t_n|^{p-2}\left(\sum_{i=1}^{k}|\bar{t}_i+t_i|^p+\sum_{i=k+1}^{\infty}|t_i|^p\right)-(p-1)|\bar{t}_n|^{p-2}\left(\sum_{i=1}^{k}|\bar{t}_i|^p\right) - (p-2)|\bar{t}_n+t_n|^{2p-2}+(p-2)|\bar{t}_n|^{2p-2}}{\left(\sum_{i=1}^{k}|\bar{t}_i+t_i|^p+\beta\right)^{\frac{p-2}{p}+1}} \right| \frac{|s_n|}{\left(\sum_{i=k+1}^{\infty}|t_i|^p\right)^{\frac{1}{p}}}$$

$$\left. + \left| \frac{\frac{p-2}{p}|\bar{t}_n+t_n|^{p-1}\operatorname{sign}(\bar{t}_n+t_n)}{\left(\sum_{i=1}^{k}|\bar{t}_i+t_i|^p+\beta\right)^{\frac{p-2}{p}+1}} \right| \frac{\sum_{i=k+1}^{\infty}|t_i|^p}{\left(\sum_{i=k+1}^{\infty}|t_i|^p\right)^{\frac{1}{p}}} \right)$$

$$\leq \lim_{\substack{v\to\theta\\ l_p}} \left( \left| \sum_{1\leq j\leq k, j\neq n}\left( -\frac{(p-2)|\bar{t}_n+t_n|^{p-1}\left|\bar{t}_j+t_j\right|^{p-1}\operatorname{sign}(\bar{t}_n+t_n)\operatorname{sign}\left(\bar{t}_j+t_j\right)}{\left(\sum_{i=1}^{k}|\bar{t}_i+t_i|^p+\beta\right)^{\frac{p-2}{p}+1}} + \frac{(p-2)|\bar{t}_n|^{p-1}\left|\bar{t}_j\right|^{p-1}\operatorname{sign}(\bar{t}_n)\operatorname{sign}\left(\bar{t}_j\right)}{\left(\sum_{i=1}^{k}|\bar{t}_i|^p\right)^{\frac{p-2}{p}+1}} \right) \right| \right.$$

$$+ \left| \frac{(p-1)|\bar{t}_n+t_n|^{p-2}\left(\sum_{i=1}^{k}|\bar{t}_i+t_i|^p+\sum_{i=k+1}^{\infty}|t_i|^p\right)-(p-1)|\bar{t}_n|^{p-2}\left(\sum_{i=1}^{k}|\bar{t}_i|^p\right) - (p-2)|\bar{t}_n+t_n|^{2p-2}+(p-2)|\bar{t}_n|^{2p-2}}{\left(\sum_{i=1}^{k}|\bar{t}_i+t_i|^p+\beta\right)^{\frac{p-2}{p}+1}} \right|$$

$$\left. + \left| \frac{\frac{p-2}{p}|\bar{t}_n+t_n|^{p-1}\operatorname{sign}(\bar{t}_n+t_n)}{\left(\sum_{i=1}^{k}|\bar{t}_i+t_i|^p+\beta\right)^{\frac{p-2}{p}+1}} \right| \left(\sum_{i=k+1}^{\infty}|t_i|^p\right)^{1-\frac{1}{p}} \right)$$

$$= \left| \sum_{1\leq j\leq k, j\neq n}\left( -\frac{(p-2)|\bar{t}_n|^{p-1}\left|\bar{t}_j\right|^{p-1}\operatorname{sign}(\bar{t}_n)\operatorname{sign}\left(\bar{t}_j\right)}{\left(\sum_{i=1}^{k}|\bar{t}_i|^p\right)^{\frac{p-2}{p}+1}} + \frac{(p-2)|\bar{t}_n|^{p-1}\left|\bar{t}_j\right|^{p-1}\operatorname{sign}(\bar{t}_n)\operatorname{sign}\left(\bar{t}_j\right)}{\left(\sum_{i=1}^{k}|\bar{t}_i|^p\right)^{\frac{p-2}{p}+1}} \right) \right|$$

$$+\left|\frac{(p-1)|\bar{t}_n|^{p-2}\left(\sum_{i=1}^{k}|\bar{t}_i|^p\right)-(p-1)|\bar{t}_n|^{p-2}\left(\sum_{i=1}^{k}|\bar{t}_i|^p\right)-(p-2)|\bar{t}_n|^{2p-2}+(p-2)|\bar{t}_n|^{2p-2}}{\left(\sum_{i=1}^{k}|\bar{t}_i|^p\right)^{\frac{p-2}{p}+1}}\right|$$

$$+\left|\frac{\frac{p-2}{p}|\bar{t}_n|^{p-1}\mathrm{sign}(\bar{t}_n)}{\left(\sum_{i=1}^{k}|\bar{t}_i|^p\right)^{\frac{p-2}{p}+1}}\right|0\Bigg)$$

$$= 0. \tag{4.9}$$

For all $n > k$, the $n$th-entry of the limit $\lim_{\substack{v\to\theta\\ l_p}}\frac{J(\bar{x}+v)-J(\bar{x})-J'(\bar{x})(v)}{\|v\|_p}$ satisfies that

$$\lim_{\substack{v\to\theta\\ l_p}}\frac{\left\|\left\{0,\ \dots\ 0,\ \frac{|s_{k+1}|^{p-1}\mathrm{sign}(s_{k+1})}{\left(\sum_{i=1}^{\infty}|\bar{t}_i+s_i|^p\right)^{\frac{p-2}{p}}},\ \frac{|s_{k+2}|^{p-1}\mathrm{sign}(s_{k+2})}{\left(\sum_{i=1}^{\infty}|\bar{t}_i+s_i|^p\right)^{\frac{p-2}{p}}},\dots\right\}\right\|_q}{\|v\|_p}$$

$$=\lim_{\substack{v\to\theta\\ l_p}}\frac{\left(\sum_{n=k+1}^{\infty}\left|\frac{|s_n|^{p-1}\mathrm{sign}(s_n)}{\left(\sum_{i=1}^{\infty}|\bar{t}_i+s_i|^p\right)^{\frac{p-2}{p}}}\right|^q\right)^{\frac{1}{q}}}{\|v\|_p}=\lim_{\substack{v\to\theta\\ l_p}}\frac{\frac{\left(\sum_{n=k+1}^{\infty}|s_n|^{(p-1)q}\right)^{\frac{1}{q}}}{\left(\sum_{i=1}^{\infty}|\bar{t}_i+s_i|^p\right)^{\frac{p-2}{p}}}}{\|v\|_p}$$

$$\leq\lim_{\substack{v\to\theta\\ l_p}}\frac{\frac{\left(\sum_{n=1}^{\infty}|s_n|^p\right)^{\frac{1}{q}}}{\left(\sum_{i=1}^{\infty}|\bar{t}_i+s_i|^p\right)^{\frac{p-2}{p}}}}{\|v\|_p}=\lim_{\substack{v\to\theta\\ l_p}}\frac{\|v\|_p^{\frac{p}{q}-1}}{\left(\sum_{i=1}^{\infty}|\bar{t}_i+s_i|^p\right)^{\frac{p-2}{p}}}=\frac{0}{\left(\sum_{i=1}^{\infty}|\bar{t}_i|^p\right)^{\frac{p-2}{p}}}$$

$= 0.$ Since $2 < p < \infty$, it implies that $\frac{p}{q}-1>0$. (4.10)

This proves case 2. By (4.9) and (4.10), we obtain that

$$\lim_{\substack{v\to\theta\\ l_p}}\frac{\|J(\bar{x}+v)-J(\bar{x})-J'(\bar{x})(v)\|_q}{\|v\|_p}=0. \qquad \square$$

## 5. The Mordukhovich Differentiability and Covering Constant for $J$ in $l_p$

Recall that $p$ and $q$ are positive numbers satisfying $1 < p, q < \infty$ and $\frac{1}{p}+\frac{1}{q}=1$. Let $\bar{x} \in l_p$ and $\bar{y} \in l_q$ with $\bar{y} = J(\bar{x})$. Suppose that $J$ is Fréchet differentiable at $\bar{x}$ with its Fréchet derivative $\nabla J(\bar{x})\colon l_p \to l_q$, which is a continuous and linear mapping. Then, by Theorem 1.38 in [28], the Mordukhovich derivative of $J$ at $(\bar{x}, \bar{y})$ exists and

$$\widehat{D}^*J(\bar{x},\bar{y})(y) = (\nabla J(\bar{x}))^*(y),\ \text{ for all } y \in l_p. \tag{5.1}$$

Here $(\nabla J(x))^*\colon (l_q)^* \to (l_p)^*$ is the conjugate operator of the operator $\nabla J(\bar{x})\colon l_p \to l_q$. However, by $(l_q)^* = l_p$ and $(l_p)^* = l_q$ again, in this special case, we have

$$\widehat{D}^*J(\bar{x},\bar{y}) = (\nabla J(x))^*: l_p \to l_q. \tag{5.2}$$

In Theorem 4.4, we only proved the Fréchet differentiability of $J$ at $\bar{x}$ with $|\mathcal{N}(\bar{x})| < \infty$. In this case, in Theorem 4.4, we find the explicit solution of $\nabla J(\bar{x})$. So, in order to find Mordukhovich derivative of $J$ at $(\bar{x},\bar{y})$ by using the connection between Fréchet and Mordukhovich derivatives in Banach spaces, we will only consider a special case that $|\mathcal{N}(\bar{x})| < \infty$.

**Theorem. 5.1**. *Let* $2 < p < \infty$. *Let* $\bar{x} = (\bar{t}_1, \bar{t}_2, \dots) \in l_p\backslash\{\theta\}$ *with* $|\mathcal{N}(\bar{x})| < \infty$, *and* $\bar{y} \in l_q$ *with* $\bar{y} = J(\bar{x})$. *Then J is Mordukhovich differentiable at* $(\bar{x},\bar{y})$ *and for any* $v = \{s_n\}_{n=1}^{\infty} \in l_p\backslash\{\theta\}$ *with ratio sequence* $\{\lambda_n\}_{n=1}^{\infty}$ *with respect to* $\bar{x}$, *one has*

$$\widehat{D}^*J(\bar{x},\bar{y})(v) = \nabla J(\bar{x})(v) = J(\bar{x})\left((p-1)\begin{pmatrix}\lambda_1 & 0 & \dots\\ 0 & \lambda_2 & 0\\ \vdots & 0 & \ddots\end{pmatrix} - \frac{(p-2)\sum_{i=1}^{\infty}\lambda_i|\bar{t}_i|^p}{\sum_{i=1}^{\infty}|\bar{t}_i|^p}\begin{pmatrix}1 & 0 & \dots\\ 0 & 1 & 0\\ \vdots & 0 & \ddots\end{pmatrix}\right). \tag{5.3}$$

*Proof*. By (5.1) and (5.2), for $2 < p < \infty$, in the Banach space $l_p$, the Mordukhovich derivative $\widehat{D}^*J(\bar{x},\bar{y})$: $l_q^* \to l_p^*$ is actually a mapping $\widehat{D}^*J(\bar{x},\bar{y})$: $l_p \to l_q$ satisfying

$$\widehat{D}^*J(\bar{x},\bar{y})(y^*) = (\nabla J(\bar{x}))^*(y^*), \text{ for any } y^* \in l_p.$$

Recall that $(\nabla J(x))^*: (l_q)^* \to (l_p)^*$ is the conjugate operator of the operator $\nabla J(\bar{x})$: $l_p \to l_q$. In Theorem 4.4, it shows that $\nabla J(\bar{x})$ is defined by a main diagonal $\infty \times \infty$ matrix. It implies that $(\nabla J(x))^* = \nabla J(x)$. Then, (5.3) is proved by (5.2) and Theorem 4.4. □

Let $\bar{x} \in l_p$ and $\bar{y} \in l_q$ with $\bar{y} = J(\bar{x})$. With respect to the explicit solution of $\widehat{D}^*J(\bar{x},\bar{y})$ given in (5.3), we can calculate the covering constant for $J$ at point $(\bar{x},\bar{y}) \in l_p \times l_q$. Once again, we only consider a special case that $|\mathcal{N}(\bar{x})| < \infty$. We review the concept of mapping $J$ in the uniformly convex and uniformly smooth Banach space $l_p$. Let $\bar{x} \in l_p$ and $\bar{y} \in l_q$ with $\bar{y} = J(\bar{x})$. By (2.6), the covering constant for $J$ at point $(\bar{x}, \bar{y})$ is given by

$$\hat{\alpha}(J,\bar{x},\bar{y}) = \sup_{\eta>0}\inf\left\{\|z^*\|_q: z^* \in \widehat{D}^*J(x,y)(w^*), x \in \mathbb{B}_{l_p}(\bar{x},\eta), y = J(x) \in \mathbb{B}_{l_q}(\bar{y},\eta), \|w^*\|_p = 1\right\}. \tag{5.4}$$

Here $\mathbb{B}_{l_p}(\bar{x},\eta)$ is the closed ball in $l_p$ centered at $\bar{x}$ with radius $\eta$ and $\mathbb{B}_{l_q}(\bar{y},\eta)$ is the closed ball in $l_q$ centered at $\bar{y}$ with radius $\eta$.

**Theorem. 5.2**. *Let* $2 < p < \infty$. *Let* $\bar{x} = (\bar{t}_1, \bar{t}_2, \dots) \in l_p\backslash\{\theta\}$ *with* $|\mathcal{N}(\bar{x})| < \infty$, *and* $\bar{y} \in l_q$ *with* $\bar{y} = J(\bar{x})$. *The covering constant* $\hat{\alpha}(J,\bar{x},\bar{y})$ *for J at point* $(\bar{x}, \bar{y})$ *satisfies that*

$$\hat{\alpha}(J,\bar{x},\bar{y}) = 0.$$

*Proof*. Let $\bar{x} = (\bar{t}_1, \bar{t}_2, \dots) \in l_p\backslash\{\theta\}$ with $|\mathcal{N}(\bar{x})| < \infty$. By $|\mathcal{N}(\bar{x})| < \infty$, the sequence (5.3) that is the solution of $\widehat{D}^*J(\bar{x},\bar{y})$ has finite number of nonzero entries. There is a positive integer $N$ which depends on $\bar{x}$, such that $\bar{t}_n = 0$, for each $n > N$. This implies that the $n$th entry of $J(\bar{x})$ is zero. And therefore, for each $n > N$,

$$\text{the } n\text{th entry of } \widehat{D}^*J(\bar{x},\bar{y})(v) \text{ is zero, for any } v = \{s_n\}_{n=1}^{\infty} \in l_p\backslash\{\theta\}. \tag{5.5}$$

For any positive integer $n$, let $u_n \in l_p$ that the $n$th entry is 1 and all other entries are 0. By (5.4) and (5.5), we estimate

$$\hat{\alpha}(J,\bar{x},\bar{y}) = \sup_{\eta>0}\inf\left\{\|z^*\|_q : z^* \in \widehat{D}^*J(x,y)(w^*), x \in \mathbb{B}_{l_p}(\bar{x},\eta), y = J(x) \in \mathbb{B}_{l_q}(\bar{y},\eta), \|w^*\|_p = 1\right\}$$

$$= \sup_{\eta>0}\inf\{\|z^*\|_q : z^* = \widehat{D}^*J(\bar{x},\bar{y})(u_n), n = 1,2,\dots\}$$

$$= 0.$$ □

## 6. Differentiation of A Duality Mapping $\mathcal{J}$ in $l_p$

### 6.1. Gâteaux and Fréchet Derivatives of the Duality Mapping $\mathcal{J}$ in $l_p$

Duality mappings has been widely applied to nonlinear analysis and analysis in Banach spaces (See [6]). In this section and the following section, we respectively define a general duality mapping $\mathcal{J}$ in $l_p$. We investigate its Gâteaux, Fréchet and Mordukhovich differentiability in $l_p$, for $3 < p < \infty$. For given point $\bar{x} \in l_p$, after the differentiability of $\mathcal{J}$ is proved, we find the explicit Gâteaux, Fréchet and Mordukhovich derivatives of the considered duality mapping.

For $1 < p < \infty$ and $\frac{1}{p} + \frac{1}{q} = 1$, a mapping $\mathcal{J}: l_p \to l_q$ is defined, for any $x = (t_1, t_2, \dots) \in l_p$, by

$$\mathcal{J}(x) = (|t_1|^{p-1}\mathrm{sign}(t_1),\ |t_2|^{p-1}\mathrm{sign}(t_2), \dots). \tag{6.1}$$

One can check that the duality mapping $\mathcal{J}: l_p \to l_q$ satisfies that $\mathcal{J}(\theta) = \theta$ and

$$\langle \mathcal{J}(x),\ x\rangle = \|J(x)\|_q^q = \|x\|_p^p, \text{ for any } x \in l_p. \tag{6.2}$$

By definition (2.1), $\mathcal{J}$ is a duality mapping in $l_p$ with gauge function $\varphi(t) = t^{\frac{p}{q}} = t^{p-1}$, for $t \in [0,\infty)$.

**Theorem 6.1**. *Let* $3 < p < \infty$. $\mathcal{J}$ *is Fréchet differentiable on* $l_p$ *and, for any* $\bar{x} = (\bar{t}_1, \bar{t}_2, \dots) \in l_p$, *we have*

$$\nabla\mathcal{J}(\bar{x}) = (p-1)\begin{pmatrix} |\bar{t}_1|^{p-2} & 0 & \dots \\ 0 & |\bar{t}_2|^{p-2} & 0 \\ \vdots & 0 & \ddots \end{pmatrix}.$$

*More precisely,* $\nabla\mathcal{J}(\bar{x}): l_p \to l_q$ *is a pointwise multiplication operator on* $l_p$ *and, for* $v = \{t_n\}_{n=1}^{\infty} \in l_p$,

$$\nabla\mathcal{J}(\bar{x})(v) = \{t_n\}_{n=1}^{\infty}(p-1)\begin{pmatrix} |\bar{t}_1|^{p-2} & 0 & \dots \\ 0 & |\bar{t}_2|^{p-2} & 0 \\ \vdots & 0 & \ddots \end{pmatrix} = \{(p-1)|\bar{t}_n|^{p-2}t_n\}_{n=1}^{\infty} \in l_q, \tag{6.3}$$

*and*

$$\|\nabla\mathcal{J}(\bar{x})(v)\|_q \le (p-1)\|\bar{x}\|_p^{p-2}\ \|v\|_p. \tag{6.4}$$

*Proof*. Let arbitrarily given $\bar{x} = (\bar{t}_1, \bar{t}_2, \dots) \in l_p$. For each $n = 1, 2, \dots$, define a function $f_n$ as follows:

$$f_n(t) = |\bar{t}_n + t|^{p-1}\mathrm{sign}(\bar{t}_n + t), \text{ for } t \in \mathbb{R}.$$

Taking the first derivative of the function $f_n$ with respect to $t$, we have

$$f_n'(t) = (p-1)|\bar{t}_n + t|^{p-2},$$

and $$f_n''(t) = (p-1)(p-2)|\bar{t}_n + t|^{p-3}\operatorname{sign}(\bar{t}_n + t), \text{ for } t \in \mathbb{R}. \tag{6.5}$$

By $3 < p < \infty$ and $\frac{1}{p} + \frac{1}{q} = 1$, for $v = \{t_n\}_{n=1}^{\infty} \in l_p \backslash \{\theta\}$, by Taylor's Approximating Theorem with order 2, by (6.1), (6.2) and (6.5), we estimate the following limit.

$$\lim_{v\to\theta} \frac{\left\| \{|\bar{t}_n + t_n|^{p-1}\operatorname{sign}(\bar{t}_n + t_n) - |\bar{t}_n|^{p-1}\operatorname{sign}(\bar{t}_n) - (p-1)|\bar{t}_n|^{p-2}t_n\}_{n=1}^{\infty} \right\|_q}{\|v\|_p}$$

$$= \lim_{v\to\theta} \frac{\left\| \{f_n(t_n) - f_n(0) - f_n'(0)t_n\}_{n=1}^{\infty} \right\|_q}{\|v\|_p}$$

$$= \lim_{v\to\theta} \frac{\left\| \left\{\frac{f_n''(s_n)}{2!}t_n^2\right\}_{n=1}^{\infty} \right\|_q}{\|v\|_p} \qquad \text{(For some } s_n \in \mathbb{R} \text{ with } |s_n| \le |t_n|)$$

$$= \frac{(p-1)(p-2)}{2} \lim_{v\to\theta} \frac{\left\| \{|\bar{t}_n + s_n|^{p-3}\operatorname{sign}(\bar{t}_n + s_n)t_n^2\}_{n=1}^{\infty} \right\|_q}{\|v\|_p}$$

$$= \frac{(p-1)(p-2)}{2} \lim_{v\to\theta} \frac{\left(\sum_{n=1}^{\infty} |\bar{t}_n + s_n|^{(p-3)q}|t_n|^{2q}\right)^{\frac{1}{q}}}{\|v\|_p}$$

$$\le \frac{(p-1)(p-2)}{2} \lim_{v\to\theta} \frac{\left( \left( \sum_{n=1}^{\infty} \left(|\bar{t}_n + s_n|^{(p-3)q}\right)^{\frac{p}{p-2q}} \right)^{\frac{p-2q}{p}} \left( \sum_{n=1}^{\infty} \left(|t_n|^{2q}\right)^{\frac{p}{2q}} \right)^{\frac{2q}{p}} \right)^{\frac{1}{q}}}{\|v\|_p} \qquad \left(\frac{p}{2q} > 1\right)$$

$$= \frac{(p-1)(p-2)}{2} \lim_{v\to\theta} \frac{\left( \left( \sum_{n=1}^{\infty} |\bar{t}_n + s_n|^{p} \right)^{\frac{p-2q}{p}} \left( \sum_{n=1}^{\infty} \left(|t_n|^{2q}\right)^{\frac{p}{2q}} \right)^{\frac{2q}{p}} \right)^{\frac{1}{q}}}{\|v\|_p} \qquad \left(\frac{(p-3)q}{p-2q} = 1\right)$$

$$= \frac{(p-1)(p-2)}{2} \lim_{v\to\theta} \frac{\left( \left( \sum_{n=1}^{\infty} |\bar{t}_n + s_n|^{p} \right)^{\frac{p-2q}{p}} \left( \sum_{n=1}^{\infty} |t_n|^{p} \right)^{\frac{2q}{p}} \right)^{\frac{1}{q}}}{\|v\|_p}$$

$$= \frac{(p-1)(p-2)}{2} \lim_{v\to\theta} \frac{\left( \left( \sum_{n=1}^{\infty} |\bar{t}_n + s_n|^{p} \right)^{\frac{p-2q}{p}} \|v\|_p^{2q} \right)^{\frac{1}{q}}}{\|v\|_p}$$

$$= \frac{(p-1)(p-2)}{2} \lim_{v\to\theta} \frac{\left( \left( \sum_{n=1}^{\infty} |\bar{t}_n + s_n|^{p} \right)^{\frac{p-2q}{p}} \right)^{\frac{1}{q}} \|v\|_p^{2}}{\|v\|_p}$$

$$= \frac{(p-1)(p-2)}{2} \lim_{v\to\theta} \left( \left( \sum_{n=1}^{\infty} |\bar{t}_n + s_n|^{p} \right)^{\frac{p-2q}{p}} \right)^{\frac{1}{q}} \|v\|_p$$

$$= \frac{(p-1)(p-2)}{2} \lim_{v\to\theta} \left( (\sum_{n=1}^{\infty} |\bar{t}_n + s_n|^p)^{\frac{1}{p}} \right)^{\frac{p-2q}{q}} \|v\|_p$$

$$\leq \frac{(p-1)(p-2)}{2} \lim_{v\to\theta} \left( (\sum_{n=1}^{\infty} (|\bar{t}_n|)^p)^{\frac{1}{p}} + (\sum_{n=1}^{\infty} (|s_n|)^p)^{\frac{1}{p}} \right)^{\frac{p-2q}{q}} \|v\|_p \qquad \left(\frac{p-2q}{q} > 0\right)$$

$$\leq \frac{(p-1)(p-2)}{2} \lim_{v\to\theta} \left( (\sum_{n=1}^{\infty} (|\bar{t}_n|)^p)^{\frac{1}{p}} + (\sum_{n=1}^{\infty} (|t_n|)^p)^{\frac{1}{p}} \right)^{\frac{p-2q}{q}} \|v\|_p \qquad (|s_n| \leq |t_n|)$$

$$= \frac{(p-1)(p-2)}{2} \lim_{v\to\theta} \left(\|\bar{x}\|_p + \|v\|_p\right)^{\frac{p-2q}{q}} \|v\|_p$$

$$= 0.$$

This proves (6.3). By (6.3), $v = \{t_n\}_{n=1}^{\infty} \in l_p$, by 3 < $p$ < ∞, we estimate

$$\|\nabla\mathcal{J}(\bar{x})(v)\|_q = \|\{(p-1)|\bar{t}_n|^{p-2} t_n\}_{n=1}^{\infty}\|_q$$

$$= (p-1)\left(\sum_{n=1}^{\infty} |\bar{t}_n|^{(p-2)q} |t_n|^q\right)^{\frac{1}{q}}$$

$$\leq (p-1)\left( \left(\sum_{n=1}^{\infty} |\bar{t}_n|^{(p-2)q\frac{p}{p-q}}\right)^{\frac{p-q}{p}} \left(\sum_{n=1}^{\infty} |t_n|^{q\frac{p}{q}}\right)^{\frac{q}{p}} \right)^{\frac{1}{q}}$$

$$= (p-1)\left( (\sum_{n=1}^{\infty} |\bar{t}_n|^p)^{\frac{p-q}{p}} (\sum_{n=1}^{\infty} |t_n|^p)^{\frac{q}{p}} \right)^{\frac{1}{q}}$$

$$= (p-1)\left(\|\bar{x}\|_p^{p-q} \|v\|_p^q\right)^{\frac{1}{q}}$$

$$= (p-1)\|\bar{x}\|_p^{p-2} \|v\|_p.$$

This proves (6.4). □

By the connection between Gâteaux and Fréchet derivatives in Banach spaces, we have the following corollary of Theorem 6.1 immediately.

**Corollary 6.2**. *Let* 3 < $p$ < ∞. $\mathcal{J}$ *is Gâteaux differentiable in* $l_p$ *and, for any* $\bar{x} = (\bar{t}_1, \bar{t}_2, \dots) \in l_p$, *the Gâteaux derivative* $\mathcal{J}'(\bar{x})$ *satisfies that*

$$\mathcal{J}'(\bar{x})(v) = \nabla\mathcal{J}(\bar{x})(v) = \{(p-1)|\bar{t}_n|^{p-2} t_n\}_{n=1}^{\infty} \in l_q, \text{ for any } v = \{t_n\}_{n=1}^{\infty} \in l_p \backslash \{\theta\}.$$

### 6.2. Mordukhovich Derivative and Covering Constant for the Duality Mapping $\mathcal{J}$ in $l_p$

In Theorem 6.1, the explicit Fréchet derivative of the duality mapping $\mathcal{J}: l_p \to l_q$ is proved for the case 3 < $p$ < ∞. In this subsection, we first apply Theorem 6.1 and Theorem 1.38 in [28] to find the explicit

Mordukhovich derivative of the mapping $\mathcal{J}: l_p \to l_q$, by which, we calculate the covering constant for $\mathcal{J}$.

**Theorem 6.3**. *Let* $3 < p < \infty$. $\mathcal{J}: l_p \to l_q$ *is Mordukhovich differentiable in* $l_p$ *and for* $\bar{x} = \{\bar{t}_n\}_{n=1}^{\infty} \in l_p$, *and the Mordukhovich derivative* $\widehat{D}^*\mathcal{J}(\bar{x}): l_p \to l_q$ *is equal to the Fréchet derivative* $\nabla J(\bar{x}): l_p \to l_q$ *with the following representation.*

$$\widehat{D}^*\mathcal{J}(\bar{x}) = \nabla\mathcal{J}(\bar{x}) = (p-1)\begin{pmatrix} |\bar{t}_1|^{p-2} & 0 & \dots \\ 0 & |\bar{t}_2|^{p-2} & 0 \\ \vdots & 0 & \ddots \end{pmatrix}.$$

*More precisely,* $\widehat{D}^*\mathcal{J}(\bar{x}) = \nabla\mathcal{J}(\bar{x}): l_p \to l_q$ *is a pointwise multiplication operator in* $l_p$ *and, for any* $v = \{t_n\}_{n=1}^{\infty} \in l_p$,

$$\widehat{D}^*\mathcal{J}(\bar{x})(v) = \nabla\mathcal{J}(\bar{x})(v) = \{t_n\}_{n=1}^{\infty}(p-1)\begin{pmatrix} |\bar{t}_1|^{p-2} & 0 & \dots \\ 0 & |\bar{t}_2|^{p-2} & 0 \\ \vdots & 0 & \ddots \end{pmatrix} = \{(p-1)|\bar{t}_n|^{p-2}t_n\}_{n=1}^{\infty} \in l_q.$$

*Furthermore, for any* $v = \{t_n\}_{n=1}^{\infty} \in l_p$ *and* $y = \{s_n\}_{n=1}^{\infty} \in l_p$ *, we have*

$$\langle \widehat{D}^*\mathcal{J}(\bar{x})(y), v\rangle = \langle y, \nabla\mathcal{J}(\bar{x})(v)\rangle = (p-1)\textstyle\sum_{n=1}^{\infty}|\bar{t}_n|^{p-2}s_n t_n.$$

*Proof*. Since $\mathcal{J}: l_p \to l_q$ and $\nabla\mathcal{J}(\bar{x}): l_p \to l_q$, by definition $\widehat{D}^*\mathcal{J}(\bar{x}): (l_q)^* \to (l_p)^*$ is a mapping from $(l_q)^*$ to $(l_p)^*$. By $(l_q)^* = l_p$ and $(l_p)^* = l_q$, it induces that $\widehat{D}^*\mathcal{J}(\bar{x}): l_p \to l_q$ is a mapping from $l_p$ to $l_q$. By Theorem 1.38 in [18], for any $y \in (l_q)^* = l_p$, we have

$$\widehat{D}^*\mathcal{J}(\bar{x})(y) = (\nabla\mathcal{J}(\bar{x}))^*(y) \in (l_p)^* = l_q.$$

Since $\nabla\mathcal{J}(\bar{x})$ is represented by a main diagonal $\infty \times \infty$ matrix, it satisfies

$$(\nabla\mathcal{J}(\bar{x}))^* = (\nabla\mathcal{J}(\bar{x}))^T = \nabla\mathcal{J}(\bar{x}).$$

This proves that $\widehat{D}^*\mathcal{J}(\bar{x}) = \nabla\mathcal{J}(\bar{x})$. For any $y = \{s_n\}_{n=1}^{\infty} \in l_p$, $v = \{t_n\}_{n=1}^{\infty} \in l_p$, we have

$$\nabla\mathcal{J}(\bar{x})(v) = \{(p-1)|\bar{t}_n|^{p-2}t_n\}_{n=1}^{\infty} \in l_q \quad \text{and} \quad \widehat{D}^*\mathcal{J}(\bar{x})(y) = \{(p-1)|\bar{t}_n|^{p-2}s_n\}_{n=1}^{\infty} \in l_q.$$

Let $\langle\cdot,\ \cdot\rangle$ denote the real canonical pairing between $l_q$ and $l_p$. This implies that

$$\langle \widehat{D}^*\mathcal{J}(\bar{x})(y), v\rangle = \langle y, \nabla\mathcal{J}(\bar{x})(v)\rangle = (p-1)\textstyle\sum_{n=1}^{\infty}|\bar{t}_n|^{p-2}s_n t_n. \qquad \square$$

By using the explicit representation of the utility-like mapping $\mathcal{J}: l_p \to l_q$, in next theorem, we calculate its covering constant.

**Theorem 6.4**. *Let* $3 < p < \infty$. *Let* $\bar{x} = \{\bar{t}_n\}_{n=1}^{\infty} \in l_p$ *and* $\bar{y} \in l_q$ *with* $\mathcal{J}(\bar{x}) = \bar{y}$. *Then*

$$\hat{\alpha}(\mathcal{J}, \bar{x}, \bar{y}) = 0.$$

*Proof*. By definition and by Theorems 4.4 and 4.5, we calculate

$$\hat{\alpha}(\mathcal{J}, \bar{x}, \bar{y}) = \sup_{\eta>0}\inf\left\{\|z\|_{l_q} : z = \widehat{D}^*\mathcal{J}(x)(w), x \in \mathbb{B}_{l_p}(\bar{x}, \eta), \mathcal{J}(x) \in \mathbb{B}_{l_q}(\bar{y}, \eta), \|w\|_{l_p} = 1\right\}$$

$$= \sup_{\eta>0} \inf\left\{\|z\|_{l_q}: z = \nabla\mathcal{J}(x)(w), x \in \mathbb{B}_{l_p}(\bar{x},\eta), \mathcal{J}(x) \in \mathbb{B}_{l_q}(\bar{y},\eta), \|w\|_{l_p} = 1\right\}$$

$$= \sup_{\eta>0} \inf\left\{\|\nabla\mathcal{J}(x)(w)\|_{l_q}: x \in \mathbb{B}_{l_p}(\bar{x},\eta), \mathcal{J}(x) \in \mathbb{B}_{l_q}(\bar{y},\eta), \|w\|_{l_p} = 1\right\}$$

$$\leq \sup_{\eta>0} \inf\left\{\|\nabla\mathcal{J}(\bar{x})(w)\|_{l_q}: \|w\|_{l_p} = 1, w = \{s_n\}_{n=1}^{\infty} \in l_p\right\}$$

$$= \sup_{\eta>0} \inf\left\{\|\{(p-1)|\bar{t}_n|^{p-2}s_n\}_{n=1}^{\infty}\|_{l_q}: \|w\|_{l_p} = 1, w = \{s_n\}_{n=1}^{\infty} \in l_p\right\}$$

$$= \sup_{\eta>0} \inf\left\{(\textstyle\sum_{n=1}^{\infty}|(p-1)|\bar{t}_n|^{p-2}s_n|^q)^{\frac{1}{q}}: \|w\|_{l_p} = 1, w = \{s_n\}_{n=1}^{\infty} \in l_p\right\}$$

$$\leq (p-1)\sup_{\eta>0} \inf\{|t_n|^{p-2}: n = 1, 2, \ldots\}$$

$$= 0.$$ □

## 7. Differentiation of A Generalized Duality Mapping $\mathfrak{J}$ in $l_p$

### 7.1. Gâteaux and Fréchet Derivatives of the Generalized Duality Mapping $\mathfrak{J}$ in $l_p$

For $3 < p < \infty$ and $\frac{1}{p} + \frac{1}{q} = 1$, it implies that $\frac{p}{q} = p - 1$. In this section, we define a generalized duality mapping $\mathfrak{J}$: $l_p \to l_q$, which is also a single-valued mapping, for any $x = (t_1, t_2, \ldots) \in l_p$, by

$$\mathfrak{J}(x) = \left\{\frac{|t_n|^{\frac{p}{q}}\text{sign}(t_n)}{1+|t_n|}\right\}_{n=1}^{\infty} = \left\{\frac{|t_n|^{p-1}\text{sign}(t_n)}{1+|t_n|}\right\}_{n=1}^{\infty} \in l_q. \tag{7.1}$$

**Lemma 7.1**. *Let* $3 < p < \infty$. *The duality mapping* $\mathfrak{J}$: $l_p \to l_q$ *satisfies that* $\mathfrak{J}(\theta) = \theta$ *and*

$$\langle\mathfrak{J}(x),\ x\rangle = \|\mathfrak{J}(x)\|_q^q \leq \|x\|_p^p, \text{ for any } x \in l_p. \tag{7.2}$$

*Proof*. By (7.1), we prove (7.2). For any $x = (x_1, x_2, \ldots) \in l_p$, by definition (7.1), we calculate

$$\langle\mathfrak{J}(x),\ x\rangle = \sum_{n=1}^{\infty}\frac{|t_n|^{p-1}\text{sign}(t_n)}{1+|t_n|}t_n = \sum_{n=1}^{\infty}\frac{|t_n|^p}{1+|t_n|} \leq \|x\|_p^p,$$

and
$$\|\mathfrak{J}(x)\|_q^q = \sum_{n=1}^{\infty}\left|\frac{|t_n|^{p-1}\text{sign}(t_n)}{1+|t_n|}\right|^q \leq \sum_{n=1}^{\infty}|t_n|^p = \|x\|_p^p.$$ □

Then, by (2.2), $\mathfrak{J}$: $l_p \to l_q$ is a generalized duality mapping on $l_p$ with gauge function $\varphi(t) = t^{p-1}$, for $t \in [0, \infty)$.

**Theorem 7.2**. *Let* $3 < p < \infty$. $\mathfrak{J}$ *is Fréchet differentiable in* $l_p$ *and, for any* $\bar{x} = (\bar{t}_1, \bar{t}_2, \ldots) \in l_p$,

$$\nabla\mathfrak{J}(\bar{x}) = \begin{pmatrix} \frac{\left((p-1)|\bar{t}_1|^{p-2}+(p-2)|\bar{t}_1|^{p-1}\right)}{(1+|\bar{t}_1|)^2} & 0 & \cdots \\ 0 & \frac{\left((p-1)|\bar{t}_2|^{p-2}+(p-2)|\bar{t}_2|^{p-1}\right)}{(1+|\bar{t}_2|)^2} & 0 \\ \vdots & 0 & \ddots \end{pmatrix}. \tag{7.3}$$

*More precisely,* $\nabla\mathfrak{J}(\bar{x})\colon l_p \to l_q$ *is a pointwise multiplication operator in* $l_p$ *and, for* $v = \{t_n\}_{n=1}^{\infty} \in l_p$*,*

$$\nabla\mathfrak{J}(\bar{x})(v) = \left\{\frac{\left((p-1)|\bar{t}_n|^{p-2}+(p-2)|\bar{t}_n|^{p-1}\right)}{(1+|\bar{t}_n|)^2} t_n\right\}_{n=1}^{\infty} \in l_q, \tag{7.4}$$

*and* $$\|\nabla\mathfrak{J}(\bar{x})(v)\|_p \leq \left((p-1)\|\bar{x}\|_p^{p-2} + (p-2)\|\bar{x}\|_p^{p-1}\right)\|v\|_p.$$

*Proof*. Let arbitrarily given $\bar{x} = (\bar{t}_1, \bar{t}_2, \ldots) \in l_p$. For each $n = 1, 2, \ldots$, define a function $f_n$ as follows:

$$f_n(t) = \frac{|\bar{t}_n+t|^{p-1}\text{sign}(\bar{t}_n+t)}{1+|\bar{t}_n+t|}, \text{ for } t \in \mathbb{R}. \tag{7.5}$$

Taking the first derivative of the function $f_n$ with respect to $t$, we have

$$f_n'(t) = \frac{(p-1)|\bar{t}_n+t|^{p-2}+(p-2)|\bar{t}_n+t|^{p-1}}{(1+|\bar{t}_n+t|)^2}, \tag{7.6}$$

and $$f_n''(t) = \frac{(p-1)(p\text{-}2)|\bar{t}_n+t|^{p-3}+2(p-1)(p\text{-}3)|\bar{t}_n+t|^{p-2}+(p-2)(p\text{-}3)|\bar{t}_n+t|^{p-1}}{(1+|\bar{t}_n+t|)^3}\text{sign}(\bar{t}_n+t), \text{ for } t \in \mathbb{R}. \tag{7.7}$$

Let $3 < p < \infty$ with $\frac{1}{p}+\frac{1}{q} = 1$, For $v = \{t_n\}_{n=1}^{\infty} \in l_p\backslash\{\theta\}$, by Taylor's Approximating Theorem with order 2, by (7.5), (7.6) and (7.7), we estimate the following limit.

$$\lim_{v\to\theta}\frac{\left\|\left\{\frac{|\bar{t}_n+t_n|^{\frac{p}{q}}\text{sign}(\bar{t}_n+t_n)}{1+|\bar{t}_n+t_n|} - \frac{|\bar{t}_n|^{\frac{p}{q}}\text{sign}(\bar{t}_n)}{1+|\bar{t}_n|} - \frac{\left((p-1)|\bar{t}_n|^{p-2}+(p-2)|\bar{t}_n|^{p-1}\right)}{(1+|\bar{t}_n|)^2}t_n\right\}_{n=1}^{\infty}\right\|_q}{\|v\|_p}$$

$$= \lim_{v\to\theta}\frac{\left\|\left\{f_n(t_n)-f_n(0)-f_n'(0)t_n\right\}_{n=1}^{\infty}\right\|_q}{\|v\|_p}$$

$$= \lim_{v\to\theta}\frac{\left\|\left\{\frac{(p-1)(p\text{-}2)|\bar{t}_n+s_n|^{p-3}+2(p-1)(p\text{-}3)|\bar{t}_n+s_n|^{p-2}+(p-2)(p\text{-}3)|\bar{t}_n+s_n|^{p-1}}{(1+|\bar{t}_n+s_n|)^3}sign(\bar{t}_n+t)\,t_n^2\right\}_{n=1}^{\infty}\right\|_q}{\|v\|_p}$$

(For some $s_n \in \mathbb{R}$ with $|s_n| \leq |t_n|$)

$$\leq \lim_{v\to\theta}\frac{(p-1)(p\text{-}2)\left(\sum_{n=1}^{\infty}\frac{|\bar{t}_n+s_n|^{(p-3)q}|t_n|^{2q}}{(1+|\bar{t}_n+s_n|)^{3q}}\right)^{\frac{1}{q}}}{\|v\|_p} + \lim_{v\to\theta}\frac{2(p-1)(p\text{-}3)\left(\sum_{n=1}^{\infty}\frac{|\bar{t}_n+s_n|^{(p-2)q}|t_n|^{2q}}{(1+|\bar{t}_n+s_n|)^{3q}}\right)^{\frac{1}{q}}}{\|v\|_p}$$

$$+ \lim_{v\to\theta} \frac{(p-2)(p\text{-}3)\left(\sum_{n=1}^{\infty}\frac{|\bar{t}_n+s_n|^{(p-1)q}|t_n|^{2q}}{(1+|\bar{t}_n+s_n|)^{3q}}\right)^{\frac{1}{q}}}{\|v\|_p}$$

$$\leq (p-1)(p\text{-}2) \lim_{v\to\theta} \frac{\left(\left(\sum_{n=1}^{\infty}\left(|\bar{t}_n+s_n|^{\left(\frac{p}{q}-2\right)q}\right)^{\frac{p}{p-2q}}\right)^{\frac{p-2q}{p}}\left(\sum_{n=1}^{\infty}\left(|t_n|^{2q}\right)^{\frac{p}{2q}}\right)^{\frac{2q}{p}}\right)^{\frac{1}{q}}}{\|v\|_p} \quad \left(\frac{p}{2q}>1\right)$$

$$+ 2(p-1)(p\text{-}3) \lim_{v\to\theta} \frac{\left(\sum_{n=1}^{\infty}|\bar{t}_n+s_n|^{(p-2)q}|t_n|^{q}\right)^{\frac{1}{q}}\|v\|_p}{\|v\|_p} + (p-2)(p\text{-}3) \lim_{v\to\theta} \frac{\left(\sum_{n=1}^{\infty}|\bar{t}_n+s_n|^{(p-1)q}\right)^{\frac{1}{q}}\|v\|_p^2}{\|v\|_p}$$

$$\leq (p-1)(p\text{-}2) \lim_{v\to\theta} \frac{\left(\left(\sum_{n=1}^{\infty}|\bar{t}_n+s_n|^{p}\right)^{\frac{p-2q}{p}}\left(\sum_{n=1}^{\infty}|t_n|^{p}\right)^{\frac{2q}{p}}\right)^{\frac{1}{q}}}{\|v\|_p}$$

$$+2(p-1)(p\text{-}3) \lim_{v\to\theta} \frac{\left(\left(\sum_{n=1}^{\infty}\left(|\bar{t}_n+s_n|^{(p-2)q}\right)^{\frac{p}{p-q}}\right)^{\frac{p-q}{p}}\left(\sum_{n=1}^{\infty}\left(|t_n|^{q}\right)^{\frac{p}{q}}\right)^{\frac{q}{p}}\right)^{\frac{1}{q}}\|v\|_p}{\|v\|_p}$$

$$+ (p-2)(p\text{-}3) \lim_{v\to\theta} \frac{\left(\sum_{n=1}^{\infty}|\bar{t}_n+s_n|^{p}\right)^{\frac{1}{q}}\|v\|_p^2}{\|v\|_p}$$

$$\leq (p-1)(p\text{-}2) \lim_{v\to\theta} \frac{\left(\|\bar{x}\|_p+\|v\|_p\right)^{\frac{p-2q}{q}}\|v\|_p^2}{\|v\|_p} +2(p-1)(p\text{-}3) \lim_{v\to\theta} \frac{\left(\|\bar{x}\|_p+\|v\|_p\right)^{\frac{p-q}{q}}\|v\|_p^2}{\|v\|_p}$$

$$+ (p-2)(p\text{-}3) \lim_{v\to\theta} \frac{\left(\|\bar{x}\|_p+\|v\|_p\right)^{\frac{p}{q}}\|v\|_p^2}{\|v\|_p} \qquad (|s_n| \leq |t_n|)$$

$= 0.$

We estimate

$$\left\|\left\{\frac{(p-1)|\bar{t}_n|^{p-2}+(p-2)|\bar{t}_n|^{p-1}}{(1+|\bar{t}_n|)^2}t_n\right\}_{n=1}^{\infty}\right\|_q$$

$$= \left(\sum_{n=1}^{\infty}\left|\frac{(p-1)|\bar{t}_n|^{p-2}+(p-2)|\bar{t}_n|^{p-1}}{(1+|\bar{t}_n|)^2}t_n\right|^q\right)^{\frac{1}{q}}$$

$$\leq \left(\sum_{n=1}^{\infty}\left|\frac{(p-1)|\bar{t}_n|^{p-2}}{(1+|\bar{t}_n|)^2}t_n\right|^q\right)^{\frac{1}{q}} + \left(\sum_{n=1}^{\infty}\left|\frac{(p-2)|\bar{t}_n|^{p-1}}{(1+|\bar{t}_n|)^2}t_n\right|^q\right)^{\frac{1}{q}}$$

$$\leq (p-1)\left(\sum_{n=1}^{\infty}\left||\bar{t}_n|^{p-2}|t_n|\right|^q\right)^{\frac{1}{q}} + (p-2)\left(\sum_{n=1}^{\infty}\left||\bar{t}_n|^{p-1}|t_n|\right|^q\right)^{\frac{1}{q}}$$

$$\leq (p-1)\left(\sum_{n=1}^{\infty}|\bar{t}_n|^{(p-2)q}|t_n|^q\right)^{\frac{1}{q}} + (p-2)\left(\sum_{n=1}^{\infty}||\bar{t}_n|^{p-1}|^q\|v\|_p^q\right)^{\frac{1}{q}}$$

$$\leq (p-1)\left(\left(\sum_{n=1}^{\infty}|\bar{t}_n|^{(p-2)q\frac{p}{p-q}}\right)^{\frac{p-q}{p}}\left(\sum_{n=1}^{\infty}|t_n|^{q\frac{p}{q}}\right)^{\frac{q}{p}}\right)^{\frac{1}{q}} + (p-2)\|v\|_p\left(\sum_{n=1}^{\infty}|\bar{t}_n|^p\right)^{\frac{1}{q}}$$

$$= (p-1)\left(\left(\sum_{n=1}^{\infty}|\bar{t}_n|^p\right)^{\frac{p-q}{p}}\left(\sum_{n=1}^{\infty}|t_n|^p\right)^{\frac{q}{p}}\right)^{\frac{1}{q}} + (p-2)\|v\|_p\left(\sum_{n=1}^{\infty}|\bar{t}_n|^p\right)^{\frac{1}{q}}$$

$$= (p-1)\left(\|\bar{x}\|_p^{p-q}\|v\|_p^q\right)^{\frac{1}{q}} + (p-2)\|v\|_p\|\bar{x}\|_p^{\frac{p}{q}}$$

$$= (p-1)\|\bar{x}\|_p^{\frac{p-q}{q}}\ \|v\|_p + (p-2)\|\bar{x}\|_p^{\frac{p}{q}}\|v\|_p$$

$$= \left((p-1)\|\bar{x}\|_p^{\frac{p-q}{q}} + (p-2)\|\bar{x}\|_p^{\frac{p}{q}}\right)\|v\|_p. \tag{7.8}$$

Finally, by (7.8), for any $v = \{t_n\}_{n=1}^{\infty} \in l_p$, we have that

$$\|\nabla\mathfrak{J}(\bar{x})(v)\|_q = \left\|\left\{\frac{(p-1)|\bar{t}_n|^{p-2}+(p-2)|\bar{t}_n|^{p-1}}{(1+|\bar{t}_n|)^2}t_n\right\}_{n=1}^{\infty}\right\|_q$$

$$\leq \left((p-1)\|\bar{x}\|_p^{\frac{p-q}{q}} + (p-2)\|\bar{x}\|_p^{\frac{p}{q}}\right)\|v\|_p.$$

Notice that $\frac{p-q}{q} = p-2$ and $\frac{p}{q} = p-1$. It proves this theorem. □

By the connection between Gâteaux and Fréchet derivatives in Banach spaces, we have the following corollary of Theorem 6.1 immediately.

**Corollary 7.3**. *Let* $3 < p < \infty$. $\mathfrak{J}$ *is Gâteaux differentiable in* $l_p$ *and, for any* $\bar{x} = (\bar{t}_1, \bar{t}_2, \dots) \in l_p$, *the Gâteaux derivative* $\mathfrak{J}'(\bar{x})$ *satisfies that for any* $v = \{t_n\}_{n=1}^{\infty} \in l_p\backslash\{\theta\}$

$$\mathfrak{J}'(\bar{x})(v) = \nabla\mathfrak{J}(\bar{x})(v) = \left\{\frac{(p-1)|\bar{t}_n|^{p-2}+(p-2)|\bar{t}_n|^{p-1}}{(1+|\bar{t}_n|)^2}t_n\right\}_{n=1}^{\infty} \in l_q.$$

### 7.2. Mordukhovich Derivative and Covering Constant for the Generalized Duality Mapping $\mathfrak{J}$ in $l_p$

Similarly to Theorems 4.4 and 4.5, in this subsection, we apply Theorem 7.2 and Theorem 1.38 in [28] to find the explicit Mordukhovich derivative of the duality mapping $\mathfrak{J}: l_p \to l_q$ and calculate its covering constant.

**Theorem 7.4**. *Let* $3 < p < \infty$. $\mathfrak{J}: l_p \to l_q$ *is Mordukhovich differentiable in* $l_p$ *and for* $\bar{x} = \{\bar{t}_n\}_{n=1}^{\infty} \in l_p$, *and the Mordukhovich derivative* $\widehat{D}^*\mathfrak{J}(\bar{x}): l_p \to l_q$ *is equal to the Fréchet derivative* $\nabla\mathfrak{J}(\bar{x}): l_p \to l_q$ *with the following representation.*

$$\widehat{D}^*\mathfrak{J}(\bar{x}) = \nabla\mathfrak{J}(\bar{x}) = \begin{pmatrix} \frac{(p-1)|\bar{t}_1|^{p-2}+(p-2)|\bar{t}_1|^{p-1}}{(1+|\bar{t}_1|)^2} & 0 & \cdots \\ 0 & \frac{(p-1)|\bar{t}_2|^{p-2}+(p-2)|\bar{t}_2|^{p-1}}{(1+|\bar{t}_2|)^2} & 0 \\ \vdots & 0 & \ddots \end{pmatrix}.$$

*More precisely,* $\widehat{D}^*\mathfrak{J}(\bar{x}) = \nabla\mathfrak{J}(\bar{x})\colon l_p \to l_q$ *is a pointwise multiplication operator in* $l_p$ *and, for any* $v = \{t_n\}_{n=1}^{\infty} \in l_p$,

$$\widehat{D}^*\mathfrak{J}(\bar{x})(v) = \nabla\mathfrak{J}(\bar{x})(v) = \left\{\frac{(p-1)|\bar{t}_n|^{p-2}+(p-2)|\bar{t}_n|^{p-1}}{(1+|\bar{t}_n|)^2} t_n\right\}_{n=1}^{\infty} \in l_q.$$

*Proof*. By Theorem 4.8, the proof is similar to the proof of Theorem 6.3 and it is omitted here. □

**Theorem 7.5**. *Let* $3 < p < \infty$. *Let* $\bar{x} = \{\bar{t}_n\}_{n=1}^{\infty} \in l_p$ *and* $\bar{y} \in l_q$ *with* $\mathfrak{J}(\bar{x}) = \bar{y}$. *Then*

$$\hat{\alpha}(\mathfrak{J}, \bar{x}, \bar{y}) = 0.$$

*Proof*. By definition and by Theorem 7.4, similarly to the proof of Theorem 6.4, we calculate

$$\hat{\alpha}(\mathfrak{J}, \bar{x}, \bar{y}) = \sup_{\eta>0} \inf\left\{\|z\|_{l_q} : z = \widehat{D}^*\mathfrak{J}(x)(w), x \in \mathbb{B}_{l_p}(\bar{x},\eta), \mathfrak{J}(x) \in \mathbb{B}_{l_q}(\bar{y},\eta), \|w\|_{l_p} = 1\right\}$$

$$= \sup_{\eta>0} \inf\left\{\|z\|_{l_q} : z = \nabla\mathfrak{J}(x)(w), x \in \mathbb{B}_{l_p}(\bar{x},\eta), \mathfrak{J}(x) \in \mathbb{B}_{l_q}(\bar{y},\eta), \|w\|_{l_p} = 1\right\}$$

$$= \sup_{\eta>0} \inf\left\{\|\nabla\mathfrak{J}(x)(w)\|_{l_q} : x \in \mathbb{B}_{l_p}(\bar{x},\eta), \mathfrak{J}(x) \in \mathbb{B}_{l_q}(\bar{y},\eta), \|w\|_{l_p} = 1\right\}$$

$$\leq \sup_{\eta>0} \inf\left\{\|\nabla\mathfrak{J}(\bar{x})(w)\|_{l_q} : \|w\|_{l_p} = 1, w = \{s_n\}_{n=1}^{\infty} \in l_p\right\}$$

$$= \sup_{\eta>0} \inf\left\{\left(\sum_{n=1}^{\infty} \left|\frac{(p-1)|\bar{t}_n|^{p-2}+(p-2)|\bar{t}_n|^{p-1}}{(1+|\bar{t}_n|)^2} s_n\right|^q\right)^{\frac{1}{q}} : \|w\|_{l_p} = 1, w = \{s_n\}_{n=1}^{\infty} \in l_p\right\}$$

$$\leq \sup_{\eta>0} \inf\left\{\frac{(p-1)|\bar{t}_n|^{p-2}+(p-2)|\bar{t}_n|^{p-1}}{(1+|\bar{t}_n|)^2} : n = 1, 2, \ldots\right\}$$

$$= 0.$$ □

## 8. Conclusion and Remarks

In this paper, we investigate the Gâteaux, Fréchet and Mordukhovich differentiability of the normalized duality mapping and other two duality mappings in the uniformly convex and uniformly smooth Banach space $l_p$, for $1 < p < \infty$. For given point $\bar{x} \in l_p$, after the differentiability is proved, we find the explicit Gâteaux, Fréchet and Mordukhovich derivatives of the considered duality mappings at $\bar{x}$, respectively. By using their Mordukhovich derivatives, we find the covering constants for these duality mappings. In this paper, the Gâteaux differentiability of the normalized duality mapping $J$ in Banach space $l_p$ has been completed solved by Theorem 3.2, which provides the explicit Gâteaux derivative of $J$. However, regarding to the more interesting Fréchet differentiability, there are some unsolved problems in this paper, which are listed below for interesting reads consideration.

**Problem 1**. Let $2 < p < \infty$. Theorem 4.4 proves that, for given $\bar{x} = (\bar{t}_1, \bar{t}_2, \dots) \in l_p\backslash\{\theta\}$, if $|\mathcal{N}(\bar{x})| < \infty$, then $J$ is Fréchet differentiable at $\bar{x}$. In this case, for any $v = \{s_n\}_{n=1}^{\infty} \in l_p\backslash\{\theta\}$ with ratio sequence $\{\lambda_n\}_{n=1}^{\infty}$ with respect to $\bar{x}$, the Fréchet derivative $\nabla J(\bar{x})(v)$ satisfies that

$$\nabla J(\bar{x})(v) = J'(\bar{x})(v) = J(\bar{x})\left((p-1)\begin{pmatrix}\lambda_1 & 0 & \dots \\ 0 & \lambda_2 & 0 \\ \vdots & 0 & \ddots\end{pmatrix} - \frac{(p-2)\sum_{i=1}^{\infty}\lambda_i|\bar{t}_i|^p}{\sum_{i=1}^{\infty}|\bar{t}_i|^p}\begin{pmatrix}1 & 0 & \dots \\ 0 & 1 & 0 \\ \vdots & 0 & \ddots\end{pmatrix}\right).$$

Under the assumption that $|\mathcal{N}(\bar{x})|$, the above element in $l_q$ has finite number of nonzero terms.

However, in general, for any $\bar{x} = (\bar{t}_1, \bar{t}_2, \dots) \in l_p\backslash\{\theta\}$, by Lemma 4.2 and part (ii) in Theorem 3.2, the Gâteaux differential operator $J'(\bar{x})(\cdot): l_p \to l_q$ represented as above is a continuous and linear operator with Lipchitz property in $l_p$. This raises the following question:

Let $2 < p < \infty$ and let arbitrarily given $\bar{x} = (\bar{t}_1, \bar{t}_2, \dots) \in l_p\backslash\{\theta\}$. Is $J$ Fréchet differentiable at $\bar{x}$ satisfying, for any $v = \{s_n\}_{n=1}^{\infty} \in l_p\backslash\{\theta\}$ with ratio sequence $\{\lambda_n\}_{n=1}^{\infty}$ with respect to $\bar{x}$,

$$\nabla J(\bar{x})(v) = J(\bar{x})\left((p-1)\begin{pmatrix}\lambda_1 & 0 & \dots \\ 0 & \lambda_2 & 0 \\ \vdots & 0 & \ddots\end{pmatrix} - \frac{(p-2)\sum_{i=1}^{\infty}\lambda_i|\bar{t}_i|^p}{\sum_{i=1}^{\infty}|\bar{t}_i|^p}\begin{pmatrix}1 & 0 & \dots \\ 0 & 1 & 0 \\ \vdots & 0 & \ddots\end{pmatrix}\right)?$$

**Problem 2**. In sections 6 and 7, we study the Gâteaux, Fréchet and Mordukhovich differentiability of the duality mapping $\mathcal{J}$ and the generalized duality mapping $\mathfrak{J}$ in $l_p$, for $3 < p < \infty$. Their Gâteaux, Fréchet and Mordukhovich derivatives are explicitly solved, respectively. We believe that is interesting to investigate the Gâteaux, Fréchet and Mordukhovich differentiability of the duality mapping $\mathcal{J}$ and the generalized duality mapping $\mathfrak{J}$ in $l_p$, for $1 < p \leq 3$.